\documentclass[a4paper, 11pt]{article}

\usepackage{amsmath, amssymb, anysize, graphicx,
	latexsym, stmaryrd, xspace,amsfonts}
\usepackage{amsthm}
\usepackage[ngerman,english]{babel}    
\usepackage[right]{eurosym}
\usepackage{csquotes}
\usepackage{array}
\usepackage[normalem]{ulem}
\usepackage{cancel}
\usepackage{tikz}
\usetikzlibrary{positioning}
\usetikzlibrary{arrows}
\usetikzlibrary{shapes}

\usepackage{helvet}
\usepackage{color}
\usepackage{empheq}

\usepackage{fancyhdr, lastpage}
\usepackage{hyperref}
\usepackage{cases}
\usepackage{mathabx}
\usepackage{xfrac}

\usepackage{wrapfig}
\usepackage{figbib}
\usepackage{cite}

\usepackage{enumerate}
\usepackage[capitalise]{cleveref}
\usepackage{verbatim}

  \usepackage{pgfplots}
\pgfplotsset{compat=newest}
\usetikzlibrary{plotmarks}
\usetikzlibrary{arrows.meta}
\usetikzlibrary{matrix}
\usetikzlibrary{pgfplots.groupplots}
\pgfplotsset{plot coordinates/math parser=false}
\pgfplotsset{try min ticks=3}
\usepackage{grffile}
\pgfplotsset{plot coordinates/math parser=false}
\newlength\figureheight
\newlength\figurewidth
\pgfplotsset{
	colormap={ECM}{
   	    rgb255=(255.000, 255.000, 255.000)
		rgb255=(247.714, 247.714, 247.714)
		rgb255=(240.429, 240.429, 240.429)
		rgb255=(233.143, 233.143, 233.143)
		rgb255=(225.857, 225.857, 225.857)
		rgb255=(218.571, 218.571, 218.571)
		rgb255=(211.286, 211.286, 211.286)
		rgb255=(204.000, 204.000, 204.000)
		rgb255=(200.821, 200.821, 200.821)
		rgb255=(197.641, 197.641, 197.641)
		rgb255=(194.462, 194.462, 194.462)
		rgb255=(191.282, 191.282, 191.282)
		rgb255=(188.103, 188.103, 188.103)
		rgb255=(184.923, 184.923, 184.923)
		rgb255=(181.744, 181.744, 181.744)
		rgb255=(178.564, 178.564, 178.564)
		rgb255=(175.385, 175.385, 175.385)
		rgb255=(172.205, 172.205, 172.205)
		rgb255=(169.026, 169.026, 169.026)
		rgb255=(165.846, 165.846, 165.846)
		rgb255=(162.667, 162.667, 162.667)
		rgb255=(159.487, 159.487, 159.487)
		rgb255=(156.308, 156.308, 156.308)
		rgb255=(153.128, 153.128, 153.128)
		rgb255=(149.949, 149.949, 149.949)
		rgb255=(146.769, 146.769, 146.769)
		rgb255=(143.590, 143.590, 143.590)
		rgb255=(140.410, 140.410, 140.410)
		rgb255=(137.231, 137.231, 137.231)
		rgb255=(134.051, 134.051, 134.051)
		rgb255=(130.872, 130.872, 130.872)
		rgb255=(127.692, 127.692, 127.692)
		rgb255=(124.513, 124.513, 124.513)
		rgb255=(121.333, 121.333, 121.333)
		rgb255=(118.154, 118.154, 118.154)
		rgb255=(114.974, 114.974, 114.974)
		rgb255=(111.795, 111.795, 111.795)
		rgb255=(108.615, 108.615, 108.615)
		rgb255=(105.436, 105.436, 105.436)
		rgb255=(102.256, 102.256, 102.256)
		rgb255=(99.077, 99.077, 99.077)
		rgb255=(95.897, 95.897, 95.897)
		rgb255=(92.718, 92.718, 92.718)
		rgb255=(89.538, 89.538, 89.538)
		rgb255=(86.359, 86.359, 86.359)
		rgb255=(83.179, 83.179, 83.179)
		rgb255=(80.000, 80.000, 80.000)
		rgb255=(75.294, 75.294, 75.294)
		rgb255=(70.588, 70.588, 70.588)
		rgb255=(65.882, 65.882, 65.882)
		rgb255=(61.176, 61.176, 61.176)
		rgb255=(56.471, 56.471, 56.471)
		rgb255=(51.765, 51.765, 51.765)
		rgb255=(47.059, 47.059, 47.059)
		rgb255=(42.353, 42.353, 42.353)
		rgb255=(37.647, 37.647, 37.647)
		rgb255=(32.941, 32.941, 32.941)
		rgb255=(28.235, 28.235, 28.235)
		rgb255=(23.529, 23.529, 23.529)
		rgb255=(18.824, 18.824, 18.824)
		rgb255=(14.118, 14.118, 14.118)
		rgb255=(9.412, 9.412, 9.412)
		rgb255=(4.706, 4.706, 4.706)
		rgb255=(0.000, 0.000, 0.000)
	},
colormap={MCCs}{
rgb255=(255.000, 245.000, 235.000)
rgb255=(254.997, 243.094, 232.067)
rgb255=(254.937, 241.233, 228.926)
rgb255=(254.832, 239.395, 225.572)
rgb255=(254.693, 237.556, 222.000)
rgb255=(254.529, 235.694, 218.205)
rgb255=(254.350, 233.785, 214.180)
rgb255=(254.163, 231.807, 209.922)
rgb255=(253.977, 229.735, 205.425)
rgb255=(253.797, 227.548, 200.684)
rgb255=(253.628, 225.222, 195.696)
rgb255=(253.475, 222.735, 190.455)
rgb255=(253.340, 220.063, 184.959)
rgb255=(253.224, 217.183, 179.202)
rgb255=(253.128, 214.072, 173.181)
rgb255=(253.049, 210.706, 166.893)
rgb255=(252.985, 207.062, 160.335)
rgb255=(252.932, 203.140, 153.528)
rgb255=(252.885, 198.992, 146.537)
rgb255=(252.845, 194.671, 139.432)
rgb255=(252.817, 190.233, 132.282)
rgb255=(252.809, 185.732, 125.152)
rgb255=(252.836, 181.219, 118.107)
rgb255=(252.915, 176.749, 111.209)
rgb255=(253.067, 172.373, 104.519)
rgb255=(253.286, 168.109, 98.063)
rgb255=(253.525, 163.927, 91.816)
rgb255=(253.737, 159.795, 85.752)
rgb255=(253.870, 155.682, 79.844)
rgb255=(253.876, 151.557, 74.064)
rgb255=(253.703, 147.391, 68.387)
rgb255=(253.302, 143.154, 62.783)
rgb255=(252.623, 138.819, 57.225)
rgb255=(251.658, 134.383, 51.706)
rgb255=(250.439, 129.865, 46.249)
rgb255=(249.002, 125.286, 40.875)
rgb255=(247.380, 120.665, 35.605)
rgb255=(245.608, 116.020, 30.466)
rgb255=(243.720, 111.370, 25.489)
rgb255=(241.752, 106.733, 20.720)
rgb255=(239.732, 102.124, 16.232)
rgb255=(237.621, 97.561, 12.060)
rgb255=(235.332, 93.064, 8.300)
rgb255=(232.779, 88.662, 5.491)
rgb255=(229.875, 84.391, 3.534)
rgb255=(226.536, 80.296, 2.245)
rgb255=(222.678, 76.428, 1.465)
rgb255=(218.218, 72.845, 1.058)
rgb255=(213.092, 69.599, 0.909)
rgb255=(207.358, 66.695, 0.949)
rgb255=(201.138, 64.102, 1.130)
rgb255=(194.551, 61.777, 1.413)
rgb255=(187.716, 59.675, 1.763)
rgb255=(180.746, 57.746, 2.152)
rgb255=(173.755, 55.940, 2.555)
rgb255=(166.851, 54.213, 2.952)
rgb255=(160.142, 52.522, 3.325)
rgb255=(153.730, 50.833, 3.662)
rgb255=(147.716, 49.115, 3.951)
rgb255=(142.200, 47.337, 4.183)
rgb255=(137.276, 45.473, 4.338)
rgb255=(133.040, 43.491, 4.357)
rgb255=(129.585, 41.352, 4.241)
rgb255=(127.000, 39.000, 4.000)
},
colormap={PCCs}{
rgb255=(247.000, 251.000, 255.000)
rgb255=(243.029, 249.104, 253.568)
rgb255=(239.374, 247.151, 252.307)
rgb255=(235.992, 245.151, 251.190)
rgb255=(232.841, 243.113, 250.193)
rgb255=(229.877, 241.044, 249.291)
rgb255=(227.056, 238.953, 248.459)
rgb255=(224.334, 236.847, 247.671)
rgb255=(221.668, 234.736, 246.904)
rgb255=(219.013, 232.626, 246.134)
rgb255=(216.326, 230.526, 245.335)
rgb255=(213.561, 228.442, 244.484)
rgb255=(210.675, 226.383, 243.557)
rgb255=(207.623, 224.355, 242.530)
rgb255=(204.360, 222.366, 241.380)
rgb255=(200.838, 220.423, 240.083)
rgb255=(197.012, 218.532, 238.616)
rgb255=(192.854, 216.681, 236.979)
rgb255=(188.376, 214.824, 235.212)
rgb255=(183.596, 212.912, 233.361)
rgb255=(178.530, 210.896, 231.472)
rgb255=(173.194, 208.730, 229.589)
rgb255=(167.601, 206.369, 227.757)
rgb255=(161.766, 203.767, 226.019)
rgb255=(155.698, 200.883, 224.418)
rgb255=(149.421, 197.716, 222.959)
rgb255=(142.978, 194.318, 221.601)
rgb255=(136.416, 190.748, 220.301)
rgb255=(129.788, 187.060, 219.013)
rgb255=(123.151, 183.308, 217.695)
rgb255=(116.574, 179.547, 216.306)
rgb255=(110.137, 175.828, 214.805)
rgb255=(103.933, 172.200, 213.154)
rgb255=(98.014, 168.667, 211.347)
rgb255=(92.371, 165.195, 209.413)
rgb255=(86.993, 161.744, 207.380)
rgb255=(81.861, 158.278, 205.275)
rgb255=(76.953, 154.762, 203.127)
rgb255=(72.239, 151.159, 200.960)
rgb255=(67.680, 147.435, 198.802)
rgb255=(63.226, 143.558, 196.676)
rgb255=(58.850, 139.535, 194.580)
rgb255=(54.544, 135.396, 192.496)
rgb255=(50.301, 131.174, 190.404)
rgb255=(46.118, 126.899, 188.285)
rgb255=(41.996, 122.601, 186.122)
rgb255=(37.942, 118.310, 183.896)
rgb255=(33.975, 114.055, 181.591)
rgb255=(30.123, 109.861, 179.188)
rgb255=(26.402, 105.728, 176.655)
rgb255=(22.819, 101.642, 173.951)
rgb255=(19.390, 97.592, 171.039)
rgb255=(16.147, 93.564, 167.879)
rgb255=(13.138, 89.547, 164.434)
rgb255=(10.431, 85.531, 160.667)
rgb255=(8.225, 81.504, 156.542)
rgb255=(6.715, 77.459, 152.026)
rgb255=(5.818, 73.386, 147.086)
rgb255=(5.446, 69.278, 141.691)
rgb255=(5.505, 65.127, 135.812)
rgb255=(5.895, 60.929, 129.424)
rgb255=(6.513, 56.677, 122.504)
rgb255=(7.251, 52.368, 115.034)
rgb255=(8.000, 48.000, 107.000)
},
colormap={AC}{
rgb255=(0.000, 68.000, 27.000)
rgb255=(0.000, 74.369, 28.400)
rgb255=(0.000, 80.353, 30.025)
rgb255=(0.000, 85.971, 31.886)
rgb255=(0.000, 91.247, 33.980)
rgb255=(0.000, 96.208, 36.295)
rgb255=(0.000, 100.885, 38.813)
rgb255=(0.000, 105.309, 41.511)
rgb255=(0.486, 109.515, 44.364)
rgb255=(4.950, 113.540, 47.346)
rgb255=(10.399, 117.421, 50.431)
rgb255=(15.766, 121.197, 53.594)
rgb255=(20.520, 124.911, 56.813)
rgb255=(24.842, 128.602, 60.063)
rgb255=(28.810, 132.316, 63.326)
rgb255=(32.461, 136.095, 66.578)
rgb255=(35.809, 139.985, 69.802)
rgb255=(38.924, 144.005, 72.984)
rgb255=(41.988, 148.120, 76.124)
rgb255=(45.170, 152.288, 79.222)
rgb255=(48.615, 156.467, 82.278)
rgb255=(52.447, 160.612, 85.292)
rgb255=(56.764, 164.679, 88.264)
rgb255=(61.643, 168.621, 91.192)
rgb255=(67.130, 172.391, 94.077)
rgb255=(73.161, 175.969, 96.935)
rgb255=(79.558, 179.374, 99.811)
rgb255=(86.174, 182.626, 102.750)
rgb255=(92.892, 185.748, 105.799)
rgb255=(99.618, 188.764, 109.003)
rgb255=(106.277, 191.698, 112.410)
rgb255=(112.804, 194.575, 116.065)
rgb255=(119.144, 197.420, 120.013)
rgb255=(125.285, 200.243, 124.253)
rgb255=(131.251, 203.033, 128.738)
rgb255=(137.062, 205.778, 133.419)
rgb255=(142.733, 208.469, 138.249)
rgb255=(148.278, 211.095, 143.178)
rgb255=(153.707, 213.643, 148.158)
rgb255=(159.030, 216.103, 153.141)
rgb255=(164.254, 218.463, 158.080)
rgb255=(169.380, 220.724, 162.955)
rgb255=(174.405, 222.892, 167.767)
rgb255=(179.326, 224.972, 172.515)
rgb255=(184.140, 226.971, 177.201)
rgb255=(188.842, 228.897, 181.824)
rgb255=(193.431, 230.757, 186.385)
rgb255=(197.901, 232.558, 190.885)
rgb255=(202.250, 234.306, 195.322)
rgb255=(206.471, 236.002, 199.688)
rgb255=(210.555, 237.641, 203.968)
rgb255=(214.493, 239.218, 208.146)
rgb255=(218.277, 240.731, 212.210)
rgb255=(221.896, 242.176, 216.143)
rgb255=(225.342, 243.548, 219.931)
rgb255=(228.605, 244.844, 223.559)
rgb255=(231.677, 246.059, 227.011)
rgb255=(234.547, 247.190, 230.273)
rgb255=(237.206, 248.234, 233.330)
rgb255=(239.645, 249.186, 236.166)
rgb255=(241.853, 250.042, 238.766)
rgb255=(243.822, 250.799, 241.115)
rgb255=(245.540, 251.453, 243.198)
rgb255=(247.000, 252.000, 245.000)
}
}
\DeclarePairedDelimiter{\floor}{\lfloor}{\rfloor}
	
\usepackage[a4paper, left=2.4cm, right=2.0cm, top=2.5cm,bottom=2cm]{geometry}
\newtheorem{theo}{Theorem}[section]
\newtheorem{lem}[theo]{Lemma}

\newtheorem{prop}[theo]{Proposition}
\newtheorem{defi}[theo]{Definition}

\renewcommand{\theequation}{\arabic{section}.\arabic{equation}}
\newcommand{\mysection}[1]{\section{#1} \setcounter{equation}{0}}

\newcounter{expcounter}			% own experiment environement
\newenvironment{experiment}
{	\refstepcounter{expcounter}
	{~\\\textbf{Experiment \theexpcounter~---~}}
}

\renewcommand{\(}{\left(}
\renewcommand{\)}{\right)}

\newcommand{\pd}{\partial}

\newcommand{\R}{\mathbb{R}}
\newcommand{\N}{\mathbb{N}}
\newcommand{\io}{\int_\Omega}
\newcommand{\nn}{\nonumber}
\newcommand{\eps}{\varepsilon}
\newcommand{\pO}{\partial\Omega}
\newcommand{\hte}{\widehat{T_\eps}}

\newcommand{\Eeps}{{\cal E_\eps}}
\newcommand{\Deps}{{\cal D_\eps}}
\newcommand{\wto}{\rightharpoonup}
	\title{Modeling multiple taxis: tumor invasion with phenotypic heterogeneity, haptotaxis, and unilateral interspecies repellence}

\author{Niklas Kolbe$^{1}$, \quad Nikolaos Sfakianakis$^{2}$, \quad Christian Stinner$^{3}$, \\
Christina Surulescu$^{4}$, \quad Jonas Lenz$^{3,5}$\\
{\small $^{1}$ Kanazawa University, Faculty of Mathematics \& Physics,}\\  
{\small Kakuma, Kanazawa 920-1192, Japan}\\
{\small $^{2}$ University of St. Andrews, School of Mathematics \& Statistics,}\\
{\small North Haugh, St. Andrews, Fife, KY16 9SS, Scotland, UK}\\
{\small $^{3}$ Technische Universit\"at Darmstadt, Fachbereich Mathematik,}\\
{\small Schlossgartenstrasse 7, 64289 Darmstadt, Germany} \\
{\small $^{4}$ Technische Universit\"{a}t Kaiserslautern, Felix-Klein-Zentrum f\"{u}r Mathematik,} \\
{\small Paul-Ehrlich-Str. 31, 67663 Kaiserslautern, Germany}\\
{\small $^{5}$ Johannes Gutenberg-Universit\"{a}t Mainz, Institut f\"{u}r Mathematik,} \\
{\small Staudingerweg 9, 55128 Mainz, Germany}\\
{\small (kolbe@staff.kanazawa-u.ac.jp, n.sfakianakis@st-andrews.ac.uk,  stinner@mathematik.tu-darmstadt.de,}\\
{\small surulescu@mathematik.uni-kl.de, j.lenz@uni-mainz.de)}
}
\begin{document}
\maketitle

%%%%%%%%%%%%%%%%%%%%%%%%%%%%%%%%%%%%%%%%%%%%%%%%%%%%%%%%%%%%
%%%%%%%%%%%%%%%%%%%%%%%%%%%%%%%%%%%%%%%%%%%%%%%%%%%%%%%%%%%%
%%%%%%%%%%%%%%%%%%%%%%%%%%%%%%%%%%%%%%%%%%%%%%%%%%%%%%%%%%%%

\begin{abstract}
We provide a short review of existing models with multiple taxis performed by (at least) one species and consider a new mathematical model for tumor invasion featuring two mutually exclusive cell phenotypes (migrating and proliferating). The migrating cells perform nonlinear diffusion and two types of taxis in response to non-diffusing cues: away from proliferating cells and up the gradient of surrounding tissue. Transitions between the two cell subpopulations are influenced by subcellular (receptor binding) dynamics, thus conferring the setting a multiscale character.  

We prove global existence of weak solutions to a simplified model version and perform numerical simulations for the full setting under several phenotype switching and motility scenarios. We also compare (via simulations) this model with the corresponding haptotaxis-chemotaxis one featuring indirect chemorepellent production and provide a discussion about possible model extensions and mathematical challenges. %{\cb To be completed after simulations}
\end{abstract}

%%%%%%%%%%%%%%%%%%%%%%%%%%%%%%%%%%%%%%%%%%%%%%%%%%%%%%%%
%%%%%%%%%%%%%%%%%%%%%%%%%%%%%%%%%%%%%%%%%%%%%%%%%%%%%%%%
%%%%%%%%%%%%%%%%%%%%%%%%%%%%%%%%%%%%%%%%%%%%%%%%%%%%%%%%

\mysection{Introduction}\label{sec:intro}

%\subsection{Interactions driving motility and tactic behavior}

The growth, development, spread, and even survival of biological organisms and  populations are tightly connected to their environment and its characteristics. A multitude of different biological and chemical species interact with each other and with their surroundings in order to perform the mentioned processes. Intra- and interpopulation reciprocities together with environmental heterogeneities have a decisive influence on the fate of individuals as well as populations. A prominent example is given by the motility and proliferation behavior of cells in response to biochemical and biophysical cues present in their extracellular space. Not only can these induce a pronounced migratory bias of the cells towards favorable signals or away from hostile ones (also known as taxis), but they can even determine their phenotype and interaction strategies. Moreover, cells are able to modify their environment according to their needs, a behavior which is obvious in higher species. 

Since such complex processes and their mutual conditioning scenarios are difficult to assess experimentally, mathematical models can help understanding the underlying biological mechanisms, test hypotheses, make predictions, and even suggest new experiments and conjectures. During the last five decades the scientific community has witnessed an unprecedented development of mathematical descriptions of problems motivated by life sciences, and to a large extent by cell biology, triggered among others by Keller and Segel's work on chemotaxis modeling \cite{Keller1971}.\\[-2ex]

%Cell migration is essential for many processes involved in tumor invasion, wound healing, development and functioning of organs and tissues, biofilm formation, etc. A prominent example is the ability of cancer cells to degrade the surrounding tissue while migrating towards blood vessels or away from unfavourable (e.g., hypoxic) regions. Mathematical models can help understanding the underlying biological mechanisms and suggest or test new therapies.

%\subsection{A short review of models with multiple taxis}

Several fine reviews of chemotaxis models are available, see e.g. \cite{BeBeTaWi,HiPa08,Horstmann,Painter2019}, addressing aspects of well-posedness, long time behavior, and patterns.  Without the pretension of exhaustiveness we provide a short review of existing models from the perspective of multiple taxis. By the latter we mean systems of PDEs and ODEs characterizing the evolution of one or several species perceiving various cues and responding to them. Of these species at least one is biasing its motion according to at least two tactic signals. \\[-2ex]

A general framework for \textit{macroscopic multiple taxis} can be formulated more precisely as 

\begin{align}
\partial_t u_i=\nabla\cdot\left(a_{i0}(u)\nabla u_i\right)-\nabla\cdot\left(\sum_{j=1}^{m-1} a_{ij}(u)\nabla b_{ij}(u)\right)+a_{im}(u),\label{EqLoc}
\end{align}
each such PDE describing the evolution of a subpopulation density $u_i$ as component of a tuple $u=(u_1,\dots,u_N)$ of $N\in\N$ variables representing cell densities, volume fractions of various non-diffusing constituents of the surrounding environment (e.g., tissue fibers), concentrations of nutrients and chemical signals, etc. Thereby, $m\in\N$ and the coefficients have the following meaning: $a_{i0}(u)\ge 0$ is the diffusion coefficient, $a_{ij}(u)$ and $b_{ij}(u)$ for $j\in\{1,\dots,m-1\}$ describe  tactic sensitivities and signal functions, respectively, and, eventually, $a_{im}(u)$ is the reaction-interaction term. For multiple taxis there should be $N\ge 3$ and at least one $i$ s.t. the flux in the middle of the right hand side of \eqref{EqLoc} has at least two nonzero terms with $j\neq i$. %for at least two $j$'s ($j\neq i$) it holds that $a_{ij}(u)\not\equiv 0$
The review \cite{CPSZ19} addresses the possibility of extending this   framework upon including nonlocalities in several different ways.\\[-2ex]

In the following we will classify the available models according to the type of taxis mixtures they involve. In view of this, we identify the following categories:

\begin{enumerate}[(i)]
	\item \label{M.i} multiple chemotaxis, meaning that the tactic population(s) follow and/or are repelled by chemical signals, all of which diffuse;
		\item \label{M.ii} mixtures of chemo- and haptotaxis, where some of the tactic signals (all of which are non-population components of the environment) are diffusing, while the others are not; 
		\item \label{M.iii} interpopulation movement responses, possibly following further diffusing or non-diffusing tactic signals.
\end{enumerate}

In all these models the tactic cues are produced and/or degraded by the described species. For the sake of conciseness we only address here settings which feature fully parabolic reaction-diffusion-taxis PDEs (possibly coupled with ODEs) and exclude parabolic-elliptic systems. We would like, however, to emphasize that the study of the latter is not of less interest, as it can provide insight into the ways of handling the (often increased) difficulties related to fully parabolic systems. Moreover, it can often deliver results about qualitative behavior of solutions which are out of reach in the fully parabolic case. \\[-2ex]

Most \textit{models with multiple chemotaxis belonging to category} (\ref{M.i}) have been introduced in the context of cell populations migrating towards and/or away from chemical cues, see \cite{Luca2003,Orme1996,Painter2000}. 
A model for microglia chemotaxis towards three chemoattractants and away from one chemorepellent was proposed in \cite{Silchenko2015}; in \cite{JLi2018b,JLi2018} the authors considered a model with attractive-repulsive chemotaxis with and respectively without (logistic) growth of the tactic population and proved global existence and uniqueness, along with boundedness and long time behavior of solutions. The setting in \cite{Kntsdttir2014} features a kind of reciprocal chemotaxis of macrophages interacting with tumor cells; thereby, the former follow a chemical produced by the latter, while these in turn follow two signals, one of which is expressed by themselves (autocrine), while the other is produced by the macrophages (paracrine). \\[-2ex]

Further multiple chemotaxis models have been considered in the framework of indirect intra- or interspieces interactions. For instance, the one in \cite{Logan1998} describes the behavior of mountain pine beetles by dividing the population into one subpopulation which is motile and another, which is sessile, producing -alone or in interaction with further actors in the environment- the diffusing signals attracting the former. 
All mentioned works in category (\ref{M.i}) contain numerical studies and (linear) stability analyses to assess availability and types of patterning which the considered models are able to describe.\\[-2ex]

%\cite{Potts2019}: interpopulation movement responses between $N\ge 2$ populations, pattern analysis, no study of existence/uniqueness, asymptotic behavior, boundedness; general framework, but results only in 1D, mutual avoidance/attraction or taxis partially influenced by other population(s), partially by a non-diffusing signal, nonlocalities in taxis are allowed (type 2 in review nonlocal)

%no analysis 

%chemo-chemo in a system of two different populations, each performing a single type of chemotaxis
%-taxis cascade, foragers-scroungers where the exploiter subpopulation follows the gradient of forager density:  \cite{Winkler2019} (nD), global existence, long-time behavior, \cite{TaoWi2019} (1D), same, but more general source terms for attractant.
%-populations competing for the same diffusing signal which they also produce: \cite{BaiWi16}, \cite{Zhang2013} (nD, global existence \& uniqueness, linear diffusion; \cite{BaiWi16} asymptotic behavior \& global boundedness, \cite{Zhang2013} global boundedness)

Concerning \textit{models in category} (\ref{M.ii}), the 1D setting introduced by Chaplain and Lolas in \cite{ChapLol05} is presumably the first to account for mixed multiple taxis in the motile behavior of one population responding to diffusive and substrate-bound signals present in its environment. Concretely, it describes the evolution of tumor cell density using the extracellular matrix (ECM) as a haptoattractant and moving up gradients of two diffusing chemicals (urokinase plasminogen activator uPA and plasminogen activator inhibitor PAI) which the cancer cells produce in order to degrade the surrounding tissue, thus enabling invasion. The work, focusing on illustrations -via numerical simulations- of the combined effects of the addressed taxis and source terms, was followed by a chemotaxis-haptotaxis  multidimensional version \cite{ChLo06}, which generated a whole avalanche of papers dealing with the mathematical analysis of several versions therewith, in terms of well-posedness, boundedness, and asymptotic behavior of solutions.\\[-2ex] 

Haptotaxis is an essential component of cell motility in  tissues  \cite{carter}, being also relevant for pattern formation  and phenotypic heterogeneity of cells moving through heterogeneous  environments \cite{ChLo06,Maini1989,mallet04}. 
Chemotaxis-haptotaxis models exhibit supplementary difficulties, as the latter kind of advection characterizes motility bias towards the gradient of an immovable signal, whose evolution is described by an ODE. This corresponds to an everywhere degenerate reaction-diffusion PDE and has no regularizing effect.\\[-2ex]

Most chemotaxis-haptotaxis models investigated from the above mentioned viewpoint of mathematical analysis are versions of the Chaplain-Lolas models \cite{ChapLol05,ChLo06}. Of these, the vast majority involves linear diffusion, with constant \cite{Cao2016,ChenTao2018,Jin2018,Ke_2018,Pang2017,PaWa18,Tao2009,TaoWang2008,Tao2015,TaoWi2019b,WangKe2016,Xiang_2019} or solution-dependent \cite{Pang2019,TaoCui10} taxis sensitivity functions. Cells migrate, however, with finite speed, thus models with nonlinear cell diffusion were considered, whereby the diffusion coefficient may \cite{li-lankeit,tao-winkler,zms} or may not \cite{Hu2017,Jia2020,LiuWang2017,Liu2016,Wang2016,Wang2016b} infer some (mild) form of degeneracy. Further (stronger) types of degeneracy were investigated for pure haptotaxis models in \cite{Winkler2018,winkler-surulescu,ZSH,ZSU}.\\[-2ex]

The model in \cite{MORALESRODRIGO2013} has a mathematical structure closely related to the Chaplain-Lolas models, but describes a different biological problem: the migration and growth of endothelial cells following concentrations of diffusing tumor angiogenic factors (TAF) and immotile  fibronectin. In \cite{HPW12} it was shown that a simplified 1D version of the chemotaxis-haptotaxis model by Chaplain and Lolas convergences for sufficiently small tissue density to a chemotaxis model with logistic source term. \\[-2ex]

Results about global boundedness are available under various scenarios for some of the above models \cite{Cao2016,ChenTao2018,Hu2017,Jia2020,Jin2018,Ke_2018,li-lankeit,LiuWang2017,Liu2016,PaWa18,Pang2019,TaoWi2019b,Wang2016,Wang2016b,WangKe2016,Xiang_2019,zms}. %Blow-up of solutions has also been proved under certain conditions \cite{TaoWang2008}. 
Whether or not the above type of models feature remodeling of the haptoattractant seems to essentially influence the long-time behavior of the solution, in particular that of tactic cells, as reviewed in \cite{BeBeTaWi}; see also \cite{Pang2019} for a result on asymptotics of a chemotaxis-haptotaxis system for tumor angiogenesis. The model recently introduced in \cite{SuWi19} and extending several versions of the Chaplain-Lolas settings \cite{ChapLol05,ChLo06} describes chemo-and haptotaxis of a cancer cell population with indirect production of both signals: they are expressed by a non-diffusing population (of cancer-associated fibroblasts) which is in turn activated ('produced') by the tactic one. This indirect signal production prevents solution blow-up -at least in 2D- even when source terms are absent in the species performing taxis.\\[-2ex]

The multiscale model with chemotaxis and haptotaxis introduced in \cite{MERAL15} and further studied in \cite{stinner_surulescu_winkler_14} belongs to this category (\ref{M.ii}) as well. It connects the macroscopic evolution of cells with the microscopic dynamics of (some of) their subcellular events, thus leading to a system of strongly coupled PDEs and ODEs with time delay, acting on two different scales and featuring nonlinear diffusion and taxis sensitivities. Global existence of weak solutions and some global (non-uniform) boundedness have been proved, the issue of uniqueness remaining open. Further multiscale models for cancer invasion facing similar problems, but involving only one type of taxis, have been proposed and studied in \cite{Meral2015,stinner-surulescu-meral,stinner-surulescu-uatay}.\\[-2ex]
%\cite{Hu2017}: nonlinear, nondegenerate diffusion, 3D, no haptoattractant reconstruction, gobal existence and uniqueness, boundedness

Yet another class of multiscale models with chemo- and haptotaxis was proposed and studied in \cite{KeSu2011,KELKEL2012} and readdressed in terms of numerical simulations in \cite{Knopoff2019} and of well-posedness under more or less restrictive assumptions in \cite{Lorenz2014,Nieto2016}. The models therein connect macroscopic descriptions of chemoattractant dynamics (reaction-diffusion equations) with mesoscopic kinetic transport equations for the evolution of cell and tissue density functions and can be seen to belong to the kinetic theory of active particles (KTAP) developed by Bellomo et al. \cite{Bell-KTAP}. For related multiscale models with multiple chemotaxis -thus falling into category (\ref{M.i})- we refer e.g., to \cite{chalub}. The way these settings account for chemo- and haptotaxis is different from their macroscopic counterparts: it can be through the choice of the turning kernel in the operator characterizing cell reorientations and/or through the transport terms w.r.t. the kinetic variables and the cell turning rates. However, under certain conditions it is possible to obtain from such kinetic descriptions in some appropriate, more or less rigorous macroscopic limit the 'usual' reaction-diffusion-advection systems with taxis - possibly with some different type(s) of diffusion. For a recent work in this direction starting from such multiscale meso-macro formulations and arriving in the (formal) limit at chemotaxis-haptotaxis systems we refer to \cite{CEKNSSW}.\\[-2ex]

A rather general framework in \textit{category} (\ref{M.iii}) has been proposed in \cite{Potts2019}: it characterizes interpopulation movement responses between $N\ge 2$ populations, with mutual avoidance/attraction or taxis partially influenced by other diffusing population(s), allowing for nonlocalities in the taxis operators. One model class therein involves multiple taxis, of which one kind is towards a non-diffusing signal (again via some spatial averaging), thus can be seen as an example belonging to category (\ref{M.ii}).

The work \cite{Vig2014} describes a population of macrophages following two kinds of bacteria, of which one is diffusing and the other non-motile. The model provides a rare example of a system with double taxis, but without diffusion of the tactic species. Well-posedness studies of such models seem to be absent; the same applies to the rigorous assessment of asymptotic behavior of their solutions. \\[-2ex]

In this note we propose a two-scale model for tumor growth and spread in a heterogeneous environment. The tumor itself is supposed to be phenotypically heterogeneous, with a subpopulaton of moving cells, while the rest are proliferating. The former perform haptotaxis towards increasing gradients of normal tissue and another taxis, away from aggregates of immotile, proliferating cells. As such, the model falls into this same category (\ref{M.iii}). More details will be provided in Section~\ref{sec:model}. \\[-2ex]

To complete this review we should also mention a further category of models which do contain transport terms biasing the population motility, however without necessarily directly involving gradients of tactic signals. These are models with nonlocal interactions, e.g. describing cell-cell and/or cell-tissue adhesions or nonlocal chemotaxis: the motile population(s) perceive(s) signals not only locally, but in a whole region (usually a ball) around that point. This is described by way of integral operators and applies not only to space, but possibly to other variables as well, although spatial nonlocalities are more frequent in the context of tactic motion. While formal passages from nonlocal models of this type to classical reaction-diffusion-taxis (so-called 'localizations') have been known for a while (see e.g. \cite{gerisch2008,Hillen2007}), a rigorous result was established only recently in \cite{KPSZ}, which also contains global existence results for a class of such nonlocal models. Each of the considered systems features single taxis, but the settings can be generalized to allow for multiple taxis - at least as far as modeling is concerned. Multiscale nonlocal models with single \cite{DTGC17} or multiple \cite{ESS-17} taxis have been considered as well, the latter also proving  global existence of weak solutions - the issues of uniqueness and boundedness remaining open. Another model which can be assigned to this category has been proposed in \cite{Chauviere2007} in the KTAP framework, however the chemotaxis described there is towards a given signal without its own dynamics, thus actually leading to single 'genuine' taxis.\\[-2ex]

%this model: haptotaxis, repellent taxis away from nondiffusing signal, multiscale
% future: cross-diffusion, p can have a term of the form $\Delta ((a_1+a_2p+a_3m)p)$.
The rest of this paper is structured as follows: in Section~\ref{sec:model} we introduce our multiscale model with double taxis and prove in Section~\ref{sec3} the global existence of a weak solution to a pure macroscopic version of it. Section~\ref{sec:numerics} provides a description of several different scenarios of the full model. In Section~\ref{sec:discussion} we provide a discussion about the model, its results, and its possible extensions. Finally, the Appendix contains information about the construction of the randomly structured ECM employed in the simulations and about the numerical method. The supplementary material therein provides information about video files of some selected simulation experiments.

\mysection{A novel model of tumor invasion with two types of taxis}\label{sec:model}

Cancer cell migration is a well-known characteristic of tumor development and an essential step in building metastases. It involves a large variety of processes and interactions, of which chemotaxis, haptotaxis, avoidance of hypoxic regions, degradation of normal cells and tissue remodeling, phenotypic heterogeneity and switch are but the most frequently described. They mainly act on the level of single cells, but are influenced by subcellular events and regulate intercellular behavior, thus contributing to the evolution of the whole tumor. \\[-2ex]

Among the many kinds of  phenotypic heterogeneity, the go-or-grow dichotomy \cite{giese-etal96,Hoek2008,Zheng2009} is perhaps closest related to migration and proliferation: it states that moving cells can become non-motile and proliferate, according to the specific signals they perceive in their environment. This behavior can deter tumor spread, enhancing instead mitosis and the associated faster cell cycling, thus potentially having great impact on the sensitivity towards therapy \cite{Giese2003,Zheng2009}. These facts have been exploited in several different models for glioma invasion \cite{ConteSurulescu,EKS16,Hatzikirou2010,Hunt2016,ZSH} and were confirmed therein (at least qualitatively) by numerical simulations.\\[-2ex]

In this work we propose a model for cancer invasion into the normal tissue (ECM and normal cells) which is based on the go-or-grow hypothesis and features two types of taxis: haptotaxis towards increasing tissue gradients and a drift accounting for interpopulation repellence, meaning that motile cells have a tendency to avoid proliferating (thus non-migrating) ones. Concretely, 
%
% Herausforderungen Heterogenität (s. Intro \cite{stinner-surulescu-uatay})\\
we consider the following (nondimensionalized) equations to describe the dynamics of migrating ($m$) and proliferating ($p$) cell densities interacting mutually and with the extracellular matrix (of density $v$):
\begin{subequations}\label{eq:model}
\begin{align}
		\pd_t m =& \nabla\cdot\( D_c\frac{1+ m\,p + m\,v + p\,v}{1+m\,(p+v)}\nabla m \) 
			-\nabla \cdot\( \chi_1(m,p,v) m\nabla v \) 
			+\nabla \cdot\( \chi_2(m,p,v) m\nabla p \)\notag\\
			&+\lambda(y,\zeta)\,p -\gamma(y,\zeta)\,m\label{eq:model-a}\\
		\pd_t p =& \mu \, p\,(1-(m+p)-v) -\lambda (y,\zeta)\,p +\gamma(y,\zeta)\,m\label{eq:model-b}\\
		\pd_t v =& -\delta\,(m+p)\,v+\mu_v\,v\,(1-(m+p)-v).\label{eq:model-c}
\end{align}
\end{subequations}

%\begin{subequations}\label{eq:model-neu}
%	\begin{align*}
%	\pd_t m =& \nabla\cdot\( D_c\frac{1+ m\,p + m\,v + p\,v}{1+m\,(p+v)}\nabla m \) 
%	-\nabla \cdot\( \chi_1(m,p,v) m\nabla v \) 
%	+\nabla \cdot\( \chi_2(m,p,v) m\nabla h \)\notag\\
%	&+\lambda(y,\zeta)\,p -\gamma(y,\zeta)\,m\notag \\
%	\pd_t p =& \mu \, p\,(1-(m+p)-v) -\lambda (y,\zeta)\,p +\gamma(y,\zeta)\,m\notag\\
%	\pd_t v =& -\delta\,h\,v+\mu_v\,v\,(1-(m+p)-v),\notag\\
%	\pd_th=&D_h\Delta h+\alpha (m+p)-\beta h
%	\end{align*}
%\end{subequations}

Thereby, the cell motility components in \eqref{eq:model-a} include nonlinear diffusion and taxis. The transitions $\lambda $ and $\gamma $ between the two cancer cell subpopulations are influenced by the amount of 
cell surface receptors bound to insoluble ligands in the tumor microenvironment. These bindings account for cell-cell (variable $\zeta $) and for cell-tissue interactions (variable $y$). The dynamics of bound receptors is given by the 'mass action kinetics' (the total amount of receptors is scaled to $1$ and is supposed to be conserved):
\begin{equation}\label{eq:microReact}
	\begin{aligned}
		1-y-\zeta + v &\underset{k_-}{\stackrel{k_+}{\rightleftharpoons}} y\\
		1-y-\zeta + m+p &\underset{k_-}{\stackrel{k_+}{\rightleftharpoons}} \zeta.
	\end{aligned}
\end{equation}
leading to the ODE system
\begin{equation}\label{eq:microODEs}
	\left.\begin{aligned}
		\pd_t y&=k_+\(1-y-\zeta\)v-k_-y\\
		\pd_t \zeta&=k_+\(1-y-\zeta\)\(m+p\)-k_-\zeta		
	\end{aligned}\right\},
\end{equation}
where we kept the same notations for the rates $k_+$ and $k_-$. In order not to further inflate the set of parameters we assume that the binding and detachment cell-cell and cell-tissue kinetics have comparable rates $k_+$ and $k_-$. \\[-2ex] 

The migrating cells 
perform haptotaxis (third term in \eqref{eq:model-a}) and a kind of repellent taxis in the direction of decreasing $p$-gradient. The latter aims at describing the fact that tumor cell motility is triggered (among others) by the cancer cell population growing and rendering the local environment hostile (e.g., by acidification due to excessive glycolysis, \cite{CCTAIT00,GPMCMCR14,UWE09}). Unlike in \cite{ESS-17,SMAMCI14,HiSu2,stinner-surulescu-meral}, here we do not account explicitly for the dynamics of an acidity variable. Instead, we assume that most of acidity is produced by the proliferating cells and accordingly orient the repellent taxis towards their negative gradient. Notice that this is not chemotaxis, as the variable $p$ infers neither diffusion nor transport. The diffusion of $m$ is nonlinear, being influenced by all model variables. The diffusion coefficient is assumed to be proportional to all mutual interactions, its effect being, however, limited by the moving cells interacting with the proliferating cells and with the surrounding tissue, i.e. we account for a 'saturation' in this interplay. Here we focus on investigating the effects of the two types of taxis, therefore we also require the diffusion coefficient to stay positive, as degenerate diffusion in connection with haptotaxis can be quite challenging by itself, see e.g. \cite{li-lankeit,winkler-surulescu,zms,ZSH,ZSU}. Both types of tumor cells degrade the tissue and limit both tissue regeneration and cancer proliferation in a logistic-like way.\\[-2ex]

The taxis sensitivity coefficients $\chi _1$ and $\chi_2$ are obtained as follows: The receptor binding dynamics (subcellular scale) is much faster than cell migration and proliferation on the population scale, therefore we assume that it equilibrates rapidly, so that we may consider just the steady states of the variables $y$ and $\zeta$. With the notation $k_D:=\frac{k_-}{k_+}$, these equilibrium states 
are obtained from \eqref{eq:microODEs}:
\begin{equation}\label{eq:microSteady}
	\left.\begin{aligned}
		\bar y&=\frac{1}{k_D+m+p+v}v\\
		\bar \zeta&=\frac{1}{k_D+m+p+v}(m+p)
	\end{aligned}\right\}.
\end{equation}

Accordingly we may consider the taxis sensitivity functions to depend on these equilibria reminding of Michaelis-Menten type kinetics, and concretely take 
\begin{equation}\label{sensitivities-first}
\chi_1(m,p,v)=\bar y,\qquad \chi_2(m,p,v)=\bar \zeta.
\end{equation} 
%The last choice omits $m$ in the numerator and accounts for the fact that we are interested in the sensitivity of migrating cells towards $p$-cells 

Alternatively we could choose 
\begin{equation}
	\left.\begin{aligned}\label{eq:microSteady2}
		\chi_1(m,p,v)&=\xi_1 \pd_v \bar y=\xi_1 \frac{k_D+m+p}{(k_D+m+p+v)^2}\\
		\chi_2(m,p,v)&=\xi_2 \pd_p \bar \zeta=\xi_2 \frac{k_D+v}{(k_D+m+p+v)^2}
	\end{aligned}\right\},
\end{equation}
in order to emphasize that the tactic sensitivity of migrating cells to tissue is basically due to the $v$-variation of receptor binding state to ECM, whereas for the sensitivity to proliferating cells only variations of $\bar \zeta $ with respect to $p$ are relevant.\\[-2ex] 
%{\cb (Bitte auch $\chi_1(m,p,v)=\bar y$ und $\chi_2(m,p,v)=\frac{p}{k_D+m+p+v}$ versuchen)}\\

Observe that considering the equilibrium states instead of the full dynamics of the receptor binding kinetics allows us to reduce the original two-scale problem to a pure macroscopic one. Thereby, the transition rates $\lambda $ and $\gamma $ are assumed as well to depend only on $\bar y$, $\bar \zeta $, so they can be written as nonlinear (but regular and bounded) functions of the macroscopic variables $m$, $p$, $v$. The model can be seen as an extension of those with only haptotaxis considered in \cite{stinner-surulescu-uatay,ZSH}. Here, however, we do not account for decay of either cell population due to therapy (as in \cite{stinner-surulescu-uatay}) and in contrast to \cite{ZSH} explicitly avoid degeneracy of diffusion for the migrating cells (by adding 1 to the numerator in the specific form of $D(m,p,v)$). Without this simplification the degeneracy occurs whenever at least two of the solution components vanish - unlike in \cite{ZSH}, which featured strong double degeneracy (for which vanishing of one involved solution component was enough).\\[-2ex]

Including in the model the two subpopulations with their phenotypic switch leads to mathematical challenges in connection with the twofold effects of moving cells: besides directly degrading the tissue they act on the one hand as a source for the proliferating cells and on the other hand contribute directly to their own decay and indirectly (via $p$) to tissue depletion. Neither the $p$-cells nor the tissue are motile, in particular no regularization by diffusion terms can be expected: this problem typically encountered when dealing with haptotaxis now extends to the interspecies repellence. 
For the proof of global existence of solutions to \eqref{eq:model} we will restrain to a version of the model which features simpler forms than \eqref{sensitivities-first} and constant transition rates $\lambda, \gamma$. The numerical simulations will be performed, however, also for the more general function choices mentioned above.

\mysection{Global existence of weak solutions}\label{sec3}
In this section we show the global existence of (at least) a weak solution to the following simplified version of \eqref{eq:model}, where
we choose $\chi_1 (m, p,v) := \frac{c_1 v}{1+v}$ and assume that $\chi_2 \equiv c_2$, $\lambda$, and $\gamma$ are constant.
In particular, the restrictions concerning $\chi_2$ are quite severe. Although some choices of the form $\chi_2 = \chi_2 (p)$
are possible as well (see the remarks after Theorem~\ref{theo3.2} and before Lemma~\ref{lem3.10}), these seem to be less appropriate from the  viewpoint of the biological problem addressed by this model than 
$\chi_2$ being constant. Hence, we consider
\begin{equation}\label{eq3.0}
	\left\{ \begin{array}{l}
	\partial _t m=\nabla \cdot \left( D(m, p, v)\nabla m \right) - \nabla \cdot \left( \frac{c_1 v}{1+v} m \nabla v 
	\right) + \nabla \cdot \left( c_2 m \nabla p \right)+ \lambda p - \gamma m,  \\[1mm]
	\partial_t p = \mu p \left( 1- (m+p) - \eta_1 v \right) - \lambda p + \gamma m, 
	  \\[1mm]
	\partial _t v= - \alpha mv - \beta pv + \mu_v v (1-v),   
	\end{array}\right.
\end{equation} 
for $(x,t) \in \Omega \times (0,\infty)$ endowed with no-flux boundary conditions
\begin{equation}\label{eq3.0b}
	D(m,p,v)\, \partial _{\nu }m 
	- \frac{c_1 v}{1+v} m\ \partial _{\nu }v +  c_2 m \ \partial_\nu p= 0,
	\qquad x\in\partial\Omega, \ t>0,
\end{equation}
and initial conditions
\begin{equation}\label{eq3.0i}
	m(x,0)=m_0(x), \quad p(x,0) = p_0 (x) \quad
	v(x,0)=v_0(x),
	\qquad x\in\Omega,
\end{equation}
where $\Omega \subset \R^n$ is a bounded domain with smooth boundary and $\nu$ denotes the outward unit normal on
$\partial \Omega$. We assume that
\begin{equation}\label{eq3.1}
  m_0\in C^0(\bar\Omega), \quad
	p_0, v_0\in W^{1,2}(\Omega) \cap C^0(\bar\Omega), \quad\mbox{with}\quad 
  m_0 \ge 0, \quad p_0 \ge 0, \quad 
	v_0 \ge 0 	\quad \mbox{in } \bar\Omega,
\end{equation}
that for any $A>0$ and $L>0$ there exist positive constants $C_1$ and $C_2$ such that 
\begin{equation}\label{eq3.2}
  \begin{array}{l}
	D \in C^3 ([0,\infty)^3) \cap W^{2,\infty} ([0,\infty) \times [0,A] \times [0,L]), \\[1mm]
  	0< C_2 \le D(m,p,v) \le C_1 \quad\mbox{for all } (m,p,v) \in [0,\infty) \times [0,A] \times [0,L], 
	\end{array}
\end{equation}
and that the parameters $c_1$, $c_2$, $\lambda$, $\gamma$, $\mu$, $\eta_1$, $\mu_v$,
$\alpha := \delta +\mu_v \eta_2$, and $\beta := \delta +\mu_v \eta_2$ are positive. \\[-2ex]

Similarly to \cite{stinner_surulescu_winkler_14,stinner-surulescu-uatay}, we will use the following concept of weak solution,
where we formally rewrite $\nabla m=2\sqrt{1+m} \cdot \nabla \sqrt{1+m}$ in view of the compactness results proved below.
\begin{defi}\label{defi3.1}
  Let $T \in (0,\infty)$. A weak solution to \eqref{eq3.0}--\eqref{eq3.0i} consists of non-negative functions
  \begin{align*}
	& m\in L^1((0,T); L^2(\Omega))
	\quad \mbox{with} \quad
	\sqrt{1+m} \in L^2 ((0,T);W^{1,2}(\Omega)) \quad\mbox{and}\\
	& \sqrt{m} \, \nabla p, \sqrt{m} \, \nabla v \in L^2(\Omega \times (0,T)), \quad
	p, v\in L^\infty(\Omega\times (0,T)) \cap L^2((0,T);W^{1,2}(\Omega)),
	\end{align*}
  such that for all $\varphi \in C_0^\infty(\bar\Omega\times [0,T))$
  \begin{align}
	&-\int_0^T \io m\partial_t \varphi - \io m_0\varphi(\cdot,0)
	= - 2\int_0^T \io D(m,p,v) \sqrt{1+m} \; \nabla \sqrt{1+m} \cdot \nabla \varphi \nn\\
	& \hspace*{+20mm}
	+ \int_0^T \io \frac{c_1 v}{1+v} m\nabla v\cdot \nabla \varphi  
	- \int_0^T \io c_2 m\nabla p\cdot \nabla \varphi+ \int_0^T \io (\lambda p - \gamma m) \varphi, \label{eqd1.1} \\
  &-\int_0^T \io p\partial_t \varphi - \io p_0\varphi(\cdot,0)
	= \int_0^T \io \Big\{ \mu p \left( 1- (m+p) - \eta_1 v \right) - \lambda p + \gamma m  \Big\} \varphi,
	\label{eqd1.2} \\
	&-\int_0^T \io v\partial_t \varphi - \io v_0\varphi(\cdot,0)
	= \int_0^T \io \Big\{ - \alpha mv - \beta pv + \mu_v v (1-v) \Big\} \varphi \label{eqd1.3} 
  \end{align}
  are fulfilled. $(m,p,v)$ is called a global weak solution to \eqref{eq3.0}--\eqref{eq3.0i}, if it is a weak solution 
	in $\Omega\times (0,T)$ for all $T>0$.
\end{defi}

The main result of this section is the global existence in the weak sense defined above.
\begin{theo}\label{theo3.2}
Let $\Omega \subset \R^n$ be a bounded domain with smooth boundary, $n \le 3$, and let \eqref{eq3.1} and \eqref{eq3.2}
be satisfied along with $0 < c_2 < \frac{4\gamma C_2}{\mu A^2}$, where $C_2 >0$ is defined in \eqref{eq3.2} 
according to $A$ and $L$ from Lemma~\ref{lem3.4} below. Then there exists a global weak solution to \eqref{eq3.0}--\eqref{eq3.0i}
in the sense of Definition~\ref{defi3.1} satisfying in addition
$$m \in L^\infty ((0,\infty), L^1 (\Omega)), \quad p,v \in L^\infty (\Omega \times (0,\infty)).$$ 
\end{theo}
This result has been proved in the Master thesis \cite{lenz19} in case of $c_1=1$ and is in fact valid for more general
coefficients $\chi_2$. Apart from the restriction $\chi_2 = \chi_2(p)$ satisfying $0 < \chi_2(p) \le c_2$ in $[0,\infty)$ 
for some $c_2$ small enough, in particular
$$\chi_2^{\prime\prime}(p) \ge 2 \frac{(\chi_2^\prime(p))^2}{\chi_2(p)}, \qquad p \ge 0,$$
has to be satisfied (see the proof of \cite[Lemma~2.12]{lenz19}). Hence, $\chi_2$ has to be convex, which is not true for the 
favorable choice of an increasing $\chi_2$
with saturation for large $p$. It remains open whether the global existence can be proved for a class of saturating 
functions $\chi_2$. 

Our proof of Theorem~\ref{theo3.2} is based on the method from \cite{stinner-surulescu-uatay} where the method from
\cite{stinner_surulescu_winkler_14} has been adapted to the case of a splitted cancer cell population in presence 
of haptotaxis. As compared to the previous settings the repellent taxis described by the term 
$+ \nabla \cdot \left( c_2 m \nabla p \right)$ is included. The main difficulty is 
that now both ODEs for $p$ and $v$ are coupled to the PDE for $m$ via a taxis mechanism. Although the new term describes 
a repellent taxis, it seems to imply further restrictions concerning global existence. The smallness condition on $c_2$ in our
method stems from the presence of $m$ in the second equation of \eqref{eq3.0} in two terms with different sign (see the comment
before Lemma~\ref{lem3.10}), while
in the third equation of \eqref{eq3.0} $m$ only appears with the favorable sign and no smallness condition on $c_1$ is needed 
here. To the best of our knowledge this is the first global existence result in the presence of attracting and repellent taxis 
terms coupled to two different ODEs.\\[-2ex]

In order to prove the global existence result, for $\eps \in (0,1)$ we first introduce the following regularized approximations
of \eqref{eq3.0}--\eqref{eq3.0i}
\begin{equation}\label{eq3.0eps}
	\left\{ \begin{array}{ll}
	\partial _t m_\eps =\nabla \cdot \left( D(m_\eps, p_\eps, v_\eps) \nabla m_\eps \right) 
	- \nabla \cdot \left( \frac{c_1 v_\eps}{1+v_\eps} m_\eps \nabla v_\eps 
	\right) & \\[1mm]
	\hspace*{+15mm} +\nabla \cdot \left( c_2 m_\eps \nabla p_\eps \right)+ \lambda p_\eps - \gamma m_\eps - \eps m_\eps^\theta, 
	&\; x \in \Omega, \, t >0, \\[1mm]
	\partial_t p_\eps = \eps \Delta p_\eps + \mu p_\eps \left( 1- (m_\eps+p_\eps) - \eta_1 v_\eps \right)
	- \lambda p_\eps + \gamma m_\eps, 
	&\; x \in \Omega, \, t >0,  \\[1mm]
	\partial _t v_\eps = \eps \Delta v_\eps - \alpha m_\eps v_\eps - \beta p_\eps v_\eps + \mu_v v_\eps (1-v_\eps), 
	&\; x \in \Omega, \, t >0,  \\[2mm]
	\partial _{\nu }m_\eps= \partial _{\nu }p_\eps=\partial_{\nu }v_\eps=0, &\; x\in\pO, \, t>0, \\[2mm]
	m_\eps(x,0)=m_{0\eps}(x), \quad p_\eps(x,0)=p_{0\eps}(x), \quad
	v_\eps(x,0)=v_{0\eps}(x), & \; x\in\Omega,
	\end{array} \right.
\end{equation}
where $\theta>\max\{2,n\}$ is a fixed parameter and the families of functions
$m_{0\eps}$, $p_{0\eps}$, and $v_{0\eps}$, $\eps \in (0,1)$, are assumed to satisfy
\begin{equation}\label{eq3_init_eps}
	\begin{array}{l} 
	m_{0\eps}, p_{0\eps}, v_{0\eps} \in C^3(\bar\Omega), \quad 
	m_{0\eps} >0 , \quad p_{0\eps} >0, \quad v_{0\eps} >0 \quad\mbox{in } \bar{\Omega}, \\[1mm]
	\partial _{\nu }m_{0\eps} =\partial _{\nu } p_{0\eps} =\partial _{\nu }v_{0\eps}
	=0 \quad \mbox{on } \partial\Omega 
	\end{array}
\end{equation}
for all $\eps \in (0,1)$ and
\begin{equation}\label{eq3_init_approx}
	m_{0\eps} \to m_0 \quad \mbox{in } C^0(\bar\Omega), \quad
	p_{0\eps} \to p_0 \quad\mbox{and}\quad
	v_{0\eps} \to v_0
	\quad \mbox{in } W^{1,2}(\Omega) \cap C^0(\bar\Omega)
\end{equation}
as $\eps\searrow 0$.\\[-2ex]

Throughout this section we assume that \eqref{eq3.1}, \eqref{eq3.2}, \eqref{eq3_init_eps}, and \eqref{eq3_init_approx}
are fulfilled. Then for each of the approximate problems \eqref{eq3.0eps} we prove the global existence of a classical
solution in Section~\ref{sec3.1}. In Section~\ref{sec3.2} we construct an entropy-type functional for \eqref{eq3.0eps}
which is quasi-dissipative in a certain sense. The latter enables us to prove appropriate compactness properties and
the existence of a global weak solution to the original problem in Section~\ref{sec3.3}.
\subsection{Global existence for the approximate problems}\label{sec3.1}
In a first step we obtain the local existence for \eqref{eq3.0eps}.
\begin{lem}\label{lem3.3}
For any $\eps \in (0,1)$ there exist $T_\eps \in (0,\infty]$ as well as positive functions $m_\eps, p_\eps, v_\eps \in
C^{2,1} (\overline{\Omega} \times [0,T_\eps))$ solving \eqref{eq3.0eps} in the classical sense in $\Omega \times (0,T_\eps)$.
Furthermore, if $T_\eps < \infty$ is fulfilled, then we have for all $\beta \in (0,1)$
\begin{equation}\label{eq3.3.1}
  \limsup\limits_{t\nearrow T_\eps}  \left\{
	\|m_\eps(\cdot,t)\|_{C^{2+\beta}(\overline{\Omega})}
	+ \|p_\eps(\cdot,t)\|_{C^{2+\beta}(\overline{\Omega})}
	+ \|v_\eps(\cdot,t)\|_{C^{2+\beta}(\overline{\Omega})} \right\} < \infty.
\end{equation}	
\end{lem}
\proof A detailed proof is given in \cite[Theorem~2.2]{lenz19}. It relies on the fixed point argument from
  \cite[Lemma~3.1]{stinner_surulescu_winkler_14} with the modifications described in \cite[Lemma~3.4]{stinner-surulescu-uatay}.
\qed

The next lemma contains several elementary estimates which are uniform with respect to $\eps \in (0,1)$.
\begin{lem}\label{lem3.4}
For any $\eps \in (0,1)$, the following estimates are satisfied:
\begin{align}
  & 0 < p_\eps (x,t) \le A := \max \left\{ \sup\limits_{\eps \in (0,1)} \| p_{0\eps} \|_{L^\infty (\Omega)}, 
		1- \frac{\lambda}{\mu}, \frac{\gamma}{\mu} \right\}, \quad x \in \bar{\Omega}, \, t \in [0,T_\eps), \label{eq3.4.1} \\
	& 0 < v_\eps (x,t) \le L := \max \left\{ \sup\limits_{\eps \in (0,1)} \| v_{0\eps} \|_{L^\infty (\Omega)}, 1
		\right\}, \quad x \in \bar{\Omega}, \, t \in [0,T_\eps), \label{eq3.4.2} \\
	& \int_\Omega m_\eps (x,t) dx \le B := \max \left\{ \sup\limits_{\eps \in (0,1)} \int_\Omega m_{0\eps}, \frac{\lambda A 
		|\Omega|}{\gamma} \right\}, \quad t \in (0, T_\eps), \label{eq3.4.3} \\
		& \eps \int_t^{t+1} \int_\Omega m_\eps^\theta (x,s) dx ds \le B + \lambda A |\Omega|, \quad t \in (0, T_\eps -1).
		\label{eq3.4.4}	
\end{align}
\end{lem}
\proof Since Lemma~\ref{lem3.3} provides the strict positivity of all solution components, the comparison principle
  applied to the second and third equation of \eqref{eq3.0eps} yields \eqref{eq3.4.1} and \eqref{eq3.4.2}, where $A$, $L$
	and $B$ are finite due to \eqref{eq3_init_approx}. Then, integrating the first equation of \eqref{eq3.0eps}, we have
	$$\frac{d}{dt} \int_\Omega m_\eps \le \lambda A |\Omega| - \gamma \int_\Omega m_\eps - \eps \int_\Omega m_\eps^\theta,
		\quad t \in (0, T_\eps),$$
	from which we deduce \eqref{eq3.4.3} and, after another integration, \eqref{eq3.4.4}.	
\qed
 
Next we obtain the global existence for the approximate problems like in \cite[Section~3.3]{stinner_surulescu_winkler_14}.
For the sake of completeness we provide an outline of the proof focusing on the influence of the new taxis term.
\begin{lem}\label{lem3.5}
For each $\eps \in (0,1)$ the solution to \eqref{eq3.0eps} exists globally in time with $T_\eps = \infty$.
\end{lem}
\proof We fix $\eps \in (0,1)$ as well as $T>0$ and define $\hte := \min \{T, T_\eps \}$. Lemma~\ref{lem3.4} implies that
  $m_\eps \in L^\theta (\Omega \times (0, \hte))$ and therefore also that
	$g_\eps := \mu p_\eps \left( 1- (m_\eps+p_\eps) - \eta_1 v_\eps \right) - \lambda p_\eps + \gamma m_\eps$ is bounded in
	$L^\theta (\Omega \times (0, \hte))$.	In view of $\theta > \max \{2,n\}$, we have  the continuous embedding 
	$W^{2,\theta} (\Omega) \hookrightarrow W^{1,\infty} (\Omega)$ and may apply well-known results on maximal
	Sobolev regularity to the second equation of \eqref{eq3.0eps} (see \cite{hieber-pruess}) to obtain
  \begin{equation}\label{eq3.5.1}
	  \int_0^{\hte} \| \nabla p_\eps (\cdot,t)\|_{L^\infty(\Omega)}^2 dt
		\le C(T) \left( 1+ \int_0^{\hte} \| p_\eps (\cdot,t)\|_{W^{2,\theta}(\Omega)}^\theta dt \right) \le C_3 (\eps,T). 
	\end{equation}
	Similarly, by choosing $C_3(\eps,T)$ large enough, we further have
	\begin{equation}\label{eq3.5.2}
	  \int_0^{\hte} \| \nabla v_\eps (\cdot,t)\|_{L^\infty(\Omega)}^2 dt \le C_3 (\eps,T). 
	\end{equation}
	Then with $A$ and $L$ as provided by Lemma~\ref{lem3.4}, we have $D(m_\eps,p_\eps, v_\eps) \ge C_2 >0$ for all 
	$(x,t) \in \Omega \times (0,\hte)$ due to \eqref{eq3.2}. Hence, multiplying the first equation in \eqref{eq3.0eps} by $m_\eps^{q-1}$, dropping 
	non-negative terms and using integration by parts, Young's inequality and Lemma~\ref{lem3.4}, we have for fixed $q>1$
	\begin{eqnarray*}
	  \frac{1}{q} \frac{d}{dt} \int_\Omega m_\eps^q &\le& - (q-1) C_2\int_\Omega m_\eps^{q-2} |\nabla m_\eps|^2 
		+ (q-1) \int_\Omega \frac{c_1 v_\eps}{1+v_\eps} m_\eps^{q-1} \nabla v_\eps \cdot \nabla m_\eps \\
		& & - (q-1) \int_\Omega c_2 m_\eps^{q-1} \nabla p_\eps \cdot \nabla m_\eps
		+ \lambda A \int_\Omega m_\eps^{q-1} \\
		&\le& \frac{(q-1)}{ C_2} \left( c_1^2 \|\nabla v_\eps (\cdot,t)\|_{L^\infty (\Omega)}^2
		+ c_2^2\|\nabla p_\eps (\cdot,t)\|_{L^\infty (\Omega)}^2\right) \int_\Omega m_\eps^q \\
		& & + \lambda A \left( |\Omega| + \int_\Omega m_\eps^q \right), \qquad t \in (0,\hte).
	\end{eqnarray*}
	Then an ODE comparison principle in conjunction with \eqref{eq3.5.1} and \eqref{eq3.5.2} implies
	\begin{equation}\label{eq3.5.3} 
	  \int_\Omega m_\eps^q (\cdot,t) \le C_4 (\eps,q,T), \qquad t \in (0,\hte),
	\end{equation}
	as well as $g_\eps \in L^\infty((0,\hte), L^q (\Omega))$ for any $q>1$ since $p_\eps$ and $v_\eps$ are bounded.	But then
	smoothing properties of the Neumann heat semigroup (see e.g. \cite[Lemma~4.1]{horstmann-winkler}) applied to the 
	second and third equation of \eqref{eq3.0eps} yield
	$$\|\nabla p_\eps (\cdot,t)\|_{L^\infty (\Omega)} + \|\nabla v_\eps (\cdot,t)\|_{L^\infty (\Omega)} 
	  \le C_5 (\eps,T), \quad t \in (0, \hte),$$
	the argument for $v_\eps$ being similar. Using this estimate and \eqref{eq3.5.3} together with parabolic H\"{o}lder and
	Schauder estimates, we finally conclude that \eqref{eq3.3.1} cannot be valid so that $T_\eps = \infty$ by
	Lemma~\ref{lem3.3} (see the proofs of \cite[Lemma~3.11]{stinner_surulescu_winkler_14} and
	\cite[Lemma~2.5]{lenz19} for details).
\qed

\subsection{An entropy-type functional}\label{sec3.2} 
An essential step towards the existence of a global weak solution to the original problem \eqref{eq3.0}--\eqref{eq3.0i}
is the following estimate stemming from an entropy-type functional. Its proof is the aim of this section.
\begin{prop}\label{prop3.6}
  Let $T>0$ and $0< c_2 < \frac{4\gamma C_2}{\mu A^2}$, where $C_2 >0$ is defined in \eqref{eq3.2} 
	according to $A$ and $L$ from Lemma~\ref{lem3.4}.  
	Then there exists a constant $C(T) >0$ such that 
	for any $\eps \in (0,1)$ the solution to \eqref{eq3.0eps} satisfies
	\begin{align}\label{eq3.6.1}
	  & \sup\limits_{t \in (0,T)} \left\{ \int_\Omega m_\eps \ln m_\eps 
		+ \int_\Omega \frac{|\nabla v_\eps|^2}{1+ v_\eps}
		+ \int_\Omega |\nabla p_\eps|^2 \right\} 
		+ \int_0^T \int_\Omega D(m_\eps, p_\eps, v_\eps) \frac{|\nabla m_\eps|^2}{m_\eps}  \nn \\
	  & + \int_0^T \io m_\eps \frac{|\nabla v_\eps|^2}{(1+v_\eps)^2}
		+ \int_0^T \io m_\eps |\nabla p_\eps|^2
    + \eps \int_0^T \int_\Omega m_\eps^\theta \ln (m_\eps +2) 
		\le C(T).
	\end{align}
\end{prop}
We prove this result by showing the existence of an entropy-type functional via a number of integral estimates.
To this end, we use the strategy from \cite[Section~3.2]{stinner-surulescu-uatay} which is an adaptation of
the strategy established in \cite[Section~4]{stinner_surulescu_winkler_14} to the setting with a splitted 
cancer cell population. The main difference here is that \eqref{eq3.0eps} contains two taxis terms which are
both coupled to ODEs. 

First we study the evolution of the first integral in \eqref{eq3.6.1} like in \cite[Lemma~3.9]{stinner-surulescu-uatay}.
\begin{lem}\label{lem3.7}
  There exists $C>0$ such that for any $\eps \in (0,1)$ and all $t>0$ we have
	\begin{align}\label{eq3.7.1}
	  & \frac{d}{dt} \int_\Omega m_\eps \ln m_\eps + \int_\Omega D(m_\eps, p_\eps, v_\eps) \frac{|\nabla m_\eps|^2}{
		m_\eps} + \frac{\eps}{2} \int_\Omega m_\eps^\theta \ln (m_\eps +2) \nn \\
		& \le \int_\Omega \frac{c_1 v_\eps}{1+v_\eps} \nabla m_\eps \cdot \nabla v_\eps 
		- \int_\Omega c_2 \nabla m_\eps \cdot \nabla p_\eps + C.
	\end{align}
\end{lem}
\proof As $m_\eps$ is positive (see Lemma~\ref{lem3.3}) we deduce from the first equation of \eqref{eq3.0eps} along with
  Lemma~\ref{lem3.4} that
	\begin{eqnarray*}
	  \frac{d}{dt} \int_\Omega m_\eps \ln m_\eps &=& \int_\Omega ( \ln m_\eps \partial_t m_\eps + \partial_t m_\eps) \\
		&=& - \int_\Omega D(m_\eps, p_\eps, v_\eps) \frac{|\nabla m_\eps|^2}{m_\eps}
		+ \int_\Omega \frac{c_1 v_\eps}{1+v_\eps} \nabla m_\eps \cdot \nabla v_\eps
		- \int_\Omega c_2 \nabla m_\eps \cdot \nabla p_\eps \\
		& & + \int_\Omega \lambda p_\eps \ln m_\eps - \int_\Omega \gamma m_\eps \ln m_\eps
		- \eps \int_\Omega m_\eps^\theta \ln m_\eps + \int_\Omega \lambda p_\eps \\
		& & - \int_\Omega \gamma m_\eps - \eps \int_\Omega m_\eps^\theta \\
		&\le& - \int_\Omega D(m_\eps, p_\eps, v_\eps) \frac{|\nabla m_\eps|^2}{m_\eps}
		+ \int_\Omega \frac{c_1 v_\eps}{1+v_\eps} \nabla m_\eps \cdot \nabla v_\eps
		- \int_\Omega c_2 \nabla m_\eps \cdot \nabla p_\eps \\
		& &	+ \lambda AB + \gamma\frac{|\Omega|}{e} 
		-  \frac{\eps}{2} \int_\Omega m_\eps^\theta \ln (m_\eps +2) + C_3 + \lambda A |\Omega| 
		\qquad\quad\mbox{for all } t>0.
	\end{eqnarray*}
  In the last step we have used
	$$\int_\Omega \lambda p_\eps \ln m_\eps \le \lambda \int_{\{m_\eps \ge 1\}} p_\eps \ln m_\eps 
	  \le \lambda A \int_{\{m_\eps \ge 1\}} m_\eps \le \lambda AB$$
	as well as $\xi \ln \xi \ge -\frac{1}{e}$ for all $\xi >0$ and that there is some $C_3 >0$ such that 
	$-\xi^\theta \ln \xi \le - \frac{1}{2} \xi^\theta \ln (\xi +2) + C_3$ holds for all $\xi >0$ (see 
	\cite[Lemma~4.2]{stinner_surulescu_winkler_14}).	
	This completes the proof.
\qed

Next, we use the third and second equation of \eqref{eq3.0eps} in order to cancel the first two terms on the
right-hand side of \eqref{eq3.7.1}. As the haptotaxis term and the equation for $v_\eps$ are the same as in
\cite{stinner-surulescu-uatay}, we may directly use the results of \cite{stinner-surulescu-uatay} for the
specific choice $\kappa_\eps(x,t) \equiv c_1$. The precise results are stated in the next two lemmas.
\begin{lem}\label{lem3.8}
  For any $\eps \in (0,1)$ we have
	\begin{eqnarray}\label{eq3.8.1}
	  \partial_t \frac{c_1 |\nabla v_\eps|^2}{1+v_\eps} 
		&\le& 2 \eps \frac{c_1}{1+v_\eps} \nabla v_\eps \cdot \nabla \Delta v_\eps 
		- \eps \frac{c_1}{(1+v_\eps)^2} |\nabla v_\eps|^2 \Delta v_\eps 
		- 2\alpha \frac{c_1 v_\eps}{1+v_\eps} \nabla m_\eps \cdot \nabla v_\eps \nn \\
		& & + \frac{\beta^2 c_1}{2 \mu_v} |\nabla p_\eps|^2 - 2\alpha c_1 m_\eps \frac{|\nabla v_\eps|^2}{(1+v_\eps)^2} 
		+ 2\mu_v c_1 \frac{|\nabla v_\eps|^2}{1+v_\eps} 
	\end{eqnarray}
	for all $x\in \Omega$, $t>0$.
\end{lem} 
\proof This is \cite[Lemma~3.10]{stinner-surulescu-uatay} with the choice $\kappa_\eps \equiv c_1$.
\qed

When \eqref{eq3.8.1} is integrated, the first two terms on the right-hand side are estimated according to the following lemma.
\begin{lem}\label{lem3.9}
  For any $T>0$ there is $C(T) >0$ such that for each $\eps\in (0,1)$ we have
  \begin{equation}\label{eq3.9.1}
	2\eps\io\frac{c_1}{1+v_\eps} \nabla v_\eps \cdot \nabla \Delta v_\eps
	- \eps \io \frac{c_1}{(1+v_\eps)^2} |\nabla v_\eps|^2 \Delta v_\eps
	\le \eps C(T) \io \frac{|\nabla v_\eps|^2}{1+v_\eps} 
  \end{equation}
	for all $t \in (0,T)$.
\end{lem}
\proof This is proved in \cite[Lemma~3.11]{stinner-surulescu-uatay} with the choice $\kappa_\eps \equiv c_1$.
\qed

In order to cancel the second term on the right-hand side of \eqref{eq3.7.1} we use ideas from the proofs of the two previous
lemmas, but now for the second equation of \eqref{eq3.0eps}. In the latter equation the switching 
term $+\gamma m_\eps$ is used to cancel the mentioned term, while the term $-\mu p_\eps m_\eps$ in the logistic proliferation
term has the opposite sign and results in an additional term containing $|\nabla m_\eps|^2$ on the right-hand side of  
\eqref{eq3.10.1}, which can only be compensated by the diffusion term in the first equation of \eqref{eq3.0eps} if $c_2$ is
sufficiently small. Only here we need the restrictions that $\chi_2$, $\lambda$, and
$\gamma$ are constant (the non-constant $\chi_2(p) := \frac{c_2}{1+p}$ 
has been investigated in \cite[Lemmas 2.11 and 2.12]{lenz19}).
\begin{lem}\label{lem3.10}
Let $\eta \in (0,2)$ be arbitrary. Then for any $\eps \in (0,1)$ we have
\begin{eqnarray}\label{eq3.10.1}
  \partial_t \left( c_2 |\nabla p_\eps|^2 \right) &\le& 2 \eps c_2 \nabla p_\eps \cdot \nabla \Delta p_\eps
	+ 2 c_2 \gamma \nabla m_\eps \cdot \nabla p_\eps + 3c_2 \mu |\nabla p_\eps|^2 \nn \\
	& & + c_2 \mu \eta_1^2 A^2 (1+L) \frac{|\nabla v_\eps|^2}{1+v_\eps} 
	+ \frac{c_2 \mu A^2}{2-\eta} \frac{|\nabla m_\eps|^2}{m_\eps} - \eta c_2 \mu m_\eps |\nabla p_\eps|^2
\end{eqnarray}
for all $x \in \Omega$, $t>0$, where $A$ and $L$ are defined in Lemma~\ref{lem3.4}.
\end{lem}
\proof As $p_\eps \in C^\infty (\Omega \times (0,\infty))$ by parabolic regularity theory, we may use the second
  equation of \eqref{eq3.0eps} as well as Lemma~\ref{lem3.4} and Young's inequality to obtain
	\begin{eqnarray*}
	  \partial_t \left( c_2 |\nabla p_\eps|^2 \right) &=& 2c_2 \nabla p_\eps \cdot \nabla (\partial_t p_\eps) \\
		&=& 2 \eps c_2 \nabla p_\eps \cdot \nabla \Delta p_\eps + 2c_2 \mu |\nabla p_\eps|^2
		- 2c_2 \mu m_\eps |\nabla p_\eps|^2 - 2c_2 \mu p_\eps \nabla m_\eps \cdot \nabla p_\eps \\
		& & - 4c_2 \mu p_\eps |\nabla p_\eps|^2 - 2c_2 \mu \eta_1 v_\eps |\nabla p_\eps|^2 
		- 2c_2 \mu \eta_1 p_\eps \nabla v_\eps \cdot \nabla p_\eps
		- 2c_2 \lambda |\nabla p_\eps|^2 \\
		& & + 2c_2 \gamma \nabla m_\eps \cdot \nabla p_\eps  \\
		&\le& 2 \eps c_2 \nabla p_\eps \cdot \nabla \Delta p_\eps + 2c_2 \gamma \nabla m_\eps \cdot \nabla p_\eps
		+ 2c_2 \mu |\nabla p_\eps|^2 - 2c_2 \mu m_\eps |\nabla p_\eps|^2 \\
		& & + 2c_2 \mu A |\nabla m_\eps \cdot \nabla p_\eps | 
		+ 2c_2 \mu \eta_1 A |\nabla v_\eps \cdot \nabla p_\eps| \\
		&\le& 2 \eps c_2 \nabla p_\eps \cdot \nabla \Delta p_\eps + 2c_2 \gamma \nabla m_\eps \cdot \nabla p_\eps
		+ 3c_2 \mu |\nabla p_\eps|^2 - \eta c_2 \mu m_\eps |\nabla p_\eps|^2 \\
		& & + \frac{c_2 \mu A^2}{2-\eta} \frac{|\nabla m_\eps|^2}{m_\eps} 
		+ c_2 \mu \eta_1^2 A^2 |\nabla v_\eps|^2 
	\end{eqnarray*}
	for all $x\in \Omega$ and $t>0$. Using the estimate $1 \le \frac{1+L}{1+v_\eps}$ from Lemma~\ref{lem3.4} in the
	last term, the claim is proved.
\qed

Similarly to the proof of Lemma~\ref{lem3.9}, we next estimate the first term on the right-hand side of \eqref{eq3.10.1}.%\eqref{eq3.9.1}.
\begin{lem}\label{lem3.11}
There is $C>0$ such that for all $\eps \in (0,1)$ we have
\begin{equation}\label{eq3.11.1}
  2 \eps c_2 \io \nabla p_\eps \cdot \nabla \Delta p_\eps \le \eps c_2 C \io |\nabla p_\eps|^2 \qquad\mbox{for all } t>0.
\end{equation}
\end{lem}
\proof We use the ideas from \cite[Lemma~3.11]{stinner-surulescu-uatay}. As $\partial \Omega$ is smooth, by 
  \cite[Lemma~4.2]{mizoguchi-souplet} there exists $c_\Omega >0$ depending only on the curvatures of $\Omega$ such that
  $\partial_\nu |\nabla p_\eps|^2 \le c_\Omega |\nabla p_\eps|^2$ on $\partial \Omega$, since 
	$p_\eps \in C^{2,1} (\overline{\Omega} \times (0,\infty))$ by Lemma~\ref{lem3.3} and $\partial_\nu p_\eps
	=0$ on $\pO$. Integrating by parts, we get
	\begin{eqnarray*}
	  2 \eps c_2 \io \nabla p_\eps \cdot \nabla \Delta p_\eps
	  &=& -2 \eps c_2 \io \left| D^2 p_\eps \right|^2 + \eps c_2 \int_{\pO} \partial_\nu |\nabla p_\eps|^2 d \sigma \\
	  &\le& -2 \eps c_2 \io \left| D^2 p_\eps \right|^2 + \eps c_2 c_\Omega \int_{\pO} |\nabla p_\eps|^2 d \sigma
		\qquad\mbox{for all } t>0.
	\end{eqnarray*}	
	Next, we fix $r \in (0, \frac{1}{2})$, and define $a := r+ 
	\frac{1}{2} \in (0,1)$. In view of the compact embedding $W^{r+ \frac{1}{2},2} (\Omega) \hookrightarrow L^2 (\partial 
	\Omega)$ 
	(see \cite[Proposition~4.22(ii) and Theorem~4.24(i)]{haroske-triebel}) and the fractional Gagliardo-Nirenberg inequality (see \cite[Lemma~2.5]{ishida-seki-yokota}) and Young's inequality, we further obtain
	\begin{eqnarray*}
	  \eps c_2 c_\Omega \int_{\pO} |\nabla p_\eps|^2 d \sigma
		&\le & \eps c_2 C_3 \| \nabla p_\eps \|_{W^{r+ \frac{1}{2},2}(\Omega)}^2 \\
		&\le& \eps c_2 C_4 \left( \| \nabla |\nabla p_\eps| \|_{L^2 (\Omega)}^{2 a}
		\| \nabla p_\eps \|_{L^2 (\Omega)}^{2(1-a)}
		+ \| \nabla p_\eps \|_{L^2 (\Omega)}^2 \right) \\
		&\le& 2 \eps c_2 \| \nabla |\nabla p_\eps| \|_{L^2 (\Omega)}^2
		+ \eps c_2 C_5 \| \nabla p_\eps \|_{L^2 (\Omega)}^2 \\
		&\le& 2 \eps c_2 \| D^2 p_\eps \|_{L^2 (\Omega)}^2
		+ \eps c_2 C_5 \| \nabla p_\eps \|_{L^2 (\Omega)}^2 \qquad\mbox{for all } t>0,
	\end{eqnarray*}
	where the latter estimate follows from 
	$\nabla |\nabla p_\eps | = \frac{D^2 p_\eps \cdot \nabla p_\eps}{|\nabla p_\eps |}$. 
	A combination with the previous estimate completes the proof.
\qed

Now we are in a position to prove Proposition~\ref{prop3.6}, if $c_2$ is sufficiently small. We use the ideas from
\cite[Lemma~4.1]{stinner_surulescu_winkler_14} and \cite[Proposition~3.8]{stinner-surulescu-uatay}.

\textbf{Proof of Proposition~\ref{prop3.6}.} \quad We fix $\eps \in (0,1)$ and $T>0$ and note that none of the constants
  $C_i$ below depends on $\eps$, while we indicate dependence on $T$.
  Combining \eqref{eq3.8.1} and \eqref{eq3.9.1}, we obtain
	\begin{eqnarray*}
	  \frac{d}{dt} \io \frac{c_1 |\nabla v_\eps|^2}{1+v_\eps} + 2\alpha c_1 \io m_\eps \frac{|\nabla v_\eps|^2}{(1+v_\eps)^2} 
		&\le& - 2\alpha \io \frac{c_1 v_\eps}{1+v_\eps} \nabla m_\eps \cdot \nabla v_\eps \nn \\
		& & + C_3(T) \left( \io |\nabla p_\eps|^2 + \io \frac{|\nabla v_\eps|^2}{1+v_\eps} \right)
	\end{eqnarray*}
  for all $t \in (0,T)$ with some $C_3(T) >0$. Multiplying this inequality by $\frac{1}{2\alpha}$ and adding it to
	\eqref{eq3.7.1}, we get $C_4 >0$ such that
	\begin{eqnarray}\label{eq3.6.2}
	  & & \hspace*{-20mm} \frac{d}{dt} \left\{ \int_\Omega m_\eps \ln m_\eps 
		+ \frac{c_1}{2\alpha} \io \frac{|\nabla v_\eps|^2}{1+v_\eps}\right\} 
		+ \int_\Omega D(m_\eps, p_\eps, v_\eps) \frac{|\nabla m_\eps|^2}{m_\eps} 
		+ c_1 \io m_\eps \frac{|\nabla v_\eps|^2}{(1+v_\eps)^2} \nn \\
		& & + \frac{\eps}{2} \int_\Omega m_\eps^\theta \ln (m_\eps +2) \nn \\
		& \le& - \int_\Omega c_2 \nabla m_\eps \cdot \nabla p_\eps 
		+ \frac{C_3(T)}{2\alpha} \left( \io |\nabla p_\eps|^2 + \io \frac{|\nabla v_\eps|^2}{1+v_\eps} \right) + C_4
	\end{eqnarray}
	for all $t \in (0,T)$. Next, let $C_2 >0$ be as defined in \eqref{eq3.2} according to $A$ and $L$ from Lemma~\ref{lem3.4}. 
	Since $c_2 < \frac{4\gamma C_2}{\mu A^2}$, there is $\eta \in (0,2)$ such that 
	$\xi := 1 - \frac{c_2 \mu A^2}{2 \gamma (2-\eta) C_2} >0$
	is fulfilled. Choosing this $\eta$ in \eqref{eq3.10.1} and combining the latter with \eqref{eq3.11.1}, we have
	\begin{eqnarray*}
	  \frac{d}{dt} \io c_2 |\nabla p_\eps|^2  + \eta c_2 \mu \io m_\eps |\nabla p_\eps|^2
		&\le& 2 c_2 \gamma \io \nabla m_\eps \cdot \nabla p_\eps 
		+ \frac{c_2 \mu A^2}{2-\eta} \io \frac{|\nabla m_\eps|^2}{m_\eps} \\
		& & + C_5 \left( \io |\nabla p_\eps|^2 +  \io \frac{|\nabla v_\eps|^2}{1+v_\eps} \right) 
	\end{eqnarray*}
	for all $t \in (0,T)$ with some $C_5 >0$. Multiplying this inequality by $\frac{1}{2\gamma}$ and adding it to \eqref{eq3.6.2},
	in view of \eqref{eq3.2} we have
	\begin{eqnarray}\label{eq3.6.3}
	  & & \hspace*{-20mm} \frac{d}{dt} \left\{ \int_\Omega m_\eps \ln m_\eps 
		+ \frac{c_1}{2\alpha} \io \frac{|\nabla v_\eps|^2}{1+v_\eps}
		+ \frac{c_2}{2\gamma} \io |\nabla p_\eps|^2 \right\} 
		+ \xi \int_\Omega D(m_\eps, p_\eps, v_\eps) \frac{|\nabla m_\eps|^2}{m_\eps} \nn \\
		& & + c_1 \io m_\eps \frac{|\nabla v_\eps|^2}{(1+v_\eps)^2} 
		+ \frac{\eta c_2 \mu}{2\gamma} \io m_\eps |\nabla p_\eps|^2
		+ \frac{\eps}{2} \int_\Omega m_\eps^\theta \ln (m_\eps +2) \nn \\
		& \le& C_6(T) \left( \io |\nabla p_\eps|^2 + \io \frac{|\nabla v_\eps|^2}{1+v_\eps} \right) + C_4
	\end{eqnarray}
	for all $t \in (0,T)$ with some $C_6 (T) >0$. Defining the non-negative functions
	\begin{align*}
	  & \Eeps (t) := \int_\Omega m_\eps \ln m_\eps 
		+ \frac{c_1}{2\alpha} \io \frac{|\nabla v_\eps|^2}{1+v_\eps}
		+ \frac{c_2}{2\gamma} \io |\nabla p_\eps|^2 + \frac{2|\Omega |}{e}, \\
		& \Deps (t) := \xi \int_\Omega D(m_\eps, p_\eps, v_\eps) \frac{|\nabla m_\eps|^2}{m_\eps}
		+ c_1 \io m_\eps \frac{|\nabla v_\eps|^2}{(1+v_\eps)^2} 
		+ \frac{\eta c_2 \mu}{2\gamma} \io m_\eps |\nabla p_\eps|^2 \\
		& \hspace*{+15mm} + \frac{\eps}{2} \int_\Omega m_\eps^\theta \ln (m_\eps +2) 
	\end{align*}
	for $t \ge 0$, we conclude from \eqref{eq3.6.3} that there exists $C_7(T)>0$ such that
	$$\frac{d}{dt} \Eeps (t) + \Deps (t) \le C_7 (T) \Eeps(t) \qquad\mbox{for all } t \in (0,T).$$
	As $\sup_{\eps \in (0,1)} \Eeps (0)$ is finite by \eqref{eq3_init_approx}, in view of the non-negativity of $\Deps$ and
	$\Eeps$ two integrations show that there is $C_8(T)>0$ such that
	$$\sup\limits_{t \in (0,T)} \Eeps (t) \le C_8 (T) \quad\mbox{and}\quad \int_0^T \Deps (t) dt \le C_8(T).$$
  This completes the proof of \eqref{eq3.6.1}. 
\qed

\subsection{Global weak solution to the original problem \eqref{eq3.0}--\eqref{eq3.0i}}\label{sec3.3}
We will derive appropriate compactness properties for the solutions to \eqref{eq3.0eps} which imply their convergence
to a global weak solution of the original problem \eqref{eq3.0}--\eqref{eq3.0i}. As a basis we will use estimate 
\eqref{eq3.6.1} stemming from the entropy-type functional $\Eeps$. The proofs mainly rely on the ideas from
\cite[Section~3.3]{stinner-surulescu-uatay} where the method from \cite[Section~5]{stinner_surulescu_winkler_14}
has been adapted to the case of a splitted cancer cell population. 

First we prove properties of $m_\eps$ similarly to \cite[Lemma~3.14]{stinner-surulescu-uatay}.
\begin{lem}\label{lem3.12}
Let $T >0$ be arbitrary. There exists a constant $C(T) >0$ such that for any $\eps \in (0,1)$ 
	\begin{equation}\label{eq3.12.1}
	  \int_0^T \| \sqrt{1+m_\eps (\cdot,t)}  \|_{W^{1,2} (\Omega)}^2 dt \le C(T)
	\end{equation}	
	is satisfied. Moreover, $(\sqrt{1+m_\eps})_{\eps \in (0,1)}$ is strongly precompact in $L^2((0,T); L^q(\Omega))$ for any
	$q \in (1,6)$ and $(m_\eps)_{\eps \in (0,1)}$ is strongly precompact in $L^1((0,T); L^2(\Omega))$.
\end{lem}
\proof In view of \eqref{eq3.2}, Lemma~\ref{lem3.4}, and Proposition~\ref{prop3.6} there is $C_3(T) >0$ such that
  \begin{eqnarray*}
	  \int_0^T \| \sqrt{1+m_\eps}  \|_{W^{1,2} (\Omega)}^2 
		&=& \int_0^T\io (1+ m_\eps) + \frac{1}{4} \int_0^T \io \frac{|\nabla m_\eps|^2}{1+ m_\eps} \\
		&\le& T(|\Omega| + B) + \frac{1}{4 C_2} \int_0^T \io D (m_\eps, p_\eps, v_\eps)
		\frac{|\nabla m_\eps|^2}{m_\eps} 
		\le C_3(T)
	\end{eqnarray*}
	and \eqref{eq3.12.1} is proved. Next, we fix $k \in \N$ such that $k > \frac{n+2}{2}$ and obtain
	\begin{equation}\label{eq3.12.2}
	  \int_0^T \| \partial_t \sqrt{1+m_\eps (\cdot,t)}  \|_{(W_0^{k,2} (\Omega))^\ast} dt \le C_4(T)
	\end{equation}
	with some $C_4(T)>0$ like in the proof of \cite[Lemma~3.14]{stinner-surulescu-uatay}. Indeed, for fixed 
	$\Psi \in C_0^\infty (\Omega)$ we get from the first equation of \eqref{eq3.0eps} and integration by parts
	\begin{eqnarray}\label{eq3.12.3}
	  & & 2 \int_0^T \io \partial_t \sqrt{1+m_\eps} \Psi \nn \\
		&=& \frac{1}{2} \int_0^T \io \frac{D(m_\eps, p_\eps, v_\eps)}{(1+m_\eps)^{\frac{3}{2}}} |\nabla m_\eps|^2
		\Psi 
		- \int_0^T \io \frac{D(m_\eps, p_\eps, v_\eps)}{(1+m_\eps)^{\frac{1}{2}}} \nabla m_\eps \cdot \nabla \Psi 
		\nn \\
		& & - \frac{1}{2} \int_0^T \io \frac{c_1 v_\eps m_\eps}{(1+v_\eps)(1+m_\eps)^{\frac{3}{2}}} \nabla m_\eps \cdot
		\nabla v_\eps \Psi 
		+ \int_0^T \io \frac{c_1 v_\eps m_\eps}{(1+v_\eps)(1+m_\eps)^{\frac{1}{2}}} \nabla v_\eps \cdot \nabla \Psi \nn \\
		& & + \frac{1}{2} \int_0^T \io \frac{c_2 m_\eps}{(1+m_\eps)^{\frac{3}{2}}} \nabla m_\eps \cdot
		\nabla p_\eps \Psi 
		- \int_0^T \io \frac{c_2 m_\eps}{(1+m_\eps)^{\frac{1}{2}}} \nabla p_\eps \cdot \nabla \Psi \nn \\
		& & + \int_0^T \io \left(\lambda p_\eps - \gamma m_\eps - \eps m_\eps^\theta \right) \frac{\Psi}{\sqrt{1+m_\eps}} \nn \\
		&\le& C_5 (T) \|\Psi\|_{W^{1,\infty}(\Omega)}.	
	\end{eqnarray}	
  Therein, all the terms except those containing $c_2$ can be estimated like in the proof of 
	\cite[Lemma~3.14]{stinner-surulescu-uatay},
	while for the remaining terms \eqref{eq3.2}, Lemma~\ref{lem3.4}, and Proposition~\ref{prop3.6} imply
	\begin{eqnarray*}
	  & & \frac{1}{2} \int_0^T \io \frac{c_2 m_\eps}{(1+m_\eps)^{\frac{3}{2}}} \nabla m_\eps \cdot
		\nabla p_\eps \Psi 
		- \int_0^T \io \frac{c_2 m_\eps}{(1+m_\eps)^{\frac{1}{2}}} \nabla p_\eps \cdot \nabla \Psi \\
		&\le& \frac{1}{4} \| \Psi\|_{L^\infty (\Omega)} \int_0^T \io \left( c_2^2 |\nabla p_\eps|^2 
			+ \frac{|\nabla m_\eps|^2}{1+m_\eps} \right) 
			+ \frac{1}{2} \| \nabla \Psi\|_{L^\infty(\Omega)} \int_0^T \io \left( m_\eps + c_2^2 |\nabla p_\eps|^2 \right)\\
		&\le& \frac{1}{4} \| \Psi\|_{W^{1,\infty} (\Omega)} \left( 3 c_2^2 \int_0^T \io  |\nabla p_\eps|^2 
			+ \frac{1}{C_2} \int_0^T \io  D(m_\eps, p_\eps, v_\eps) \frac{|\nabla m_\eps|^2}{m_\eps} + 2TB \right)	\\
		&\le& C_6 (T) \|\Psi\|_{W^{1,\infty}(\Omega)}.	
	\end{eqnarray*}
	Since the embedding $W_0^{k,2}(\Omega) \hookrightarrow W^{1,\infty}(\Omega)$ is continuous due to $k> \frac{n+2}{2}$, 
	\eqref{eq3.12.2} follows from \eqref{eq3.12.3} and
	$$\int_0^T \| \partial_t \sqrt{1+m_\eps (\cdot,t)}  \|_{(W_0^{k,2} (\Omega))^\ast} dt
	= \int_0^T \sup\limits_{\Psi \in C_0^\infty (\Omega), \|\Psi\|_{W_0^{k,2} (\Omega)} \le 1} 
	\io \partial_t \sqrt{1+m_\eps(\cdot,t)} \Psi.$$
	For $q \in (1,6)$  we deduce the strong precompactness of $(\sqrt{1+m_\eps})_{\eps \in (0,1)}$ in
	$L^2((0,T); L^q (\Omega))$ from the Aubin-Lions Lemma (see e.g. Theorem~2.3 and Remark~2.1 in Chapter III of
	\cite{temam}), since \eqref{eq3.12.1} and \eqref{eq3.12.2} yield the boundedness of 
	$(\sqrt{1+m_\eps})_{\eps \in (0,1)}$ in $L^2 ((0,T); W^{1,2} (\Omega))$ and of
	$(\partial_t \sqrt{1+m_\eps})_{\eps \in (0,1)}$ in $L^1 ((0,T); (W_0^{k,2} (\Omega))^\ast)$, while 
	$W^{1,2} (\Omega) \hookrightarrow L^q (\Omega)$ is compact and $L^q (\Omega) \hookrightarrow (W_0^{k,2} (\Omega))^\ast$
	is continuous due to $n \le 3$ and $k > \frac{n+2}{2}$. Finally, the strong precompactness of $(m_\eps)_{\eps \in (0,1)}$ 
	in $L^1((0,T); L^2(\Omega))$ is a consequence of the case $q=4$ and $m_\eps \ge 0$.
\qed

For $p_\eps$ and $v_\eps$ we have the following result.
\begin{lem}\label{lem3.13}
  Let $T >0$ be arbitrary. There exists a constant $C(T) >0$ such that for any $\eps \in (0,1)$ 
	\begin{equation}\label{eq3.13.1}
	  \sup\limits_{t \in (0,T)} \left\{ \io |\nabla p_\eps (\cdot,t)|^2 
		+ \io |\nabla v_\eps (\cdot,t)|^2 \right\} \le C(T)
	\end{equation}	
	is fulfilled. Furthermore, $(p_\eps)_{\eps \in (0,1)}$ and $(v_\eps)_{\eps \in (0,1)}$ 
	are strongly precompact in $L^2(\Omega \times (0,T))$.
\end{lem}
\proof \eqref{eq3.13.1} is a consequence of Proposition~\ref{prop3.6}, since Lemma~\ref{lem3.4}
  implies $|\nabla v_\eps|^2 \le (1+L) \frac{|\nabla v_\eps|^2}{1+v_\eps}$. Then \eqref{eq3.13.1}
	and Lemma~\ref{lem3.4} yield the boundedness of $(p_\eps)_{\eps \in (0,1)}$ and $(v_\eps)_{\eps \in (0,1)}$
	in $L^2 ((0,T); W^{1,2} (\Omega))$, while the boundedness of the respective time derivatives in
	$L^1 ((0,T); (W_0^{k,2} (\Omega))^\ast)$ for $k > \frac{n+2}{2}$ can be shown similarly to \eqref{eq3.12.3}. 	
  Finally, like in the proof of Lemma~\ref{lem3.12} we obtain the claimed compactness properties from the
	Aubin-Lions Lemma.	
\qed 

Finally, we are in a position to prove the existence of a global weak solution to \eqref{eq3.0}--\eqref{eq3.0i}
by relying on the ideas of the proofs of \cite[Theorem~1.1]{stinner_surulescu_winkler_14} and
\cite[Theorem~3.2]{stinner-surulescu-uatay}.\\

\textbf{Proof of Theorem~\ref{theo3.2}.} \quad 
  Due to Lemmas \ref{lem3.4}, \ref{lem3.12}, and \ref{lem3.13} there exist non-negative functions $m$, $p$, and $v$ satisfying the regularity
  properties stated in Definition~\ref{defi3.1} and claimed in Theorem~\ref{theo3.2} such that along an appropriate sequence
  $\eps = \eps_j \searrow 0$ as $j \to \infty$ for any $T>0$
	\begin{equation}\label{eq3.2.1}
	  \begin{array}{l}
	  l_\eps \to l \qquad\mbox{strongly in } L^2(\Omega \times (0,T)) \mbox{ and a.e. in } \Omega \times (0,\infty),
		\quad\mbox{for } l \in \{\sqrt{1+m},p,v \}, \\[2mm]
		m_\eps \to m \qquad\mbox{strongly in } L^1((0,T);L^2(\Omega)) \mbox{ and a.e. in } \Omega \times (0,\infty), \\[2mm]
		\nabla \sqrt{1+m_\eps} \wto \nabla \sqrt{1+m}, \quad \nabla p_\eps \wto \nabla p, \quad\mbox{and}\quad \nabla v_\eps \wto \nabla v
		\quad\mbox{weakly in } L^2(\Omega \times (0,T)), \\[2mm]
		\sqrt{m_\eps} \, \nabla p_\eps \wto \sqrt{m} \, \nabla p \quad\mbox{and}\quad 
		\sqrt{m_\eps} \, \nabla v_\eps \wto \sqrt{m} \, \nabla v \qquad\mbox{weakly in } L^2(\Omega \times (0,T))
	  \end{array}
  \end{equation}
  are fulfilled. Here the last two convergences are consequences of $\sqrt{m_\eps} \to \sqrt{m}$ strongly in $L^2(\Omega \times (0,T))$ (and a.e.) 
	as well as $\nabla p_\eps \wto \nabla p$ and $\nabla v_\eps \wto \nabla v$ weakly in $L^2(\Omega \times (0,T))$, 
	since we obtain from Proposition~\ref{prop3.6} and Lemma~\ref{lem3.4} that
	$$\int_0^T \io m_\eps |\nabla p_\eps|^2 + \int_0^T \io m_\eps |\nabla v_\eps|^2 \le 
	\int_0^T \io m_\eps |\nabla p_\eps|^2 
	+ (1+L)^2 \int_0^T \io m_\eps \frac{|\nabla v_\eps|^2}{(1+v_\eps)^2} \le C_3(T)$$
	holds for all $\eps \in (0,1)$ with some $C_3(T) >0$. 
	Fixing $T>0$ and $\varphi \in C_0^\infty (\overline{\Omega} \times [0,T))$, the first equation in \eqref{eq3.0eps} along with integration by parts 
	implies that
        \begin{align}\label{eq3.2.2}
          & \hspace*{-5mm} -\int_0^T \io m_\eps \partial_t \varphi - \io m_{0\eps}\varphi(\cdot,0) \nn \\ 
	  & = - 2\int_0^T \io D(m_\eps,p_\eps,v_\eps) \sqrt{1+m_\eps} \nabla \sqrt{1+m_\eps} \cdot \nabla \varphi 
	  + \int_0^T \io \frac{c_1 v_\eps}{1+v_\eps} \sqrt{m_\eps} \sqrt{m_\eps} \, \nabla v_\eps\cdot \nabla \varphi \nn \\
	  & - \int_0^T \io c_2 \sqrt{m_\eps} \sqrt{m_\eps} \, \nabla p_\eps \cdot \nabla \varphi
	  + \int_0^T \io (\lambda p_\eps - \gamma m_\eps)\varphi - \eps  \int_0^T \io m_\eps^\theta \varphi
        \end{align}
        for all $\eps \in (0,1)$. Passing to the limit $\eps = \eps_j \searrow 0$, we deduce from \eqref{eq3.2.1}, 
	\eqref{eq3.2}, and \eqref{eq3_init_approx} that each of the
	terms in \eqref{eq3.2.2} except the last one converges to the respective term of \eqref{eqd1.1}.
	Therein we used the strong convergences $D(m_\eps,p_\eps,v_\eps) \sqrt{1+m_\eps} \to D(m,p,v) \sqrt{1+m}$ and $\frac{c_1 v_\eps}{1+v_\eps} 
	\sqrt{m_\eps} \to \frac{c_1 v}{1+v} \sqrt{m}$ in $L^2(\Omega \times (0,T))$, which are consequences of
	\cite[Lemma~5.10]{stinner_surulescu_winkler_14} as \eqref{eq3.2} and Lemma~\ref{lem3.4} yield $0 \le D(m_\eps,p_\eps,v_\eps) \le C_1$ and 
	$0 \le \frac{c_1 v_\eps}{1+v_\eps} \le c_1$.
	
	The last term in \eqref{eq3.2.2} converges to zero as $\eps \searrow 0$ like in the proof of \cite[Theorem~3.2]{stinner-surulescu-uatay}. 
	To this end, let $C_4(T)$ denote the constant from Proposition~\ref{prop3.6} and let $\eta >0$ be arbitrary. Then we fix $S>0$ large enough 
	such that $\frac{C_4(T)}{\ln(S+2)} \le \frac{\eta}{2}$ and deduce from Proposition~\ref{prop3.6} that
	\begin{eqnarray*}
	  \eps  \int_0^T \io m_\eps^\theta &=& \eps  \int_0^T \io \chi_{\{m_\eps \le S\} } m_\eps^\theta  
		+ \eps  \int_0^T \io \chi_{\{m_\eps > S\} } m_\eps^\theta \\
		&\le& \eps T |\Omega| S^\theta + \frac{\eps}{\ln (S+2)} \int_0^T \io m_\eps^\theta \ln (m_\eps +2)
		\le \frac{\eta}{2} + \frac{\eta}{2}
	\end{eqnarray*}
	for all $\eps \in (0, \eps_0)$ with $\eps_0 <1$ such that $\eps_0 T |\Omega| S^\theta \le \frac{\eta}{2}$, which shows the claimed convergence to zero.
        \eqref{eqd1.2} and \eqref{eqd1.3} can be verified in a similar way in view of
	\eqref{eq3.2.1}, \eqref{eq3_init_approx}, and \eqref{eq3.13.1}.
\qed

%------------------------------------------------------
\mysection{Numerical simulations}\label{sec:numerics}

We study numerically several variants of the model \eqref{eq:model}. In particular, we consider
constant and variable phenotypic-switch transition rates (from proliferating to migrating 
cells and vice versa), we compare our interspecies repellence approach with an acidity driven cell migration model, and investigate the effects of degenerate diffusion, anoikis cell apoptosis evasion, and ECM remodeling by tumor cells. %and of the cancer
%resistance to the anoikis apoptotic process.

All numerical experiments are conducted over the square domain $\Omega =[-2,2]\times[-2,2]$ on which we prescribe homogeneous Neumann boundary conditions. We augment the models with initial conditions corresponding to a biological time where the development of the tumor has
already started. We assume in particular that as the experiments start, the proliferating cells (PCCs) have
already given rise to a certain amount of migrating cells (MCCs). We also study the effect of tissue anisotropy and consider for the initial ECM density two
different scenarios: one where the matrix exhibits discontinuities in the form of stripes, and one
where the ECM is continuous and randomly structured. Namely, the ECM-with-stripes of initial
conditions are given by 
\begin{equation}\label{eq:ICstrps}
\left.\begin{aligned}
p_0(x_1,x_2) &= \begin{cases} 
0, & \text{if } |x_1|<0.05\text{ or } |x_2-x_1|<0.1,\\  
e^{-(x_1^2+x_2^2)/\varepsilon},  &\text{otherwise} \end{cases}\\
m_0(x_1,x_2) &=  0.0\,5p_0(x_1,x_2)\\
v_0(x_1,x_2) &= 1-m_0(x_1,x_2)-p_0(x_1,x_2)
\end{aligned}\right\}
\end{equation}
for $(x_1,x_2)\in\Omega$  and $\varepsilon=0.3$. 

The randomly-structured-ECM initial conditions are given by
\begin{equation}\label{eq:ICrnd}
\left. \begin{aligned}
p_0(x_1,x_2) &= e^{-(x_1^2+x_2^2)/\varepsilon}\\
m_0(x_1,x_2) &=  0.05 \, p_0(x_1,x_2)\\
%				v_0(x_1,x_2) &= \begin{cases} 1, & \text{if }|x_2-x_1|>0.1 \\ 1-p_0(x_1,x_2),
%&\text{otherwise}\end{cases}
\end{aligned}\right\}
\end{equation}
for $(x_1,x_2)\in\Omega$ and $\varepsilon=0.3$. The ECM in this case is a structured random matrix
produced by an inductive process described in Appendix~\ref{app:ECM}; this was first introduced in
\cite{nam.2020}. We augment the model \eqref{eq:model} with the sensitivity functions
\eqref{eq:microSteady2} and employ in all experiments, unless otherwise stated, the parameters
$D_c=0.001$, $\xi_1 = 0.4$, $\xi_2=0.1$, $k_D=1$, $\mu=0.1$, $\mu_v = 0.15$ and $\delta=0.3$. The
experiment-specific parameters are discussed in the corresponding experiment descriptions. All parameter choices are in accordance to those in \cite{ChLo06}. 

We refrain from a detailed presentation of the numerical method that we use but note that it is a
2nd order Implicit-Explicit Runge-Kutta Finite Volume (IMEX-RK FV)  scheme and refer to
Appendix~\ref {app:numerics} and \cite{kolbe_study_2016, gerisch_numerical_2017, nn.2016} for more
details. The implementations of the algorithms and the numerical simulations were conducted in
MATLAB \cite{MATLAB.2019}.

\begin{figure}[t] %exp 1 stripes
	\figurewidth=\linewidth
	\def \expSetup {TODO1A_ECM_stripes}
	\begin{tikzpicture}
	\begin{groupplot}[
	/tikz/mark size=1.5pt,
	group style={
		group name=my plots,
		group size=4 by 3,
		horizontal sep=1cm,      % <-- default: 1cm
		vertical sep=0.7cm,        % <-- default: 1cm
	},
	xmin=-2,
	xmax=2,
	ymin=-2,
	ymax=2,
	ylabel shift = 2 em,
	xtick = \empty,
	ytick= \empty,
	ticklabel style = {font=\scriptsize},
	axis line style = thick,
	axis background/.style={fill=white},
	width=.25\figurewidth,
	height=.25\figurewidth, 
	]
	\nextgroupplot[ylabel = MCCs, title = {$t=0$}]
	\addplot [forget plot] graphics [xmin=-2, xmax=2, ymin=-2, ymax=2] {\expSetup1};
	
	\nextgroupplot[title = {$t=3.3$}]
	\addplot [forget plot] graphics [xmin=-2, xmax=2, ymin=-2, ymax=2] {\expSetup2};
	
	\nextgroupplot[title = {$t=6.7$}]
	\addplot [forget plot] graphics [xmin=-2, xmax=2, ymin=-2, ymax=2] {\expSetup3};
	
	\nextgroupplot[
	title = {$t=10$},
	point meta min=0,
	point meta max=0.05,
	colormap name=MCCs,
	colorbar,
	colorbar style={at={(1.2,1)},anchor=north west}
	]
	\addplot [forget plot] graphics [xmin=-2, xmax=2, ymin=-2, ymax=2] {\expSetup4};
	\nextgroupplot[ylabel = PCCs]
	\addplot [forget plot] graphics [xmin=-2, xmax=2, ymin=-2, ymax=2] {\expSetup5};
	\nextgroupplot
	\addplot [forget plot] graphics [xmin=-2, xmax=2, ymin=-2, ymax=2] {\expSetup6};
	\nextgroupplot
	\addplot [forget plot] graphics [xmin=-2, xmax=2, ymin=-2, ymax=2] {\expSetup7};
	\nextgroupplot[
	point meta min=0,
	point meta max=1,
	colormap name=PCCs,
	colorbar,
	colorbar style={at={(1.2,1)},anchor=north west}
	]
	\addplot [forget plot] graphics [xmin=-2, xmax=2, ymin=-2, ymax=2] {\expSetup8};
	\nextgroupplot[ylabel = ECM]
	\addplot [forget plot] graphics [xmin=-2, xmax=2, ymin=-2, ymax=2] {\expSetup9};
	\nextgroupplot
	\addplot [forget plot] graphics [xmin=-2, xmax=2, ymin=-2, ymax=2] {\expSetup10};
	\nextgroupplot
	\addplot [forget plot] graphics [xmin=-2, xmax=2, ymin=-2, ymax=2] {\expSetup11};
	\nextgroupplot[
	point meta min=0,
	point meta max=1,
	colormap name=ECM,
	colorbar,
	colorbar style={at={(1.2,1)},anchor=north west}
	]
	\addplot [forget plot] graphics [xmin=-2, xmax=2, ymin=-2, ymax=2] {\expSetup12};
	\end{groupplot}
	\end{tikzpicture}%
	
	\caption{Simulation results of \textbf{Experiment~\ref{exp:cnstEMT} --- constant phenotypic switch
		rates} with the ECM-with-stripes initial conditions \eqref{eq:ICstrps}. In this and the rest of the
	simulation results, the densities of MCCs and PCCs are represented with isolines colored according to
	the displayed colorbars. The density of the ECM is visualized by a variable-intensity color that follows
	the corresponding colorbar. The MCCs, in their taxis-biased random motion, follow the gradients of the ECM and accordingly their density increases over the stripes of the ECM. The ECM is depleted by
	both cell subpopulations, which also limit its reconstruction. The
	PCCs obey a logistic-type growth and fill the free space left by the ECM and MCCs; they moreover undergo
	phenotypic transitions back-and-forth to MCCs according to the PMT and MPT rates $\lambda$ and $\gamma$.}\label{fig:exp1_stripes}	
\end{figure}
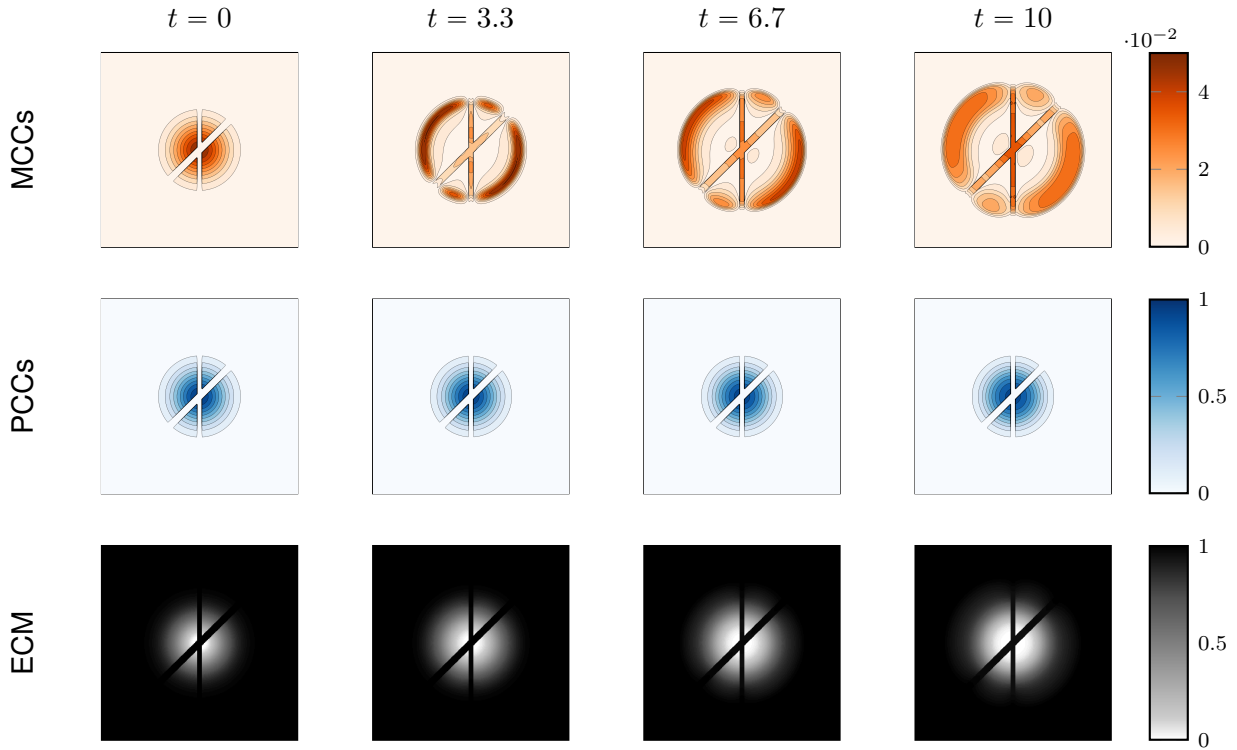

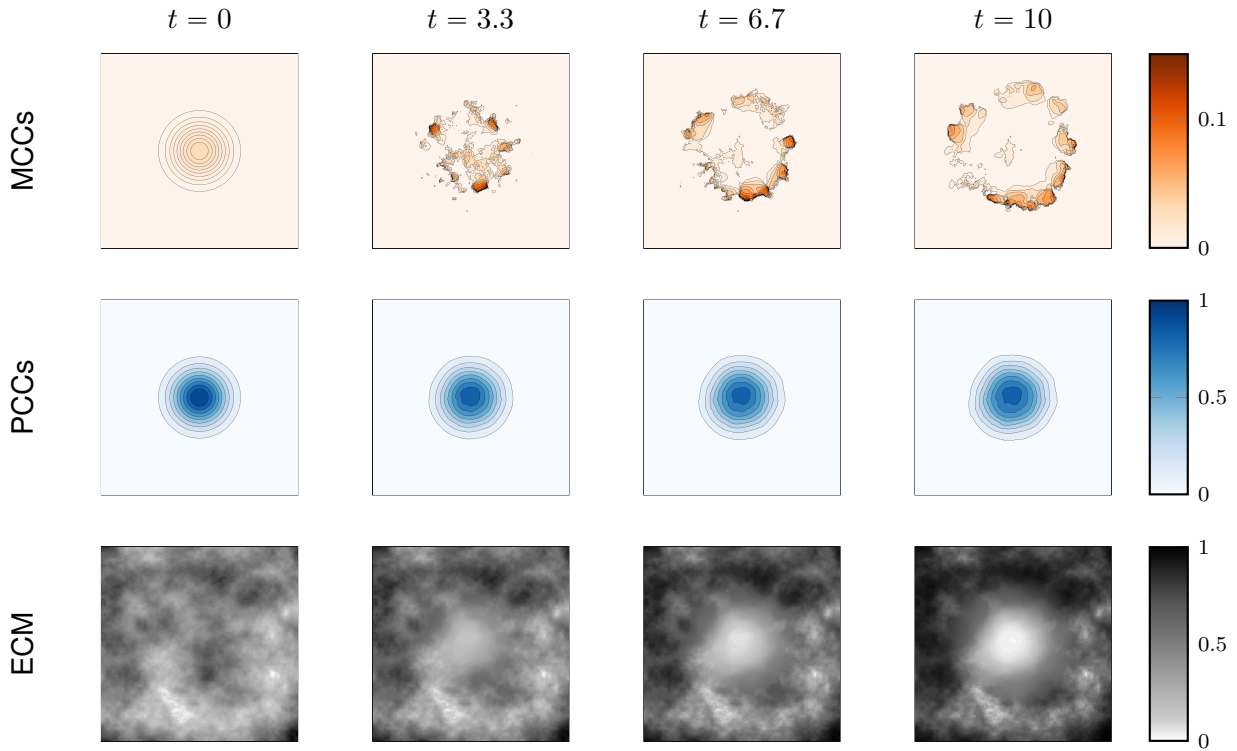
\begin{figure}[t] %exp1 rnd
	\figurewidth=\linewidth
	\def \expSetup {TODO1A_ECM_rnd}
	\begin{tikzpicture}
	\begin{groupplot}[
	/tikz/mark size=1.5pt,
	group style={
		group name=my plots,
		group size=4 by 3,
		horizontal sep=1cm,      % <-- default: 1cm
		vertical sep=0.7cm,        % <-- default: 1cm
	},
	xmin=-2,
	xmax=2,
	ymin=-2,
	ymax=2,
	ylabel shift = 2 em,
	xtick = \empty,
	ytick= \empty,
	ticklabel style = {font=\scriptsize},
	axis line style = thick,
	axis background/.style={fill=white},
	width=.25\figurewidth,
	height=.25\figurewidth, 
	]
	\nextgroupplot[ylabel = MCCs, title = {$t=0$}]
	\addplot [forget plot] graphics [xmin=-2, xmax=2, ymin=-2, ymax=2] {\expSetup1};
	\nextgroupplot[title = {$t=3.3$}]
	\addplot [forget plot] graphics [xmin=-2, xmax=2, ymin=-2, ymax=2] {\expSetup2};
	\nextgroupplot[title = {$t=6.7$}]
	\addplot [forget plot] graphics [xmin=-2, xmax=2, ymin=-2, ymax=2] {\expSetup3};
	\nextgroupplot[
	title = {$t=10$},
	point meta min=0,
	point meta max=0.15,
	colormap name=MCCs,
	colorbar,
	colorbar style={at={(1.2,1)},anchor=north west}
	]
	\addplot [forget plot] graphics [xmin=-2, xmax=2, ymin=-2, ymax=2] {\expSetup4};
	\nextgroupplot[ylabel = PCCs]
	\addplot [forget plot] graphics [xmin=-2, xmax=2, ymin=-2, ymax=2] {\expSetup5};
	\nextgroupplot
	\addplot [forget plot] graphics [xmin=-2, xmax=2, ymin=-2, ymax=2] {\expSetup6};
	\nextgroupplot
	\addplot [forget plot] graphics [xmin=-2, xmax=2, ymin=-2, ymax=2] {\expSetup7};
	\nextgroupplot[
	point meta min=0,
	point meta max=1,
	colormap name=PCCs,
	colorbar,
	colorbar style={at={(1.2,1)},anchor=north west}
	]
	\addplot [forget plot] graphics [xmin=-2, xmax=2, ymin=-2, ymax=2] {\expSetup8};
	\nextgroupplot[ylabel = ECM]
	\addplot [forget plot] graphics [xmin=-2, xmax=2, ymin=-2, ymax=2] {\expSetup9};
	\nextgroupplot
	\addplot [forget plot] graphics [xmin=-2, xmax=2, ymin=-2, ymax=2] {\expSetup10};
	\nextgroupplot
	\addplot [forget plot] graphics [xmin=-2, xmax=2, ymin=-2, ymax=2] {\expSetup11};
	\nextgroupplot[
	point meta min=0,
	point meta max=1,
	colormap name=ECM,
	colorbar,
	colorbar style={at={(1.2,1)},anchor=north west}
	]
	\addplot [forget plot] graphics [xmin=-2, xmax=2, ymin=-2, ymax=2] {\expSetup12};
	\end{groupplot}
	\end{tikzpicture}%
	
	\caption{Simulation results of \textbf{Experiment~\ref{exp:cnstEMT} --- constant phenotypic switch
		rates} with the randomly-structured ECM initial conditions \eqref{eq:ICrnd}. Through their
	haptotaxis-biased  random migration, the MCCs identify the higher ECM density regions and
	accordingly invade the surrounding environment. The ECM and PCCs exhibit similar behavior as in the
	ECM-with-stripes case shown in Figure~\ref{fig:exp1_stripes};  the ECM is depleted by the action of both
	MCCs and PCCs while the PCCs fill the space left by the ECM and MCCs.}\label{fig:exp1_rnd}	
\end{figure}

\begin{experiment}\textbf{Constant phenotypic switch rates.}\label{exp:cnstEMT}
	In this first experiment we consider constant proliferative-to-migrative (PMT)  and
	migrative-to-proliferative (MPT) phenotype transition rates 
	 $\lambda(y,\zeta)=0.01$ and
	$\gamma(y,\zeta)=0.002$, respectively. 
	The system is augmented, in turn, with the sets of initial conditions \eqref{eq:ICstrps} and \eqref{eq:ICrnd}
	and the corresponding simulation results are presented in Figures \ref{fig:exp1_stripes} and
	\ref{fig:exp1_rnd}, respectively. In both cases, the dynamics of the solutions are driven by the irregularities of
	the ECM and the interactions between the two cell types. We note that the density of the MCCs 
	increases at interfaces between higher and lower ECM density; this is due to haptotaxis. The ECM is degraded by both types of cells and, due to the particular balance between the matrix degradation and remodeling processes, the ECM infers a gradual depletion. The PCCs proliferate and due to the logistic type constraint they obey, they fill
	the ``free space'' left by the MCCs and the ECM. 
\end{experiment}

\begin{figure} % exp 2   rates
	\includegraphics[scale = 0.66]{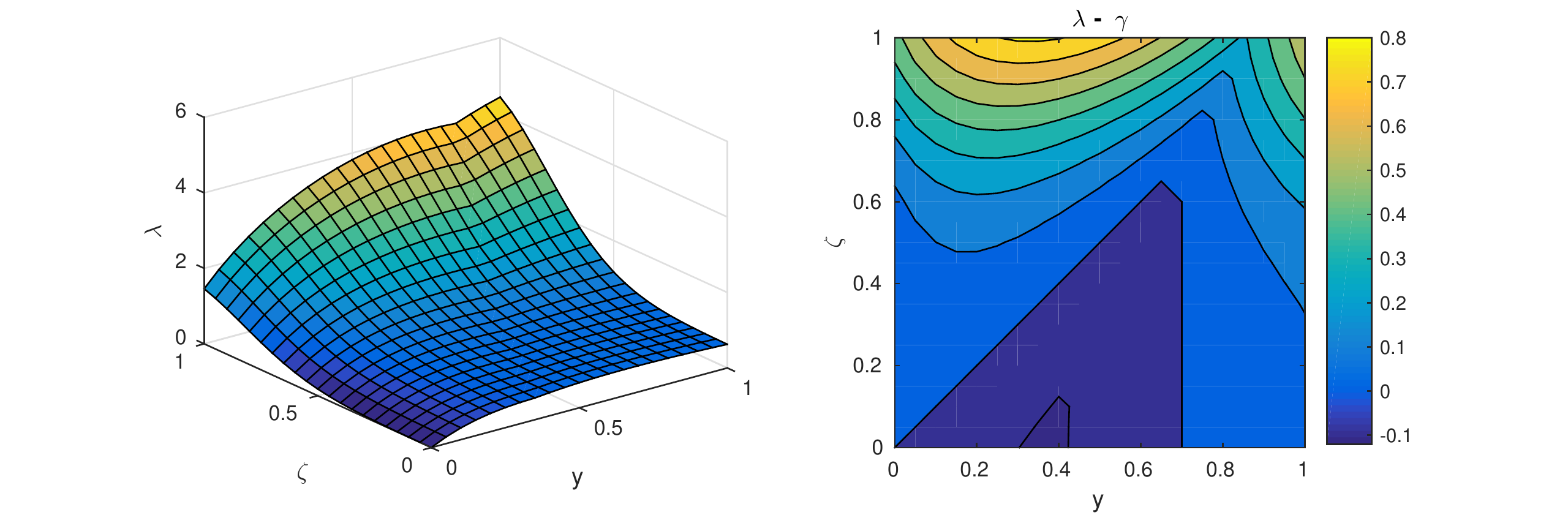}
	\caption{The PMT rate $\lambda$ and its relation to the MPT rate $\gamma$ with respect to the
		amount of occupied cell-cell receptors ($\zeta$) and cell-tissue receptors ($y$) in case
		$\gamma_0=1$.}\label{fig:rates}
\end{figure}

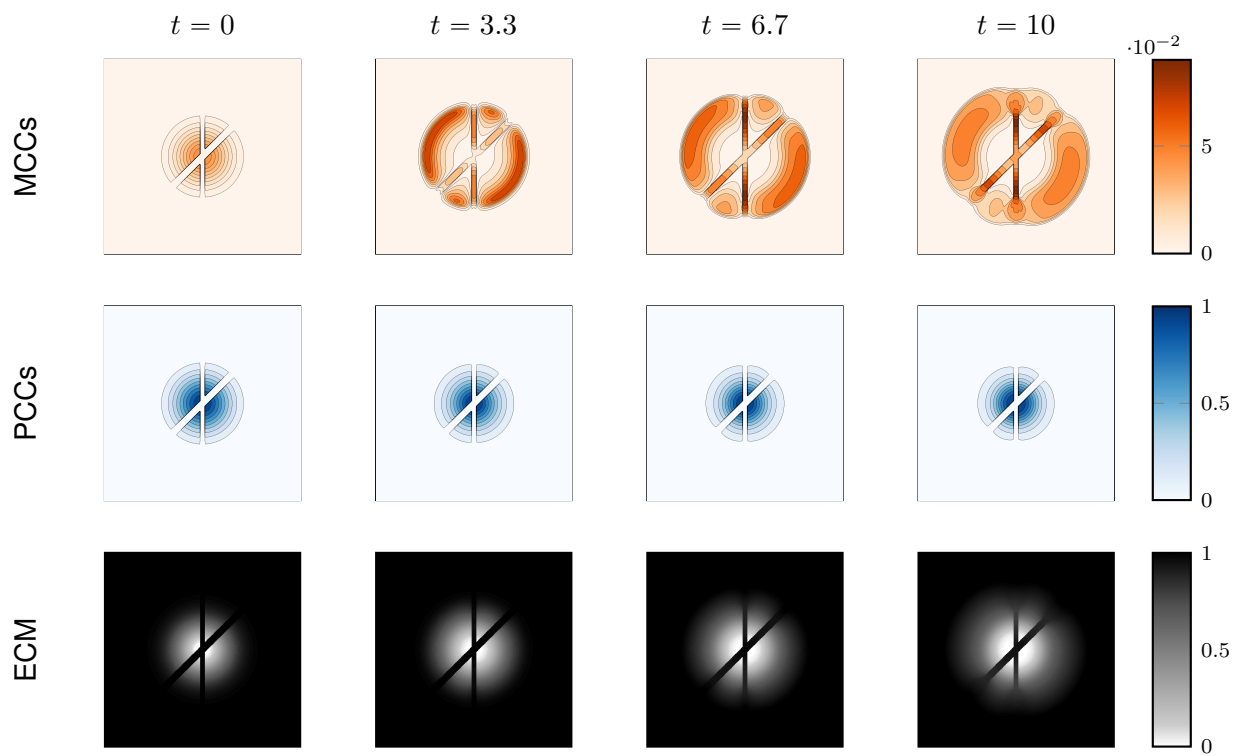
\begin{figure}[t] % exp 2 - stripes
	\figurewidth=\linewidth
	\def \expSetup {TODO1B_ECM_stripes}
	\begin{tikzpicture}
	\begin{groupplot}[
	/tikz/mark size=1.5pt,
	group style={
		group name=my plots,
		group size=4 by 3,
		horizontal sep=1cm,      % <-- default: 1cm
		vertical sep=0.7cm,        % <-- default: 1cm
	},
	xmin=-2,
	xmax=2,
	ymin=-2,
	ymax=2,
	ylabel shift = 2 em,
	xtick = \empty,
	ytick= \empty,
	ticklabel style = {font=\scriptsize},
	axis line style = thick,
	axis background/.style={fill=white},
	width=.25\figurewidth,
	height=.25\figurewidth, 
	]
	\nextgroupplot[ylabel = MCCs, title = {$t=0$}]
	\addplot [forget plot] graphics [xmin=-2, xmax=2, ymin=-2, ymax=2] {\expSetup1};
	\nextgroupplot[title = {$t=3.3$}]
	\addplot [forget plot] graphics [xmin=-2, xmax=2, ymin=-2, ymax=2] {\expSetup2};
	\nextgroupplot[title = {$t=6.7$}]
	\addplot [forget plot] graphics [xmin=-2, xmax=2, ymin=-2, ymax=2] {\expSetup3};
	\nextgroupplot[
	title = {$t=10$},
	point meta min=0,
	point meta max=0.09,
	colormap name=MCCs,
	colorbar,
	colorbar style={at={(1.2,1)},anchor=north west}
	]
	\addplot [forget plot] graphics [xmin=-2, xmax=2, ymin=-2, ymax=2] {\expSetup4};
	\nextgroupplot[ylabel = PCCs]
	\addplot [forget plot] graphics [xmin=-2, xmax=2, ymin=-2, ymax=2] {\expSetup5};
	\nextgroupplot
	\addplot [forget plot] graphics [xmin=-2, xmax=2, ymin=-2, ymax=2] {\expSetup6};
	\nextgroupplot
	\addplot [forget plot] graphics [xmin=-2, xmax=2, ymin=-2, ymax=2] {\expSetup7};
	\nextgroupplot[
	point meta min=0,
	point meta max=1,
	colormap name=PCCs,
	colorbar,
	colorbar style={at={(1.2,1)},anchor=north west}
	]
	\addplot [forget plot] graphics [xmin=-2, xmax=2, ymin=-2, ymax=2] {\expSetup8};
	\nextgroupplot[ylabel = ECM]
	\addplot [forget plot] graphics [xmin=-2, xmax=2, ymin=-2, ymax=2] {\expSetup9};
	\nextgroupplot
	\addplot [forget plot] graphics [xmin=-2, xmax=2, ymin=-2, ymax=2] {\expSetup10};
	\nextgroupplot
	\addplot [forget plot] graphics [xmin=-2, xmax=2, ymin=-2, ymax=2] {\expSetup11};
	\nextgroupplot[
	point meta min=0,
	point meta max=1,
	colormap name=ECM,
	colorbar,
	colorbar style={at={(1.2,1)},anchor=north west}
	]
	\addplot [forget plot] graphics [xmin=-2, xmax=2, ymin=-2, ymax=2] {\expSetup12};
	\end{groupplot}
	\end{tikzpicture}%
	
	\caption{Simulation results of \textbf{Experiment~\ref{exp:dynEMT} --- dynamic phenotypic switch
		rates} using the ECM-with-stripes initial conditions \eqref{eq:ICstrps}. When comparing with the
	corresponding simulation in Experiment~\ref{exp:cnstEMT} (constant phenotypic transition rates) shown in
	Figure \ref{fig:exp1_stripes}, the MCCs can infer higher densities at sites with larger ECM gradients, and lower ones where cell-tissue interfaces are less sharp (e.g., at the center of fiber strands crossing), thus allowing for less cells to move beyond the main fiber tracts.   %exhibit similar behaviour 
	%and are somewhat more uniform 	especially in regions where the ECM is more depleted.
}\label{fig:exp2_stripes}	
\end{figure}

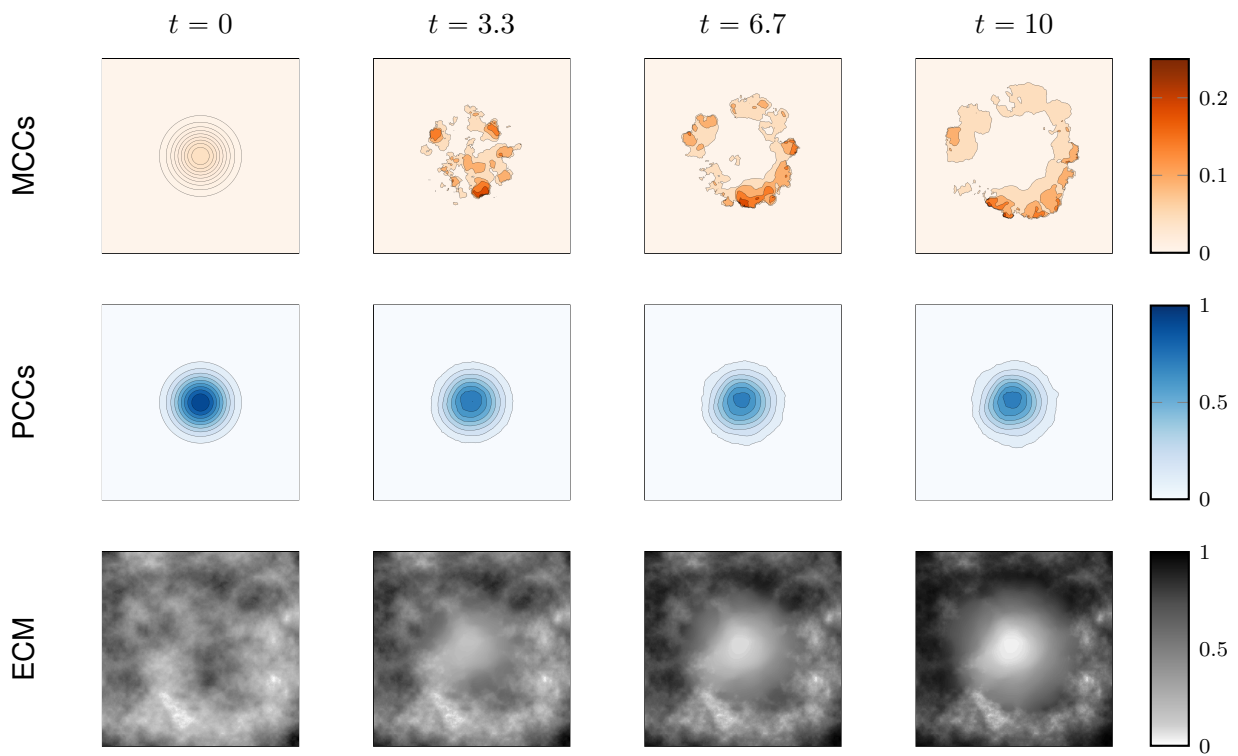
\begin{figure}[t] % exp 2 - rnd
	\figurewidth=\linewidth
	\def \expSetup {TODO1B_ECM_rnd}
	\begin{tikzpicture}
	\begin{groupplot}[
	/tikz/mark size=1.5pt,
	group style={
		group name=my plots,
		group size=4 by 3,
		horizontal sep=1cm,      % <-- default: 1cm
		vertical sep=0.7cm,        % <-- default: 1cm
	},
	xmin=-2,
	xmax=2,
	ymin=-2,
	ymax=2,
	ylabel shift = 2 em,
	xtick = \empty,
	ytick= \empty,
	ticklabel style = {font=\scriptsize},
	axis line style = thick,
	axis background/.style={fill=white},
	width=.25\figurewidth,
	height=.25\figurewidth, 
	]
	\nextgroupplot[ylabel = MCCs, title = {$t=0$}]
	\addplot [forget plot] graphics [xmin=-2, xmax=2, ymin=-2, ymax=2] {\expSetup1};
	\nextgroupplot[title = {$t=3.3$}]
	\addplot [forget plot] graphics [xmin=-2, xmax=2, ymin=-2, ymax=2] {\expSetup2};
	\nextgroupplot[title = {$t=6.7$}]
	\addplot [forget plot] graphics [xmin=-2, xmax=2, ymin=-2, ymax=2] {\expSetup3};
	\nextgroupplot[
	title = {$t=10$},
	point meta min=0,
	point meta max=0.25,
	colormap name=MCCs,
	colorbar,
	colorbar style={at={(1.2,1)},anchor=north west}
	]
	\addplot [forget plot] graphics [xmin=-2, xmax=2, ymin=-2, ymax=2] {\expSetup4};
	\nextgroupplot[ylabel = PCCs]
	\addplot [forget plot] graphics [xmin=-2, xmax=2, ymin=-2, ymax=2] {\expSetup5};
	\nextgroupplot
	\addplot [forget plot] graphics [xmin=-2, xmax=2, ymin=-2, ymax=2] {\expSetup6};
	\nextgroupplot
	\addplot [forget plot] graphics [xmin=-2, xmax=2, ymin=-2, ymax=2] {\expSetup7};
	\nextgroupplot[
	point meta min=0,
	point meta max=1,
	colormap name=PCCs,
	colorbar,
	colorbar style={at={(1.2,1)},anchor=north west}
	]
	\addplot [forget plot] graphics [xmin=-2, xmax=2, ymin=-2, ymax=2] {\expSetup8};
	\nextgroupplot[ylabel = ECM]
	\addplot [forget plot] graphics [xmin=-2, xmax=2, ymin=-2, ymax=2] {\expSetup9};
	\nextgroupplot
	\addplot [forget plot] graphics [xmin=-2, xmax=2, ymin=-2, ymax=2] {\expSetup10};
	\nextgroupplot
	\addplot [forget plot] graphics [xmin=-2, xmax=2, ymin=-2, ymax=2] {\expSetup11};
	\nextgroupplot[
	point meta min=0,
	point meta max=1,
	colormap name=ECM,
	colorbar,
	colorbar style={at={(1.2,1)},anchor=north west}
	]
	\addplot [forget plot] graphics [xmin=-2, xmax=2, ymin=-2, ymax=2] {\expSetup12};
	\end{groupplot}
	\end{tikzpicture}%
	
	\caption{Simulation results of \textbf{Experiment~\ref{exp:dynEMT} --- dynamic phenotypic switch
		rates} with the randomly-structured  ECM initial conditions \eqref{eq:ICrnd}. Comparing with the
	simulation of Experiment~\ref{exp:cnstEMT} (constant phenotypic transition), shown in
	Figure~\ref{fig:exp1_rnd}, the MCCs' invasion is here slightly  more cohesive, allowing for fewer, but larger local maxima. The density of the PCCs is thereby slightly
	lower than in Figure~\ref{fig:exp1_rnd}.}\label{fig:exp2_rnd}	
\end{figure}

\begin{experiment}\textbf{Dynamic phenotypic switch rates.}\label{exp:dynEMT}
In the next step we consider dynamic PMT and MPT rates $\lambda$ and $\gamma$ that depend on the
amount of cell surface receptors bound to insoluble ligands in the tumor microenvironment. 

We choose the MPT rate as a sum between a Gamma-distribution with respect to the amount of receptors
bound to tissue and a bivariate normal distribution with respect to the variables accounting (via
receptors) for cell-cell and cell-tissue interactions on the cellular level:
\begin{equation}\label{eq:sim-gamma}
\gamma(y, \zeta) = \frac{\gamma_0 \, b^p}{\Gamma(p)}y^{p-1}e^{-by} + \frac{\gamma_0}{2 \pi \sigma_y
	\sigma_\zeta}	
e^{-\frac{(y-\mu_y)^2}{2 \sigma_y^2} - \frac{(\zeta - \mu_\zeta)^2}{2\sigma_\zeta^2}},
\end{equation}
where we take $\gamma_0=0.1$, $b=2$, $p=2$, $\mu_y=\mu_\zeta=2$, $\sigma_y = 0.5$ and
$\sigma_\zeta=0.3.$ This choice ensures that $\gamma $ is largest when many cell receptors are
occupied, thus limiting motility and correspondingly enhancing proliferation, due to the dichotomous
behavior of the two phenotypes. Thereby, cell-tissue interactions play a more prominent role than
intercellular attachment, as cells can also migrate as collectives.The PMT rate is modeled on the
other hand by
\begin{equation}\label{eq:sim-lambda}
\lambda(y, \zeta) = \frac{2 \gamma}{1 + e^{-m}}, \quad m = \max\{\zeta - y,~y -y_\text{ref}\}.
\end{equation} 
Here we assume that $y_\text{ref}=0.6$. In this model the PMT rate exceeds the MPT rate whenever
there are more cell-cell than cell-tissue receptors occupied or when the amount of cell-tissue interactions 
exceeds a reference amount $y_\text{ref}$ of correspondingly bound receptors. If on the other hand the amount of occupied cell-tissue
receptors is small, but larger than that of cell-cell receptor occupancy, the transition to the
proliferating type is dominant. In the case of no cell-cell interactions, both transition
rates are close to zero. These correlations are illustrated in Figure~\ref{fig:rates}. 

As in the previous Experiment~\ref{exp:cnstEMT}, we consider both sets of initial conditions
\eqref{eq:ICstrps} and \eqref{eq:ICrnd}. The corresponding simulation results are shown in
Figures \ref{fig:exp2_stripes} and \ref{fig:exp2_rnd}, respectively. The comparison with Experiment~\ref{exp:cnstEMT} reveals that the dynamic behavior of the MCC invasion and the proliferation of the PCCs is similar, however the current regime can lead to higher densities of MCCs at sites with higher tissue gradients, while lowering the cell density at less sharp cell-tissue interfaces (as e.g., at the crossing of tissue strands in Figures \ref{fig:exp1_stripes} and \ref{fig:exp2_stripes}). Thus, dynamic PMT and MPT rates accentuate the effect of cell-tissue and cell-cell interactions, as they contribute to the source terms of migrating/proliferating cells and hence to deciding the fate of cells. The accumulation of cells at either side of the mentioned crossing of tissue strands is no longer visible in this scenario were the phenotypic switch depends on the local availability of cells and tissue. The influence of dynamic transition rates on the tactic behavior is more indirect, as the sensitivity functions $\chi_1$ and $\chi_2$ involve the sum $m+p$ rather than one of the two cell densities, while the two phenotypes also act on $v$ in a joint manner. Figure \ref{fig:exp2_rnd} shows in a more pregnant way a similar trend: the cells have a higher cohesion tendency, which leads to a correspondingly localized tissue degradation. %{\cred There is one issue with the initial ECM picture: it looks differently in Figures \ref{fig:exp1_rnd} and \ref{fig:exp2_rnd}, so we cannot really compare them, can we? Same applies to all subsequent pictures involving the random ECM: the IC should be the same. }
%less ``aggressive'' Experiment~\ref{exp:dynEMT} despite the similar dynamical behaviour of the ECM.
\end{experiment}

% -------------- exp3
\begin{experiment}\textbf{Acidity driven migration.}\label{exp:acidity}
	In this experiment we replace the repellence between the PCCs and MCCs by repellent pH-taxis. This is more indirect: both cell phenotypes produce acidity, which they try to avoid, as it is detrimental to their functioning; in particular to their proliferation. The interspecies repellence is thus replaced by chemotaxis away from high acidity (i.e. low pH). The haptotaxis is indirectly affected: the tissue degradation is now due to acidity instead of direct contacts with cancer cells. To describe these effects we introduce a variable $h$ representing acidity in terms of proton concentration and
	consider the modified model
	\begin{subequations}\label{eq:modelacidity}
		\begin{align}
		\pd_t m =& \nabla\cdot\( D_c\frac{1+ m\,p + m\,v + p\,v}{1+m\,(p+v)}\nabla m \) 
		-\nabla \cdot\( \chi_1(m,p,v) m\nabla v \) 
		+\nabla \cdot\( \chi_2(m,p,v)  m\nabla h \)\notag\\
		&+\lambda(y,\zeta)\,p -\gamma(y,\zeta)\,m\\
		\pd_t p =& {\mu(h)} \, p\,(1-(m+p)-v) -\lambda (y,\zeta)\,p +\gamma(y,\zeta)\,m\\
		\pd_t v =& -\delta\, h\,v+\mu_v\,v\,(1-(m+p)-v) \label{eq:ac:v}\\
		\pd_t h = & D_h \Delta h + \alpha(m+ p) - \beta h. \label{eq:ac:h}
		\end{align}
	\end{subequations}
	We choose the PMT and MPT phenotypic transition rates $\lambda$ and $\gamma$ as in
	Experiment~\ref{exp:dynEMT}, pH-sensitivity $\chi_2(m,p,v)$ from \eqref{eq:microSteady2} and use the
	parameters $D_h=0.07$, $\alpha=0.55$, $\beta=0.05$ and $\delta=0.2$ in equations \eqref{eq:ac:v} and
	\eqref{eq:ac:h}. Here, again, we were guided in the parameter assessment by the choice made in \cite{ChLo06} for the chemotactic  signal \footnote{there it was uPA protease, which should be comparable with acidity; the production of proteases and acidity by tumor cells are known to be tightly related}. As for all other model components we prescribe homogeneous Neumann boundary conditions for the acidity $h$, moreover we set its initial concentration proportional to the initial PCC density, by $h_0(x)= 0.2 + p_0(x)$ and consider acidity-dependent cell proliferation in the form $\mu(h)=\mu_0(h-1)_+$, where we set $\mu_0=0.1$. Thereby, $h$ is a nondimensional quantity, thus the proliferation rate $\mu $ actually refers to a comparison of the rescaled acidity $\tilde h:=hh_T$ with the threshold proton concentration $h_T=10^{-6.4}$ (which when exceeded is known to be lethal even for tumor cells \cite{Webb2011}). In the plots we visualize the acidity %we rescale the nondimensional quantity $h$ in \eqref{eq:ac:h} with respect to the threshold proton concentration $h_T=10^{-6.4}$ (which when exceeded is known to be lethal even for tumor cells \cite{Webb2011}) and consider acidity 
	in terms of pH levels, i.e. $-\log_{10}(hh_T)$.
	
%	{\cred nkls: I kept Christinas remarks for reference.
%        
%        Perhaps we would see more if we include the detrimental effect on proliferation, i.e. replace the $p$-equation with
%$$\pd_t p =\mu (h)\, p\,(1-(m+p)-v) -\lambda (y,\zeta)\,p +\gamma(y,\zeta)\,m$$	
%where e.g. $\mu (h)=\mu_0(h-h_T)_+$, with $h_T$ some critical acidity threshold for tumor cells, thus we could take $h_T=10^{-6.4}$, and $\mu_0$ as previously for $\mu$. If that $\mu _0$ is too large or too small we can try something between $0.05-2$. If we do not acount for $h$-dependent proliferation we don't actually change anything in the $m$- and $p$-equations -except the indirect effect of tissue degradation and pH-taxis. 
%
%What should actually be done when modeling more precisely is to consider yet another subcellular variable $\xi $ accounting for transmembrane units occupied with protons and phenotypic switch rates depending on $\xi $ as well. We don't do that, however - I will only mention the possibility, here or in the Discussion section. 
%
%We need the pH to be between 6.4 (very acidic) and 7.2 (alkaline), with corresponding h-values.}
	
	The numerical simulation results corresponding to the initial conditions \eqref{eq:ICstrps}
	and\eqref{eq:ICrnd} are presented in Figures \ref{fig:exp3_stripes} and \ref{fig:exp3_rnd}.  Comparing these simulations with the corresponding ones in Experiment~\ref{exp:dynEMT} we see that the MCCs invade a wider area and form smaller aggregates, while the proliferation of tumor cells is reduced, simultaneously with a stronger tissue degradation  throughout the domain.
\end{experiment}

% -------------- exp4
\begin{experiment}\textbf{Degenerate diffusion.}\label{exp:degenerate}
	In this experiment we consider a model variant in which the diffusion coefficient in
	\eqref{eq:diffusion}, 
	$$D_c\frac{1+mp+mv+pv}{1+mv+pv},$$
	is adjusted to the degenerate version
	\begin{equation}\label{diff-coeff-degen}
	D_c\frac{mv+mp+pv}{1+mv+pv},
	\end{equation}
	which allows, among others, for very small diffusivities at low MCC densities. Moreover, as mentioned in Section \ref{sec:model}, this quantity nullifies whenever two of the model variables become zero simultaneously.  Apart from this modification, we employ the same PMT and MPT rates and parameters as 
	in Experiment~\ref{exp:dynEMT}. The numerical simulation results corresponding to the initial
	conditions \eqref{eq:ICstrps} and \eqref{eq:ICrnd} are presented in Figures \ref{fig:exp4_stripes}
	and \ref{fig:exp4_rnd}, respectively. These simulations are compared with the corresponding ones in
	Experiment~\ref{exp:dynEMT}. We note that the present case leads to MCCs forming very localized, relatively high aggregates which either merge into regions of high density or remain localized (and large), while the PCCs are almost not affected. Such behavior is better visible in Figure~\ref{fig:exp4_stripes} and in the closeup in Figure~\ref{fig:exp4_stripes_closeup}; one can observe small ``islands'' of MCCs which are ``trapped'' there by the fact that one of $m$, $p$, and $v$ is (almost) completely degraded, while another one has a relatively low density. A similar behavior can be seen when comparing Figures \ref{fig:exp2_rnd} and \ref{fig:exp4_rnd}: the degenerate case allows for higher, more localized MCC densities, mainly near the invasion front.
		
	%We also note a higher 	fragmentation in the support of the MCCs in the degenerate diffusion case. On the other hand the 	evolution of the PCCs and the ECM is almost identical between the two experiments.
\end{experiment}

\begin{figure}[t] % exp5 - stripes
	\figurewidth=\linewidth
	\def \expSetup {TODO4_ECM_stripes}
	\begin{tikzpicture}
	\begin{groupplot}[
	/tikz/mark size=1.5pt,
	group style={
		group name=my plots,
		group size=4 by 3,
		horizontal sep=1cm,      % <-- default: 1cm
		vertical sep=0.7cm,        % <-- default: 1cm
	},
	xmin=-2,
	xmax=2,
	ymin=-2,
	ymax=2,
	ylabel shift = 2 em,
	xtick = \empty,
	ytick= \empty,
	ticklabel style = {font=\scriptsize},
	axis line style = thick,
	axis background/.style={fill=white},
	width=.25\figurewidth,
	height=.25\figurewidth, 
	]
	\nextgroupplot[ylabel = MCCs, title = {$t=0$}]
	\addplot [forget plot] graphics [xmin=-2, xmax=2, ymin=-2, ymax=2] {\expSetup1};
	\nextgroupplot[title = {$t=3.3$}]
	\addplot [forget plot] graphics [xmin=-2, xmax=2, ymin=-2, ymax=2] {\expSetup2};
	\nextgroupplot[title = {$t=6.7$}]
	\addplot [forget plot] graphics [xmin=-2, xmax=2, ymin=-2, ymax=2] {\expSetup3};
	\nextgroupplot[
	title = {$t=10$},
	point meta min=0,
	point meta max=0.09,
	colormap name=MCCs,
	colorbar,
	colorbar style={at={(1.2,1)},anchor=north west}
	]
	\addplot [forget plot] graphics [xmin=-2, xmax=2, ymin=-2, ymax=2] {\expSetup4};
	\nextgroupplot[ylabel = PCCs]
	\addplot [forget plot] graphics [xmin=-2, xmax=2, ymin=-2, ymax=2] {\expSetup5};
	\nextgroupplot
	\addplot [forget plot] graphics [xmin=-2, xmax=2, ymin=-2, ymax=2] {\expSetup6};
	\nextgroupplot
	\addplot [forget plot] graphics [xmin=-2, xmax=2, ymin=-2, ymax=2] {\expSetup7};
	\nextgroupplot[
	point meta min=0,
	point meta max=1,
	colormap name=PCCs,
	colorbar,
	colorbar style={at={(1.2,1)},anchor=north west}
	]
	\addplot [forget plot] graphics [xmin=-2, xmax=2, ymin=-2, ymax=2] {\expSetup8};
	\nextgroupplot[ylabel = ECM]
	\addplot [forget plot] graphics [xmin=-2, xmax=2, ymin=-2, ymax=2] {\expSetup9};
	\nextgroupplot
	\addplot [forget plot] graphics [xmin=-2, xmax=2, ymin=-2, ymax=2] {\expSetup10};
	\nextgroupplot
	\addplot [forget plot] graphics [xmin=-2, xmax=2, ymin=-2, ymax=2] {\expSetup11};
	\nextgroupplot[
	point meta min=0,
	point meta max=1,
	colormap name=ECM,
	colorbar,
	colorbar style={at={(1.2,1)},anchor=north west}
	]
	\addplot [forget plot] graphics [xmin=-2, xmax=2, ymin=-2, ymax=2] {\expSetup12};
	\end{groupplot}
	
	\end{tikzpicture}%
	
	\caption{Simulation results of \textbf{Experiment~\ref{exp:remod} --- ECM remodeling by cancer
			cells} with the ECM-with-stripes initial conditions \eqref{eq:ICstrps}. Compared to
		Experiment~\ref{exp:dynEMT} (dynamic phenotypic transition with self-remodeling of the matrix)
		shown in Figure~\ref{fig:exp2_stripes} we see only a slight impact of the cell reconstruction of tissue; the results are almost identical, maybe with a slightly higher
		concentration of the MCCs towards the invasion front and higher ECM degradation in the inner part of the tumor.}\label{fig:exp5_stripes}	
\end{figure}
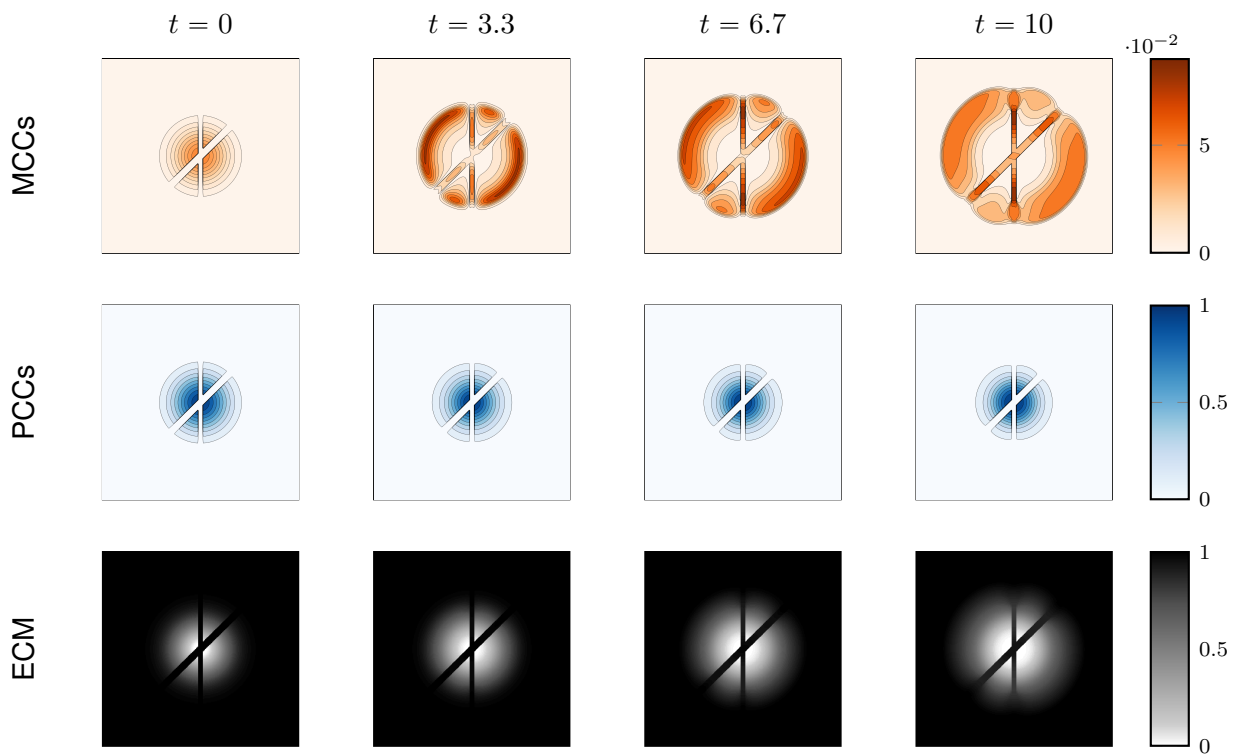

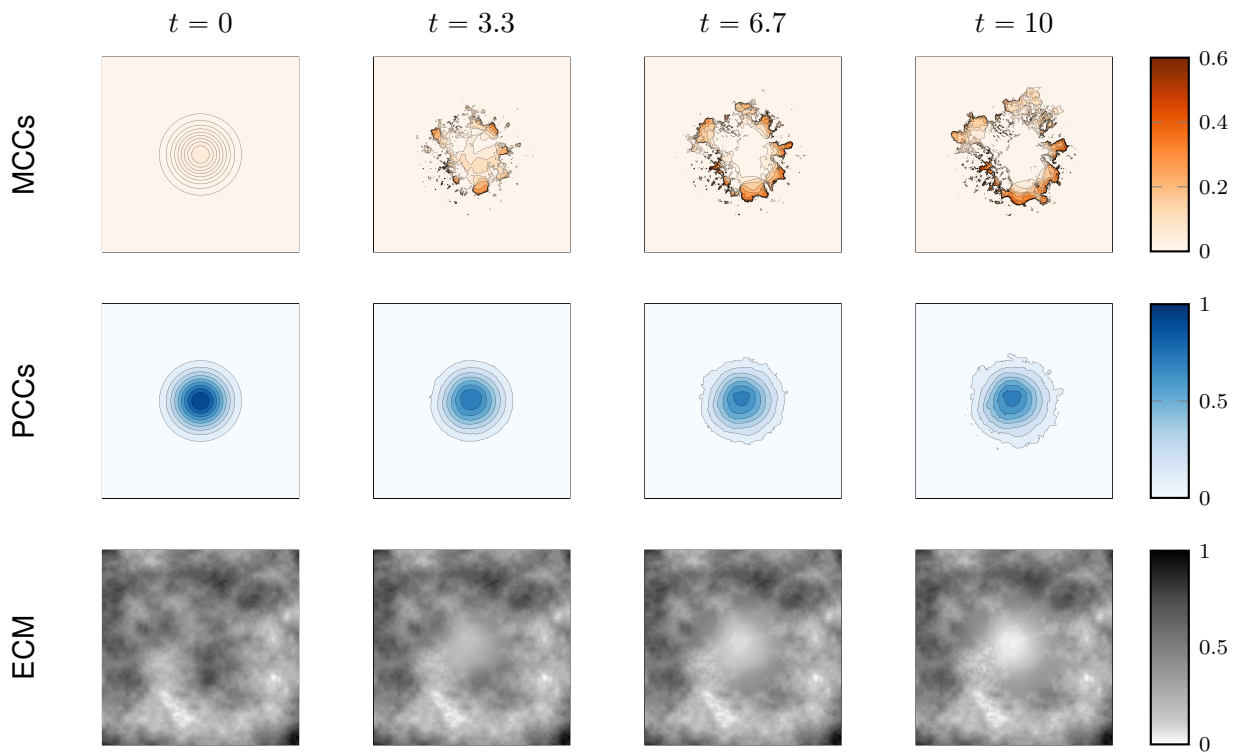
\begin{figure}[t] %exp5 - rnd 
	\figurewidth=\linewidth
	\def \expSetup {TODO4_ECM_rnd}
	\begin{tikzpicture}
	\begin{groupplot}[
	/tikz/mark size=1.5pt,
	group style={
		group name=my plots,
		group size=4 by 3,
		horizontal sep=1cm,      % <-- default: 1cm
		vertical sep=0.7cm,        % <-- default: 1cm
	},
	xmin=-2,
	xmax=2,
	ymin=-2,
	ymax=2,
	ylabel shift = 2 em,
	xtick = \empty,
	ytick= \empty,
	ticklabel style = {font=\scriptsize},
	axis line style = thick,
	axis background/.style={fill=white},
	width=.25\figurewidth,
	height=.25\figurewidth, 
	]
	\nextgroupplot[ylabel = MCCs, title = {$t=0$}]
	\addplot [forget plot] graphics [xmin=-2, xmax=2, ymin=-2, ymax=2] {\expSetup1};
	\nextgroupplot[title = {$t=3.3$}]
	\addplot [forget plot] graphics [xmin=-2, xmax=2, ymin=-2, ymax=2] {\expSetup2};
	\nextgroupplot[title = {$t=6.7$}]
	\addplot [forget plot] graphics [xmin=-2, xmax=2, ymin=-2, ymax=2] {\expSetup3};
	\nextgroupplot[
	title = {$t=10$},
	point meta min=0,
	point meta max=0.6,
	colormap name=MCCs,
	colorbar,
	colorbar style={at={(1.2,1)},anchor=north west}
	]
	\addplot [forget plot] graphics [xmin=-2, xmax=2, ymin=-2, ymax=2] {\expSetup4};
	\nextgroupplot[ylabel = PCCs]
	\addplot [forget plot] graphics [xmin=-2, xmax=2, ymin=-2, ymax=2] {\expSetup5};
	\nextgroupplot
	\addplot [forget plot] graphics [xmin=-2, xmax=2, ymin=-2, ymax=2] {\expSetup6};
	\nextgroupplot
	\addplot [forget plot] graphics [xmin=-2, xmax=2, ymin=-2, ymax=2] {\expSetup7};
	\nextgroupplot[
	point meta min=0,
	point meta max=1,
	colormap name=PCCs,
	colorbar,
	colorbar style={at={(1.2,1)},anchor=north west}
	]
	\addplot [forget plot] graphics [xmin=-2, xmax=2, ymin=-2, ymax=2] {\expSetup8};
	\nextgroupplot[ylabel = ECM]
	\addplot [forget plot] graphics [xmin=-2, xmax=2, ymin=-2, ymax=2] {\expSetup9};
	\nextgroupplot
	\addplot [forget plot] graphics [xmin=-2, xmax=2, ymin=-2, ymax=2] {\expSetup10};
	\nextgroupplot
	\addplot [forget plot] graphics [xmin=-2, xmax=2, ymin=-2, ymax=2] {\expSetup11};
	\nextgroupplot[
	point meta min=0,
	point meta max=1,
	colormap name=ECM,
	colorbar,
	colorbar style={at={(1.2,1)},anchor=north west}
	]
	\addplot [forget plot] graphics [xmin=-2, xmax=2, ymin=-2, ymax=2] {\expSetup12};
	\end{groupplot}
	\end{tikzpicture}%
	
	\caption{Simulation results of \textbf{Experiment~\ref{exp:remod} --- ECM remodeling by cancer
			cells} on a randomly-structured ECM \eqref{eq:ICrnd}.  When compared with Experiment~\ref{exp:dynEMT}
		(dynamic phenotypic transition with self-remodeling of the matrix), shown in
		Figure~\ref{fig:exp2_rnd}, it is clear that the cell reconstruction of the ECM leads to a  more fragmented invasion of the MCCs invasion and to higher concentrations along the propagating fronts. We moreover see that the PCCs exhibit a non-smooth boundary/periphery in their support, and that the reconstruction of the ECM is localized where the MCCs are located.}\label{fig:exp5_rnd}	
\end{figure}

\begin{experiment}\textbf{ECM remodeling by cancer cells.}\label{exp:remod}
	In this experiment we consider a model variant in which the MCCs are responsible for the remodeling
	of the ECM rather than the tissue itself. Indeed, there is abundant biological evidence that tumor cells are able to restructure the ECM in a manner favorable to their migration, see e.g. \cite{Malandrino2018,Xiong2016} and references therein. We therefore adjust \eqref{eq:model-c} to
	\begin{equation}
	\pd_t v = -\delta\,(m+p)\,v+\mu_v\,m\,(1-(m+p)-v).
	\end{equation}
	Except for the matrix reconstruction rate, which is adapted to $\mu_v = 0.5$, we consider the same
	parameters as in Experiment~\ref{exp:dynEMT}. The numerical simulation results are presented in
	Figures \ref{fig:exp5_stripes} and \ref{fig:exp5_rnd}. 	
	Comparing with the ECM-with-stripes initial conditions of Experiment~\ref{exp:dynEMT} (dynamic
	phenotypic switch with self-remodeling of the matrix), shown in Figure~\ref{fig:exp2_stripes}, the
	results are quite similar - with some higher MCC densities in the former case. Compared to the randomly-structured ECM initial conditions, shown in
	Figure~\ref{fig:exp2_rnd}, however, we see the effect of cancer cell remodeling of the ECM through a more fragmented support of the MCCs and higher concentrations closer to the propagating fronts.
	Similarly, the domain occupied by PCCs exhibits more fractal boundaries. Overall, the tumor will have a more infiltrative spread, with less sharp margins; such behavior is typical for glioblastoma multiforme, see e.g., \cite{Matsukado1961}.
\end{experiment}

\begin{figure}[t] % exp6 - rnd
	\figurewidth=\linewidth
	\def \expSetup {TODO5_ECM_rnd}
	\begin{tikzpicture}
	\begin{groupplot}[
	/tikz/mark size=1.5pt,
	group style={
		group name=my plots,
		group size=4 by 3,
		horizontal sep=1cm,      % <-- default: 1cm
		vertical sep=0.7cm,        % <-- default: 1cm
	},
	xmin=-2,
	xmax=2,
	ymin=-2,
	ymax=2,
	ylabel shift = 2 em,
	xtick = \empty,
	ytick= \empty,
	ticklabel style = {font=\scriptsize},
	axis line style = thick,
	axis background/.style={fill=white},
	width=.25\figurewidth,
	height=.25\figurewidth, 
	]
	\nextgroupplot[ylabel = MCCs, title = {$t=0$}]
	\addplot [forget plot] graphics [xmin=-2, xmax=2, ymin=-2, ymax=2] {\expSetup1};
	\nextgroupplot[title = {$t=3.3$}]
	\addplot [forget plot] graphics [xmin=-2, xmax=2, ymin=-2, ymax=2] {\expSetup2};
	\nextgroupplot[title = {$t=6.7$}]
	\addplot [forget plot] graphics [xmin=-2, xmax=2, ymin=-2, ymax=2] {\expSetup3};
	\nextgroupplot[
	title = {$t=10$},
	point meta min=0,
	point meta max=0.18,
	colormap name=MCCs,
	colorbar,
	colorbar style={at={(1.2,1)},anchor=north west}
	]
	\addplot [forget plot] graphics [xmin=-2, xmax=2, ymin=-2, ymax=2] {\expSetup4};
	\nextgroupplot[ylabel = PCCs]
	\addplot [forget plot] graphics [xmin=-2, xmax=2, ymin=-2, ymax=2] {\expSetup5};
	\nextgroupplot
	\addplot [forget plot] graphics [xmin=-2, xmax=2, ymin=-2, ymax=2] {\expSetup6};
	\nextgroupplot
	\addplot [forget plot] graphics [xmin=-2, xmax=2, ymin=-2, ymax=2] {\expSetup7};
	\nextgroupplot[
	point meta min=0,
	point meta max=1,
	colormap name=PCCs,
	colorbar,
	colorbar style={at={(1.2,1)},anchor=north west}
	]
	\addplot [forget plot] graphics [xmin=-2, xmax=2, ymin=-2, ymax=2] {\expSetup8};
	\nextgroupplot[ylabel = ECM]
	\addplot [forget plot] graphics [xmin=-2, xmax=2, ymin=-2, ymax=2] {\expSetup9};
	\nextgroupplot
	\addplot [forget plot] graphics [xmin=-2, xmax=2, ymin=-2, ymax=2] {\expSetup10};
	\nextgroupplot
	\addplot [forget plot] graphics [xmin=-2, xmax=2, ymin=-2, ymax=2] {\expSetup11};
	\nextgroupplot[
	point meta min=0,
	point meta max=1,
	colormap name=ECM,
	colorbar,
	colorbar style={at={(1.2,1)},anchor=north west}
	]
	\addplot [forget plot] graphics [xmin=-2, xmax=2, ymin=-2, ymax=2] {\expSetup12};
	\end{groupplot}
	\end{tikzpicture}%
	
	\caption{Simulation results of \textbf{Experiment~\ref{exp:anoikis} --- anoikis effect} on the
	randomly structured  ECM initial conditions \eqref{eq:ICrnd}. Compared to 
	Experiment~\ref{exp:dynEMT}, shown in Figure~\ref{fig:exp2_rnd}, the effect of anoikis becomes visible: the tumor pattern is more heterogeneous (mainly due to the evolution of MCCs) with correspondingly lower PCC 
	density in regions with stronger degraded ECM.}\label{fig:exp6_rnd}	
\end{figure}
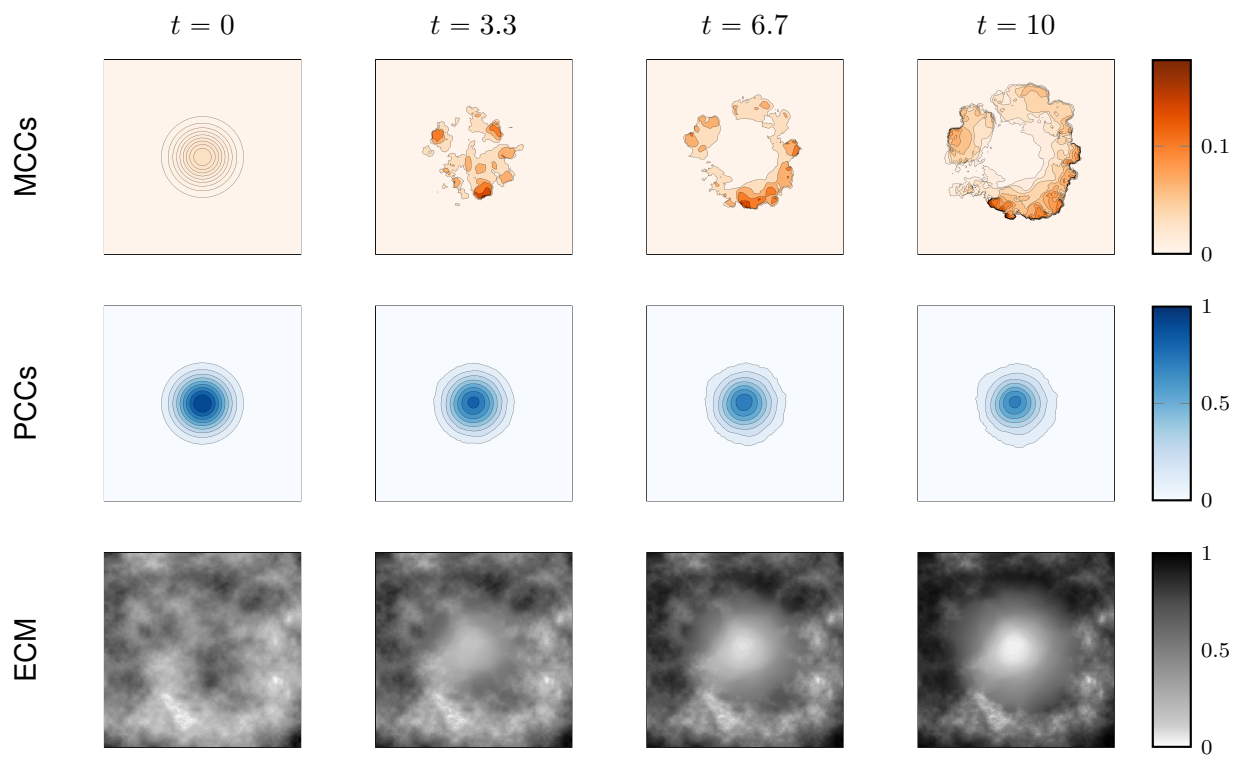

\begin{experiment}\textbf{Anoikis effect.}\label{exp:anoikis} 
	Cell--matrix interactions are essential for cell survival. When detaching from the
	ECM, the cell cycle is arrested and a particular form of programmed cell death, known as
	\textit{anoikis}, is initiated. Cancer cells can, however, escape anoikis in the invasion process
	\cite{anoikis.1994, HAWEI11}. 	
	In this last experiment we account for this aspect by replacing \eqref{eq:model-b} with
	\begin{equation}
	\pd_t p = \mu pv (1-(m+p)-v) -\lambda (y,\zeta)\,p +\gamma(y,\zeta)\,m.
	\end{equation}
	In this approach the first term ensures proliferation only when the cells are in contact with the
	ECM. %\comm{I hope I got this right, one could also deduce that this is a model for anoikis because
	%there is only proliferation when there is ECM}
	The parameters chosen here are the same as in Experiment~\ref{exp:dynEMT} with the corresponding
	simulation results for the two types of initial conditions \eqref{eq:ICstrps} and \eqref{eq:ICrnd}
	being presented in Figures \ref{fig:exp6_stripes} and \ref{fig:exp6_rnd}, respectively. Comparing with the
	simulations of Experiment~\ref{exp:dynEMT}, we note that when the ECM-with-stripes initial
	condition \eqref{eq:ICstrps} is considered the results are almost identical (compare Figure~\ref{fig:exp2_stripes}).
	On the other hand, when the randomly-structured initial condition \eqref{eq:ICrnd} is taken, the effect of
	anoikis evasion of the cancer cells becomes clear through the higher tumor heterogeneity with more fractured patterns and correspondingly lower PCC density in regions with stronger degraded ECM.
	%invasion and lower MCC density in regions of sparse ECM.
\end{experiment}

%----------------------------------------------
\section{Discussion}\label{sec:discussion}

In this note we reviewed existing models with multiple taxis upon classifying them in three categories, cf. Section \ref{sec:intro}. We then introduced in Section \ref{sec:model} a novel multiscale model with double taxis belonging to the last of those categories. It describes tumor invasion, accounting for two mutually exclusive cell phenotypes: moving/proliferating. The motile cells are supposed to perform nonlinear diffusion and haptotaxis and be repelled by the proliferating ones. The phenotypic switch rates depend on subcellular dynamics, more precisely on the amount of cell receptors occupied during cell-cell and cell-tissue interactions; this -together with the dependency of the tactic sensitivity functions on the (equilibria of) subcellular variables- confers the model its multiscale character. In Section \ref{sec3} we provided a global existence result concerning weak solutions to a simplified, pure macroscopic version of the model. Several issues remain open from the viewpoint of mathematical analysis:
\begin{itemize}
	\item Well-posedness of the full multiscale model \eqref{eq:model}, with diffusion and taxis coefficients which are more general than those considered in the analysis of Section \ref{sec3}. In particular, sensitivity functions like those in \eqref{sensitivities-first} or \eqref{eq:microSteady2} should be allowed. Moreover, the diffusion coefficient should be allowed to degenerate, e.g. upon choosing it in the form \eqref{diff-coeff-degen} instead of that given in \eqref{eq:model-a}. This type of degeneracy requires at least two solution components to become zero and seems thus to be milder than e.g., the double degeneracy handled in \cite{ZSH,ZSU} for pure haptotaxis models (with and without phenotypic heterogeneity of the same kind as here). The double taxis towards/away from immotile signals, however, can create supplementary challenges. 
	\item The boundedness of solutions remains to be proved even for the simplified model version considered in Section \ref{sec3} and is correspondingly more challenging for the more general versions. Same applies to the issue of asymptotic behavior of solutions. All scenarios simulated here involved solutions which remained bounded for all times, however another choice of source terms and motility coefficients and/or a different parameter regime might possibly lead to singular structure formation.
	\item Comparison (in terms of qualitative behavior analysis) between the model \eqref{eq:model} with haptotaxis and interspecies repellence and the corresponding formulation with haptotaxis and repellent pH-taxis, like that investigated numerically in Experiment \ref{exp:acidity}. Notice that considering \eqref{eq:modelacidity} instead of \eqref{eq:model} changes the model category from (\ref{M.iii}) into (\ref{M.ii}), according to the classification in Section \ref{sec:intro}.  The model \eqref{eq:modelacidity} features chemotaxis away from a signal (acidity) which is partially produced by the tactic species and partially by the non-diffusing one, whereby the growth of the latter is, in turn, influenced by that very signal, who is, moreover, degrading the other tactic cue. This leads to intricate couplings, even in a single-scale framework with linear diffusion. The question arises whether in spite of these, the haptotaxis-chemotaxis model in category (\ref{M.ii}) has better analytical tractability than its original counterpart in category (\ref{M.iii}).
	\item Analysis of pattern formation for different model versions: the numerical simulations showed the occurrence of various patterns, which seem to be influenced by several factors investigated in the simulation scenarios: dynamic rates of phenotypic switch, degenerating diffusion coefficients, form of the growth/degradation terms, initial space distribution of the underlying tissue, direct (interspecies)/indirect (acidity-mediated) repellence. 
\end{itemize}
Finally, we would also like to address some issues related to modeling:
\begin{itemize}
	\item The way we built our models in this paper was rather heuristic, mainly relying on the balance of fluxes (diffusion, haptotaxis, interspecies repellence) and source terms for proliferation/decay/phenotype transitions. It would be interesting to find a way allowing a more or less formal deduction of such or related models from dynamics on lower scales (single cells, populations of cells structured according to several variables) and more basic principles. Attempts in this direction have been made in the context of cancer invasion e.g. in \cite{Chauviere2007,ConteSurulescu,Corbin2018,MMCIITINF14,EHKS,EHS,EKS16,HPW12,Hunt2016,KPSZ,Nieto2016} or for chemotaxis in other contexts e.g. in \cite{chalub,HiPa08}, see also \cite{BeBeTaWi} and references therein, to name but a few. 
	\item We already mentioned the influence of the tissue structure on the tumor patterns. Here we considered for simple comparison purposes two types of initial conditions for the tissue evolution. For certain real-world problems concerning tumor invasion the underlying tissue distribution can be assessed, e.g. patient-specific brain tissue reconstruction from diffusion tensor imaging (a variant of MRI) data. Understanding and predicting the patterns of tumor growth and spread are of major importance, as they may facilitate diagnosis: to keep the example of brain tumors, patterns observed on histological slides are used to grade the tumor and thus to make a survival prognosis, see e.g. \cite{kleih}.
	\item In Experiment \ref{exp:acidity} we modeled the influence of acidity on tumor evolution, which was more indirect than the interspecies repellence. A more precise modeling would have been to consider subcellular dynamics for the interaction of the tumor cells with acidity. e.g. by characterizing the occupancy of certain transmembrane units (specific receptors, ion channels, etc.) and let the phenotypic switch rates depend on such variables. This idea was employed in \cite{ConteSurulescu,Corbin2018} for a more careful deduction of acidity-mediated brain tumor development and patterning leading to macroscopic PDEs with another type of diffusion, but with similar taxis terms.
\end{itemize}

%TODO: 
%\begin{itemize}
%	\item Review all figure captions (did it for all scenarios, modulo acidity-triggered behavior).
%	\item Spell check for either American or British English.
	%\item Comment on the structure of tissue being relevant for the tumor patterns in each of the scenarios. Tissue assessment available via DTI (e.g., for reconstruction of brain structure). Importance of patterns: may facilitate diagnosis (e.g., grading of brain tumors).
%	\item Comment on multiscale modeling of response to acidity.
%		\item {\cred Please ensure that all numerics references are cited in the text.}
%\end{itemize}

%A sentence about the multiple taxis problems being difficult for the analysis, therefore the
%preliminary study of parabolic-elliptic systems can provide insight into the ways of handling the
%difficulties (perhaps better placed where we said that we ignore such studies in order not to extend
%too much the review).
%Also test what would happen when the tissue is built up by (moving) cells (and not by itself).

%What would be the difference to a model with 'genuine' pH-taxis and haptotaxis (still in the
%go-or-grow setting)? That would mean indirect signal production... 

\section*{Acknowledgments}
N. Kolbe gratefully acknowledges the support by the International Research Fellowship of the
Japanese Society for the Promotion of Science. C. Surulescu was partially supported by the Federal Ministry of Education and Research BMBF, project \textit{GlioMaTh} 05M2016.

\appendix \label{sec:appendix}
\renewcommand{\theequation}{\Alph{section}.\arabic{equation}}
\setcounter{equation}{0}
%--------------------------------
\section{Construction of the randomly structured ECM.}\label{app:ECM}
\begin{figure}
	\centering
	\begin{tabular}{ccc}
		\includegraphics[width=0.20\linewidth]{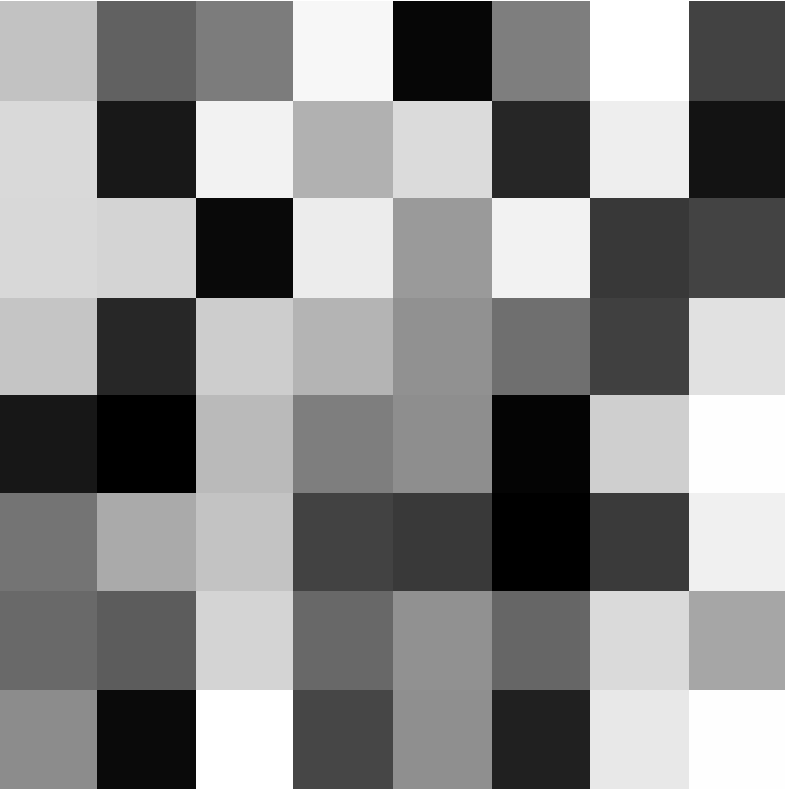}&
		\includegraphics[width=0.20\linewidth]{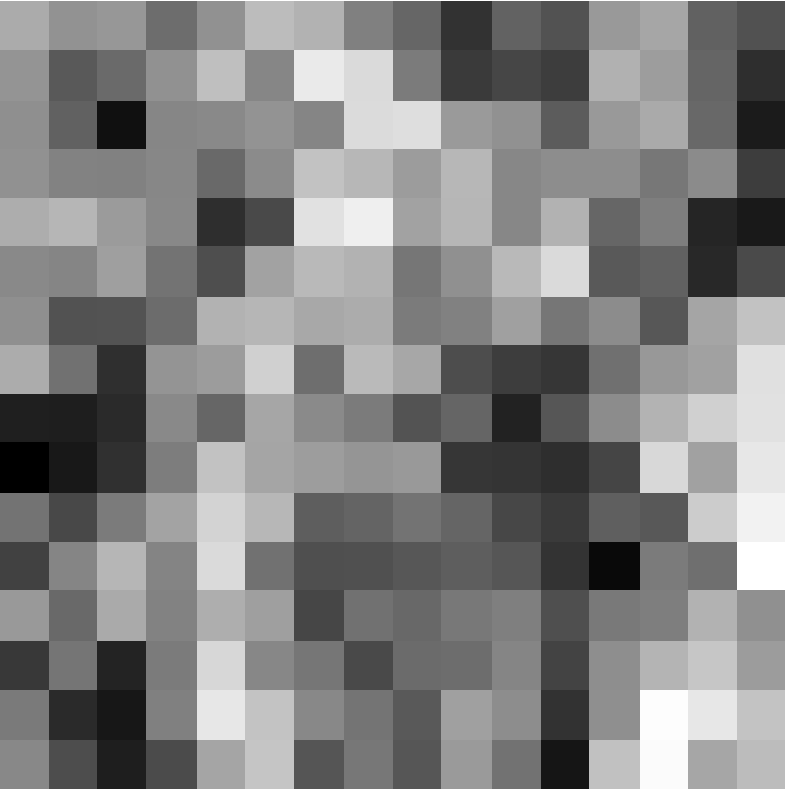}
		&	\includegraphics[width=0.20\linewidth]{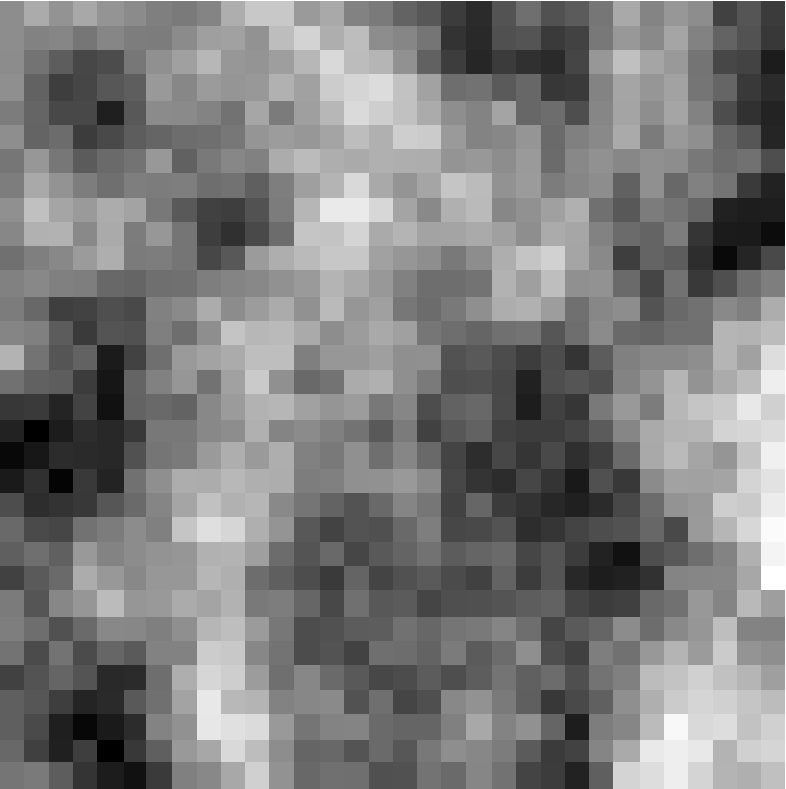}\\
		(a) $8\times8$&(b) $16\times16$&(c) $32\times32$\\[1ex]
		\includegraphics[width=0.20\linewidth]{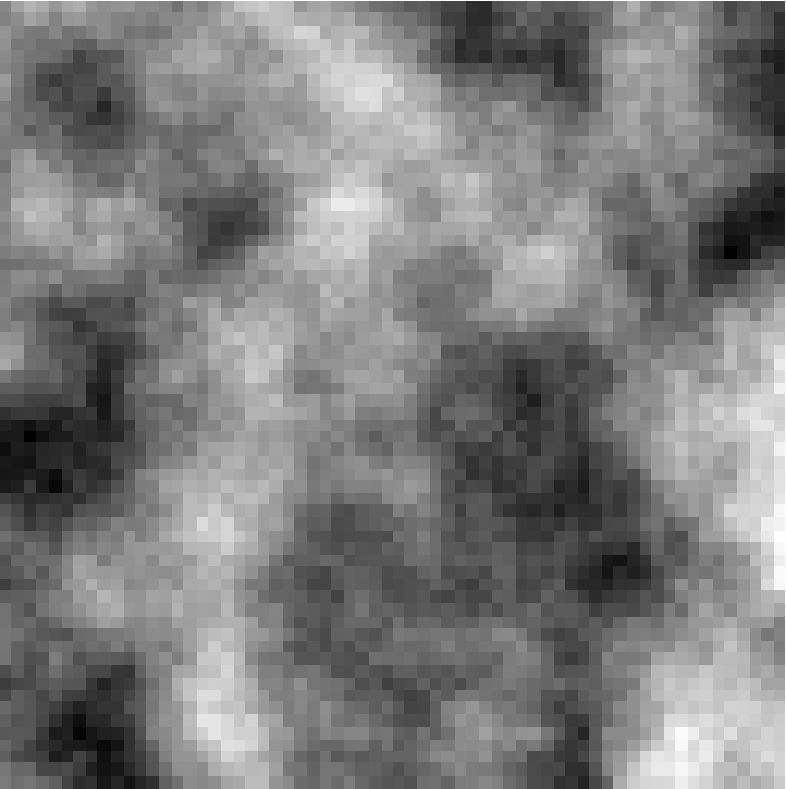}&
		\includegraphics[width=0.20\linewidth]{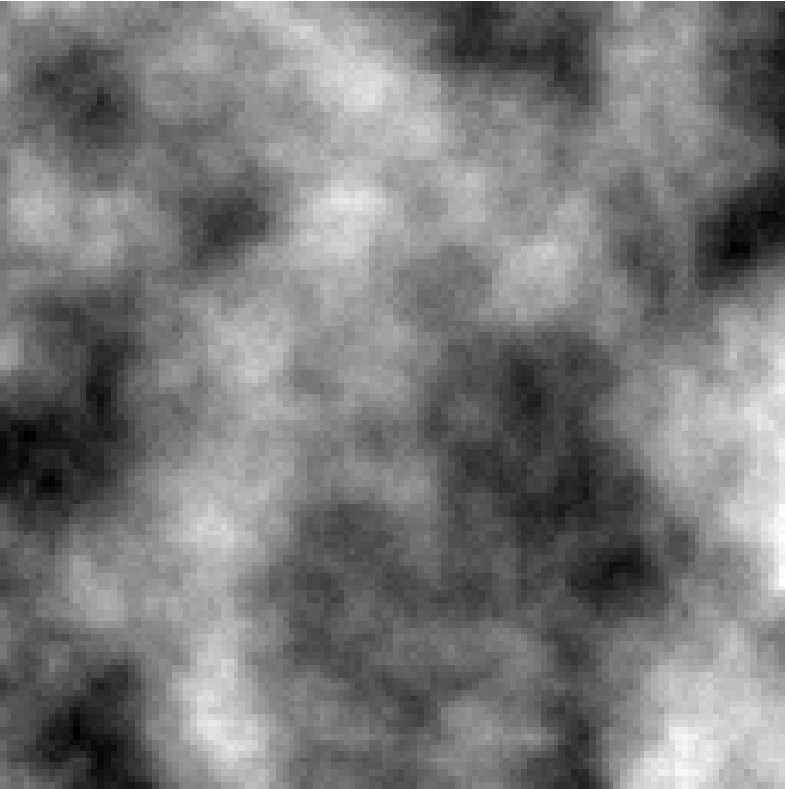}
		&	\includegraphics[width=0.20\linewidth]{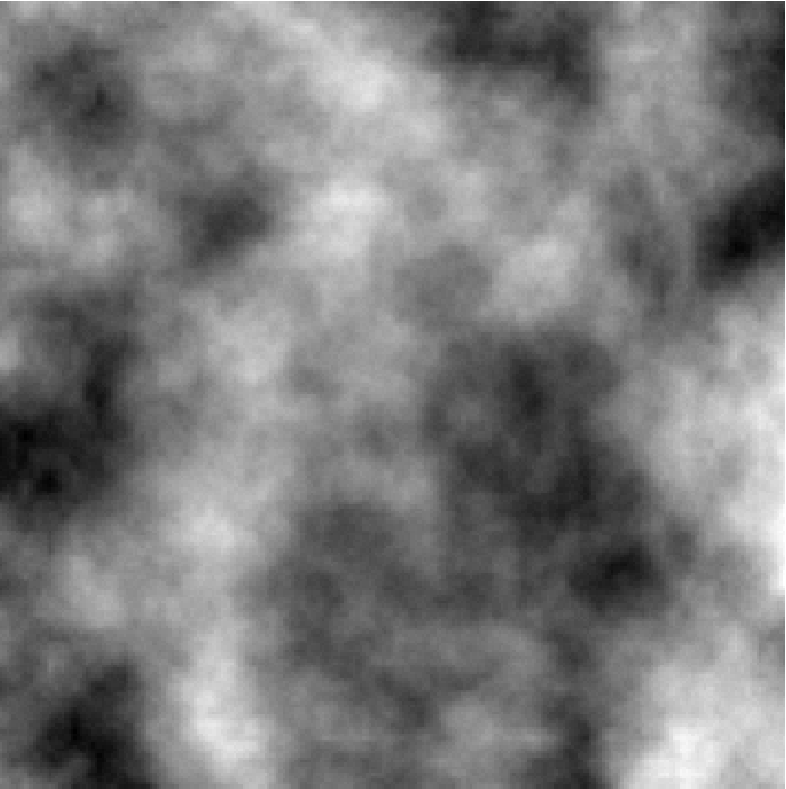}\\
		(d) $64\times64$&(e) $128\times128$&(f) $256\times256$\\
	\end{tabular}
	\caption{Construction of the randomly-structured ECM with a sequence of grid refinements steps. The
		first stage of this process is the construction of a random $8\times8$ grid (top left panel) with
		values normally distributed in $[0,1]$. This grid is progressively refined to the final (for this
		case) resolution of $256\times256$ (bottom right panel). At every refinement step the number of
		computational cells is doubled along each dimension and the new values are obtained by a) averaging
		the values of the neighboring  cells of the coarser grid, and b) adding some random and normally
		distributed noise. Periodic interpolations are employed at the ``boundary'' the discretization
		domain. It can be clearly seen that the coarse structure of the ECM that was randomly chosen in the
		$8\times8$ matrix is still visible in the refined $256\times 256$ grid.}
	\label{fig:ecmrefine}
\end{figure}

The randomly constructed matrix that we use in our numerical experiments is defined inductively. For
the sake of simplicity we describe here the one-dimensional case, over the domain $[0,1)$, and refer
to Figure~\ref{fig:ecmrefine} for a graphical representation of the same process in two-dimensional
domains. 

At first, a coarse approximation of the ECM is decided by setting the number of major ``hills'' and
``valleys'' in the density of the matrix. Should this number be $8$ (to coincide with
Figure\ref{fig:ecmrefine}), the first approximation to the ECM is set 
\[
\sum_{i=1}^{8}  c^{(8)}_i \mathcal X_{C^{(8)}_i}(x), \quad x\in [0,1),
\]
where $C^{(8)}_i$, with $|C^{(8)}_i|=\Delta x^{(8)}$ for $i=1,...,8$, represent the $8$
computational cells of the uniform discretization of $[0,1)$, and where the coefficients $c_i^{(8)}$
are uniformly distributed random numbers within $[0,1)$.	When we globally refine (by bisection), the
domain $[0,1)$ is discretized by $16$ equivalent computational cells $C^{(16)}_i$, $i=1,...,16$.
Accordingly, the ECM is approximated by the simple function
\[
\sum_{i=1}^{16} c^{(16)}_i \mathcal X_{C^{(16)}_i}(x), \quad x\in [0,1).
\]
The new coefficients $c^{(16)}_i$ interpolate---with the addition of some random noise---between the
previous values, i.e. 
\[
c^{(16)}_i = \left(1 + 0.002\left(r^{(16)}_i-0.5\right)\right)\frac{ c^{(8)}_{\floor{i/2}} +
	c^{(8)}_{\floor{i/2}+1} }{2},\quad i=1\ldots 16,
\]
where $\floor{\cdot}$ represents the Gauss floor function, and where $r^{(16)}_i$ are uniformly
distributed random numbers within $[0,1)$. The first and last coefficients, $c^{(16)}_1$ and
$c^{(16)}_{16}$, are computed periodically with respect to the $c^{(8)}_{\cdot.}$ values. The
rescaling factor $0.002$ is chosen so that the multiplicative randomness/noise is adjusted to
$0.1\%$ of the interpolated value. A similar refinement process is iterated until the desired
resolution of the ECM is reached. Then values of the ECM density are rescaled within the biological
range of a minimum and maximum ECM density.

%--------------------------------
\mysection{Numerical method}\label{app:numerics}

To numerically solve  the system  we first rewrite model \eqref{eq:model} for  $\mathbf x\in \Omega$
and $ \mathbf w(\mathbf x,t) := ( m(\mathbf x,t),$ $p(\mathbf x, t),$ $v(\mathbf x,t))$ %$a(\mathbf x,t) )$ 
in the form 
\begin{equation}\label{eq:w.model}
\mathbf w_t (\mathbf x,t)= A(\mathbf w(\mathbf x,t)) + D(\mathbf w(\mathbf x,t)) + R(\mathbf
w(\mathbf x,t)),
\end{equation}
where
\begin{subequations}
	\begin{align}
	A(\mathbf w) &= 
	\begin{pmatrix} 
	-\nabla \cdot \left(\chi_1(m,p,v)\,m\,\nabla v \right)+\nabla \cdot \left(\chi_2(m,p,v)\,m\,\nabla
	p \right),\ 0, \ 0
	\end{pmatrix}^T,
	\\
	R(\mathbf w) &= 
	\begin{pmatrix} 
	\lambda(y,\zeta)p-\gamma(y,\zeta)m\\
	\mu p(1-(m+p)-v)-\lambda(y,\zeta)p+\gamma(y,\zeta)m\\
	-\delta (m+p)v + \mu_v\, v (1-(m+p)-v)\\
	\end{pmatrix},
	\\
	D(\mathbf w) &= 
	\begin{pmatrix} 
	\nabla\cdot\left( D_c\frac{1+mp+mv+pv}{1+m(p+v)}\nabla m \right),\  0,\ 0
	\end{pmatrix}^T,\label{eq:diffusion}
	\end{align}
\end{subequations}
denote the advection, reaction, and diffusion operators respectively. The modified models considered
in Section~\ref{sec:numerics} can be rewritten analogously.

For the discretization of the spatial grid, we consider a uniform grid of diameter $h$, and denote
by $\mathbf w_h(t)$ the piecewise constant \textit{Finite Volume} (FV) approximation of the exact
solution $\mathbf w$, that satisfies the semi-discrete numerical scheme
\begin{equation}\label{eq:app.gen.scheme}
\partial_t \mathbf w_h = \mathcal A(\mathbf w_h) + \mathcal R(\mathbf w_h) + \mathcal D(\mathbf
w_h).
\end{equation}
The operators $\mathcal A_h$, $\mathcal R_h$, and $\mathcal D_h$ are \textit{discrete
	approximations} of the advection, reaction, and diffusion operators $A$, $R$, and $D$ in
\eqref{eq:w.model} respectively. For the diffusion terms we use central differences and
\textit{central upwind} numerical fluxes for the advection terms. At the interfaces of the
computational cells we use values of $\mathbf w_h$ reconstructed by the \textit{minimized-central}
(MC) limiter, \cite{van_leer_towards_1979}.

The semi-discrete scheme \eqref{eq:app.gen.scheme} is then solved with \textit{Implicit-Explicit
	Runge-Kutta} (IMEX-RK) numerical method, \cite{pareschi_implicitexplicit_2005}, that is based on the
\textit{time splitting} of \eqref{eq:app.gen.scheme}, in \textit{explicit} and \textit{implicit}
terms, in the form
\begin{equation}\label{eq:app.gen.IMEX}
\partial_t \mathbf w_h = \mathcal I(\mathbf w_h) + \mathcal E(\mathbf w_h).
\end{equation}
In the typical case, and also in the current paper, the advection terms are treated explicitly, the
diffusion terms implicitly, and the reaction terms explicitly. For the implicit part of the scheme
we employ a diagonally implicit RK method and an explicit RK for the explicit part
\begin{equation}\label{eq:app.IMEXRK}
\left\{\begin{array}{ll}
\displaystyle
\mathbf W_i^\ast = \mathbf w_h^n + \tau_n \sum_{j=1}^{i-2}\bar a_{i,j}\mathbf E_j + \tau_n \bar
a_{i,i-1}\mathbf E_{i-1},&\quad i=1,\dots ,s\\
\displaystyle
\mathbf W_i = \mathbf W_i^\ast + \tau_n \sum_{j=1}^{i-1} a_{i,j}\mathbf I_j + \tau_n a_{i,i}\mathbf
I_i,&\quad i=1,\dots, s\\
\displaystyle
\mathbf w_h^{n+1} = \mathbf w_h^n + \tau_n \sum _{i=1}^s\bar b_i \mathbf E_i + \tau_n
\sum_{i=1}^sb_i\mathbf I_i
\end{array}\right..
\end{equation}
Here $s$ represents the number of stages of the IMEX method, $\mathbf E_i=\mathcal E(\mathbf W_i)$,
$I_i=\mathcal I(\mathbf W_i)$, $i=1\dots s$, $\{\bar b,\, \bar A\}$, $\{b,\, A\}$ are respectively
the coefficients for the explicit and the implicit part of the scheme. We employ the L-stable and
stiffly accurate scheme ARK3(2)4L[2]SA from \cite{kennedy_additive_2003} with Butcher tableau given
in Table~\ref{F:tbl:IMEX}.

To solve the linear systems in \eqref{eq:app.IMEXRK} we use the \textit{iterative biconjugate
	gradient stabilized Krylov subspace} method \cite{van_der_vorst_bi-cgstab_1992}.

\renewcommand{\arraystretch}{1.5}
\begin{table}[t]
	\caption{Butcher tableau for the explicit (upper) and the implicit (lower) parts of the third
		order IMEX scheme ARK3(2)4L[2]SA we use in \eqref{eq:app.IMEXRK}, see also
		\cite{kennedy_additive_2003}.}
	\label{F:tbl:IMEX}
	\centering
	\scriptsize
	\begin{tabular}{c|cccc}
		$0$&&&&\\[0.5em]
		$\frac{1767732205903}{2027836641118}$&$\frac{1767732205903}{2027836641118}$&&&\\[0.5em]
		
		$\frac{3}{5}$&$\frac{5535828885825}{10492691773637}$&$\frac{788022342437}{10882634858940}$&&\\[0.5em]
		
		$1$&$\frac{6485989280629}{16251701735622}$&$-\frac{4246266847089}{9704473918619}$&$\frac{10755448449292}{10357097424841}$&\\[0.5em]
		\hline\\[-0.5em]
		&$\frac{1471266399579}{7840856788654}$ & $-\frac{4482444167858}{7529755066697}$ &
		$\frac{11266239266428}{11593286722821}$ & $\frac{1767732205903}{4055673282236}$ 
	\end{tabular}
	\\[0.5em]
	\begin{tabular}{c|cccc}
		$0$&0&&&\\[0.5em]
		
		$\frac{1767732205903}{2027836641118}$&$\frac{1767732205903}{4055673282236}$&$\frac{1767732205903}{4055673282236}$&&\\[0.5em]
		
		$\frac{3}{5}$&$\frac{2746238789719}{10658868560708}$&$-\frac{640167445237}{6845629431997}$&$\frac{1767732205903}{4055673282236}$&\\[0.5em]
		
		$1$&$\frac{1471266399579}{7840856788654}$&$-\frac{4482444167858}{7529755066697}$&$\frac{11266239266428}{11593286722821}$&$\frac{1767732205903}{4055673282236}$\\[0.5em]
		\hline\\[-0.5em]
		&$\frac{1471266399579}{7840856788654}$ & $-\frac{4482444167858}{7529755066697}$ &
		$\frac{11266239266428}{11593286722821}$ & $\frac{1767732205903}{4055673282236}$ 
	\end{tabular}
\end{table}
%\clearpage
\newcommand{\noopsort}[1]{}
\addcontentsline{toc}{section}{References}
%\nocite{*} %not cited works are not shown
\bibliographystyle{plain}
\bibliography{literature}

\begin{figure}[t] % exp3 - stripes
	
	\figurewidth=\linewidth
	\def \expSetup {TODO2varProlif_ECM_stripes}
	\begin{tikzpicture}
	\begin{groupplot}[
	/tikz/mark size=1.5pt,
	group style={
		group name=my plots,
		group size=4 by 4,
		horizontal sep=1cm,      % <-- default: 1cm
		vertical sep=0.7cm,        % <-- default: 1cm
	},
	xmin=-2,
	xmax=2,
	ymin=-2,
	ymax=2,
	ylabel shift = 2 em,
	xtick = \empty,
	ytick= \empty,
	ticklabel style = {font=\scriptsize},
	axis line style = thick,
	axis background/.style={fill=white},
	width=.25\figurewidth,
	height=.25\figurewidth, 
	]
	\nextgroupplot[ylabel = MCCs, title = {$t=0$}]
	\addplot [forget plot] graphics [xmin=-2, xmax=2, ymin=-2, ymax=2] {\expSetup1};
	\nextgroupplot[title = {$t=3.3$}]
	\addplot [forget plot] graphics [xmin=-2, xmax=2, ymin=-2, ymax=2] {\expSetup2};
	\nextgroupplot[title = {$t=6.7$}]
	\addplot [forget plot] graphics [xmin=-2, xmax=2, ymin=-2, ymax=2] {\expSetup3};
	\nextgroupplot[
	title = {$t=10$},
	point meta min=0,
	point meta max=0.0866,
	colormap name=MCCs,
	colorbar,
	colorbar style={at={(1.2,1)},anchor=north west}
	]
	\addplot [forget plot] graphics [xmin=-2, xmax=2, ymin=-2, ymax=2] {\expSetup4};
	\nextgroupplot[ylabel = PCCs]
	\addplot [forget plot] graphics [xmin=-2, xmax=2, ymin=-2, ymax=2] {\expSetup5};
	\nextgroupplot
	\addplot [forget plot] graphics [xmin=-2, xmax=2, ymin=-2, ymax=2] {\expSetup6};
	\nextgroupplot
	\addplot [forget plot] graphics [xmin=-2, xmax=2, ymin=-2, ymax=2] {\expSetup7};
	\nextgroupplot[
	point meta min=0,
	point meta max=1,
	colormap name=PCCs,
	colorbar,
	colorbar style={at={(1.2,1)},anchor=north west}
	]
	\addplot [forget plot] graphics [xmin=-2, xmax=2, ymin=-2, ymax=2] {\expSetup8};
	\nextgroupplot[ylabel = ECM]
	\addplot [forget plot] graphics [xmin=-2, xmax=2, ymin=-2, ymax=2] {\expSetup9};
	\nextgroupplot
	\addplot [forget plot] graphics [xmin=-2, xmax=2, ymin=-2, ymax=2] {\expSetup10};
	\nextgroupplot
	\addplot [forget plot] graphics [xmin=-2, xmax=2, ymin=-2, ymax=2] {\expSetup11};
	\nextgroupplot[
	point meta min=0,
	point meta max=1,
	colormap name=ECM,
	colorbar,
	colorbar style={at={(1.2,1)},anchor=north west}
	]
	\addplot [forget plot] graphics [xmin=-2, xmax=2, ymin=-2, ymax=2] {\expSetup12};
	\nextgroupplot[ylabel = pH level]
	\addplot [forget plot] graphics [xmin=-2, xmax=2, ymin=-2, ymax=2] {\expSetup13};
	\nextgroupplot
	\addplot [forget plot] graphics [xmin=-2, xmax=2, ymin=-2, ymax=2] {\expSetup14};
	\nextgroupplot
	\addplot [forget plot] graphics [xmin=-2, xmax=2, ymin=-2, ymax=2] {\expSetup15};
	\nextgroupplot[
	point meta min=6.4,
	point meta max=7.2,
	colormap name=AC,
	colorbar,
	colorbar style={at={(1.2,1)},anchor=north west}
	]
	\addplot [forget plot] graphics [xmin=-2, xmax=2, ymin=-2, ymax=2] {\expSetup16};
	\end{groupplot}
	\end{tikzpicture}%
	
	\caption{Simulation results of \textbf{Experiment~\ref{exp:acidity} -- acidity driven migration}
		with the  ECM-with-stripes initial conditions \eqref{eq:ICstrps}. In addition to the MCCs, PCCs, and ECM, we also visualize here the pH levels. When comparing with Experiment~\ref{exp:dynEMT} (dynamic phenotypic switch without acidity), Figure~\ref{fig:exp2_stripes}, the effect of the acidity can be seen in the more extensive spread of MCCs, due to chemorepellence by a self-diffusing signal, along with reduced proliferation due to hypoxia, and enhanced ECM degradation throughout the domain.}\label{fig:exp3_stripes}	
\end{figure}
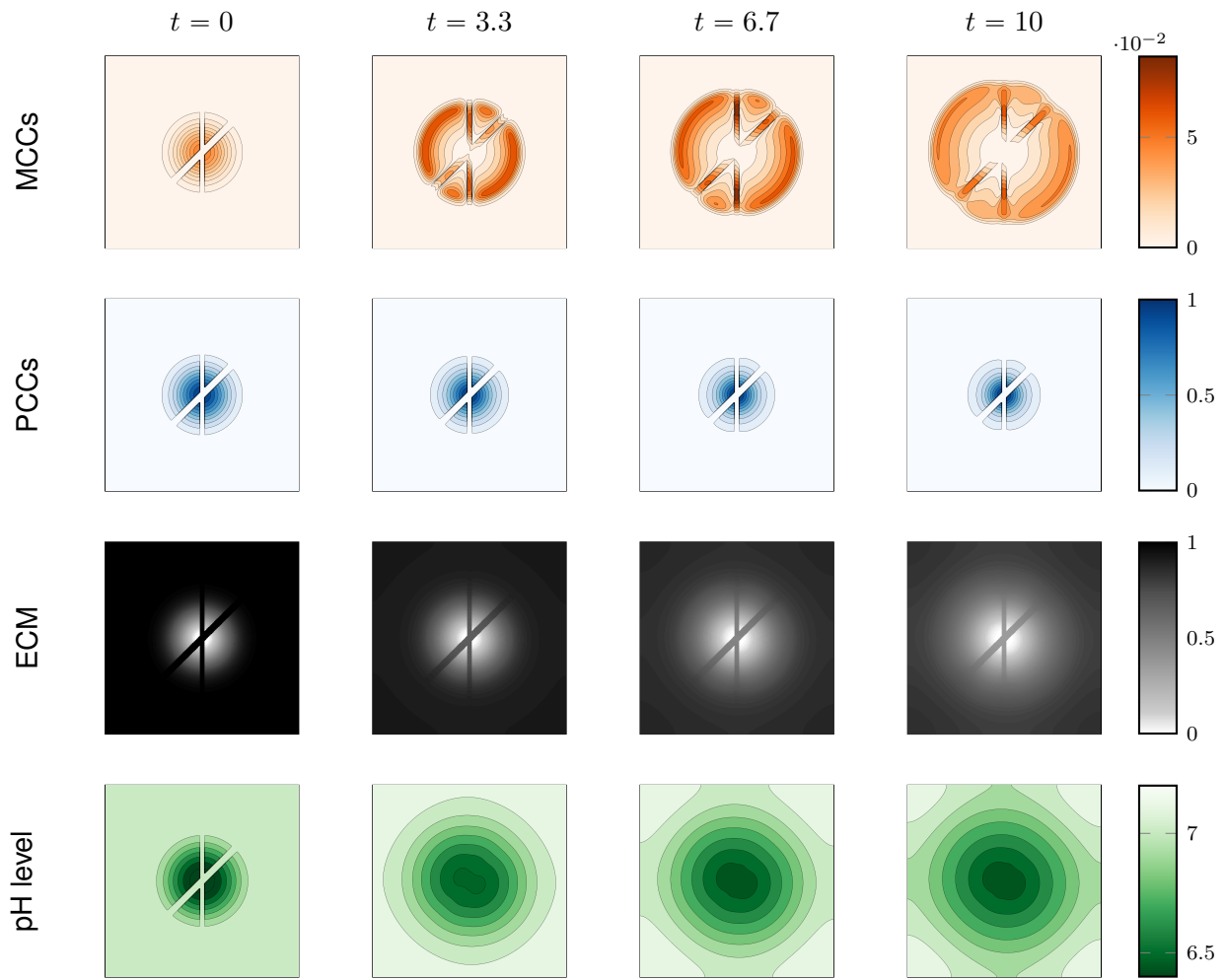

\begin{figure}[t] % exp3 - rnd
	\figurewidth=\linewidth
	\def \expSetup {TODO2varProlif_ECM_rnd}
	\begin{tikzpicture}
	\begin{groupplot}[
	/tikz/mark size=1.5pt,
	group style={
		group name=my plots,
		group size=4 by 4,
		horizontal sep=1cm,      % <-- default: 1cm
		vertical sep=0.7cm,        % <-- default: 1cm
	},
	xmin=-2,
	xmax=2,
	ymin=-2,
	ymax=2,
	ylabel shift = 2 em,
	xtick = \empty,
	ytick= \empty,
	ticklabel style = {font=\scriptsize},
	axis line style = thick,
	axis background/.style={fill=white},
	width=.25\figurewidth,
	height=.25\figurewidth, 
	]
	\nextgroupplot[ylabel = MCCs, title = {$t=0$}]
	\addplot [forget plot] graphics [xmin=-2, xmax=2, ymin=-2, ymax=2] {\expSetup1};
	\nextgroupplot[title = {$t=3.3$}]
	\addplot [forget plot] graphics [xmin=-2, xmax=2, ymin=-2, ymax=2] {\expSetup2};
	\nextgroupplot[title = {$t=6.7$}]
	\addplot [forget plot] graphics [xmin=-2, xmax=2, ymin=-2, ymax=2] {\expSetup3};
	\nextgroupplot[
	title = {$t=10$},
	point meta min=0,
	point meta max=0.372,
	colormap name=MCCs,
	colorbar,
	colorbar style={at={(1.2,1)},anchor=north west}
	]
	\addplot [forget plot] graphics [xmin=-2, xmax=2, ymin=-2, ymax=2] {\expSetup4};
	\nextgroupplot[ylabel = PCCs]
	\addplot [forget plot] graphics [xmin=-2, xmax=2, ymin=-2, ymax=2] {\expSetup5};
	\nextgroupplot
	\addplot [forget plot] graphics [xmin=-2, xmax=2, ymin=-2, ymax=2] {\expSetup6};
	\nextgroupplot
	\addplot [forget plot] graphics [xmin=-2, xmax=2, ymin=-2, ymax=2] {\expSetup7};
	\nextgroupplot[
	point meta min=0,
	point meta max=1,
	colormap name=PCCs,
	colorbar,
	colorbar style={at={(1.2,1)},anchor=north west}
	]
	\addplot [forget plot] graphics [xmin=-2, xmax=2, ymin=-2, ymax=2] {\expSetup8};
	\nextgroupplot[ylabel = ECM]
	\addplot [forget plot] graphics [xmin=-2, xmax=2, ymin=-2, ymax=2] {\expSetup9};
	\nextgroupplot
	\addplot [forget plot] graphics [xmin=-2, xmax=2, ymin=-2, ymax=2] {\expSetup10};
	\nextgroupplot
	\addplot [forget plot] graphics [xmin=-2, xmax=2, ymin=-2, ymax=2] {\expSetup11};
	\nextgroupplot[
	point meta min=0,
	point meta max=1,
	colormap name=ECM,
	colorbar,
	colorbar style={at={(1.2,1)},anchor=north west}
	]
	\addplot [forget plot] graphics [xmin=-2, xmax=2, ymin=-2, ymax=2] {\expSetup12};
	\nextgroupplot[ylabel = pH level]
	\addplot [forget plot] graphics [xmin=-2, xmax=2, ymin=-2, ymax=2] {\expSetup13};
	\nextgroupplot
	\addplot [forget plot] graphics [xmin=-2, xmax=2, ymin=-2, ymax=2] {\expSetup14};
	\nextgroupplot
	\addplot [forget plot] graphics [xmin=-2, xmax=2, ymin=-2, ymax=2] {\expSetup15};
	\nextgroupplot[
	point meta min=6.4,
	point meta max=7.2,
	colormap name=AC,
	colorbar,
	colorbar style={at={(1.2,1)},anchor=north west}
	]
	\addplot [forget plot] graphics [xmin=-2, xmax=2, ymin=-2, ymax=2] {\expSetup16};
	\end{groupplot}
	\end{tikzpicture}%
	
	\caption{Simulation results of \textbf{Experiment~\ref{exp:acidity} -- acidity driven migration}
	with the random-structured ECM initial conditions \eqref{eq:ICrnd}. Remarks analogous to those made in Figure~\ref{fig:exp3_stripes} apply here as well, when correspondingly comparing  with Experiment~\ref{exp:dynEMT} (dynamic phenotypic switch without acidity), Figure~\ref{fig:exp2_rnd}.}\label{fig:exp3_rnd}	
\end{figure}
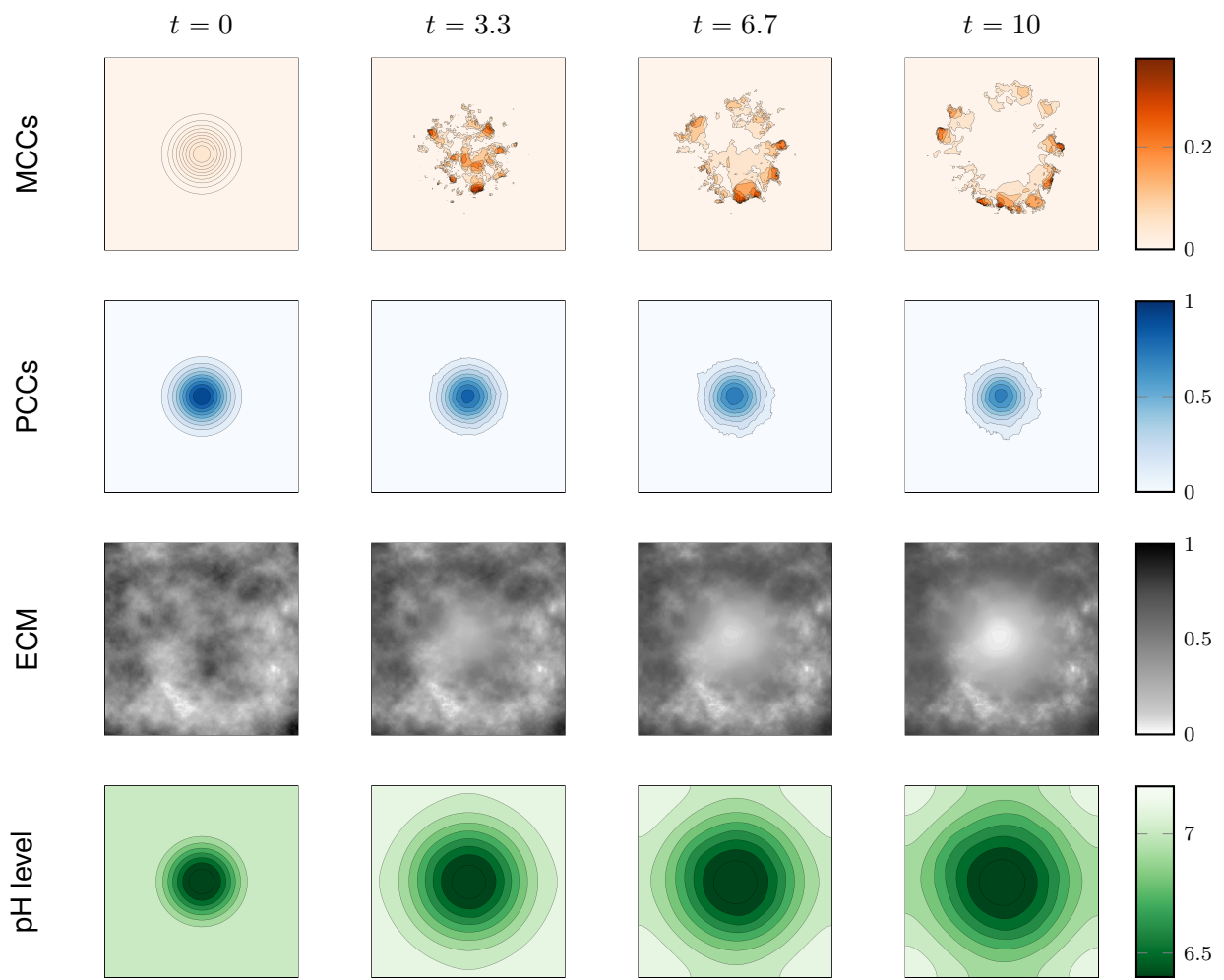

\begin{figure}[t] % epx4 - stripes
	\figurewidth=\linewidth
	\def \expSetup {TODO3_ECM_stripes}
	\begin{tikzpicture}
	\begin{groupplot}[
	/tikz/mark size=1.5pt,
	group style={
		group name=my plots,
		group size=4 by 3,
		horizontal sep=1cm,      % <-- default: 1cm
		vertical sep=0.7cm,        % <-- default: 1cm
	},
	xmin=-2,
	xmax=2,
	ymin=-2,
	ymax=2,
	ylabel shift = 2 em,
	xtick = \empty,
	ytick= \empty,
	ticklabel style = {font=\scriptsize},
	axis line style = thick,
	axis background/.style={fill=white},
	width=.25\figurewidth,
	height=.25\figurewidth, 
	]
	\nextgroupplot[ylabel = MCCs, title = {$t=0$}]
	\addplot [forget plot] graphics [xmin=-2, xmax=2, ymin=-2, ymax=2] {\expSetup1};
	\nextgroupplot[title = {$t=3.3$}]
	\addplot [forget plot] graphics [xmin=-2, xmax=2, ymin=-2, ymax=2] {\expSetup2};
	\nextgroupplot[title = {$t=6.7$}]
	\addplot [forget plot] graphics [xmin=-2, xmax=2, ymin=-2, ymax=2] {\expSetup3};
	\nextgroupplot[
	title = {$t=10$},
	point meta min=0,
	point meta max=0.0999,
	colormap name=MCCs,
	colorbar,
	colorbar style={at={(1.2,1)},anchor=north west}
	]
	\addplot [forget plot] graphics [xmin=-2, xmax=2, ymin=-2, ymax=2] {\expSetup4};
	\nextgroupplot[ylabel = PCCs]
	\addplot [forget plot] graphics [xmin=-2, xmax=2, ymin=-2, ymax=2] {\expSetup5};
	\nextgroupplot
	\addplot [forget plot] graphics [xmin=-2, xmax=2, ymin=-2, ymax=2] {\expSetup6};
	\nextgroupplot
	\addplot [forget plot] graphics [xmin=-2, xmax=2, ymin=-2, ymax=2] {\expSetup7};
	\nextgroupplot[
	point meta min=0,
	point meta max=1,
	colormap name=PCCs,
	colorbar,
	colorbar style={at={(1.2,1)},anchor=north west}
	]
	\addplot [forget plot] graphics [xmin=-2, xmax=2, ymin=-2, ymax=2] {\expSetup8};
	\nextgroupplot[ylabel = ECM]
	\addplot [forget plot] graphics [xmin=-2, xmax=2, ymin=-2, ymax=2] {\expSetup9};
	\nextgroupplot
	\addplot [forget plot] graphics [xmin=-2, xmax=2, ymin=-2, ymax=2] {\expSetup10};
	\nextgroupplot
	\addplot [forget plot] graphics [xmin=-2, xmax=2, ymin=-2, ymax=2] {\expSetup11};
	\nextgroupplot[
	point meta min=0,
	point meta max=1,
	colormap name=ECM,
	colorbar,
	colorbar style={at={(1.2,1)},anchor=north west}
	]
	\addplot [forget plot] graphics [xmin=-2, xmax=2, ymin=-2, ymax=2] {\expSetup12};
	\end{groupplot}
	\end{tikzpicture}%
	
	\caption{Simulation results of \textbf{Experiment~\ref{exp:degenerate} --- degenerate diffusion},
		with the ECM-with-stripes initial conditions \eqref{eq:ICrnd}. Compared to the non-degenerate
		diffusion in Experiment~\ref{exp:dynEMT}, shown in Figure~\ref{fig:exp2_stripes}, we note that there is a similar extent of tumor spread, however with MCCs forming very localized, relatively large aggregates (see also e.g., the closeup in Figure~\ref{fig:exp4_stripes_closeup}), while the PCC density remains almost the same.
		}\label{fig:exp4_stripes}	
\end{figure}
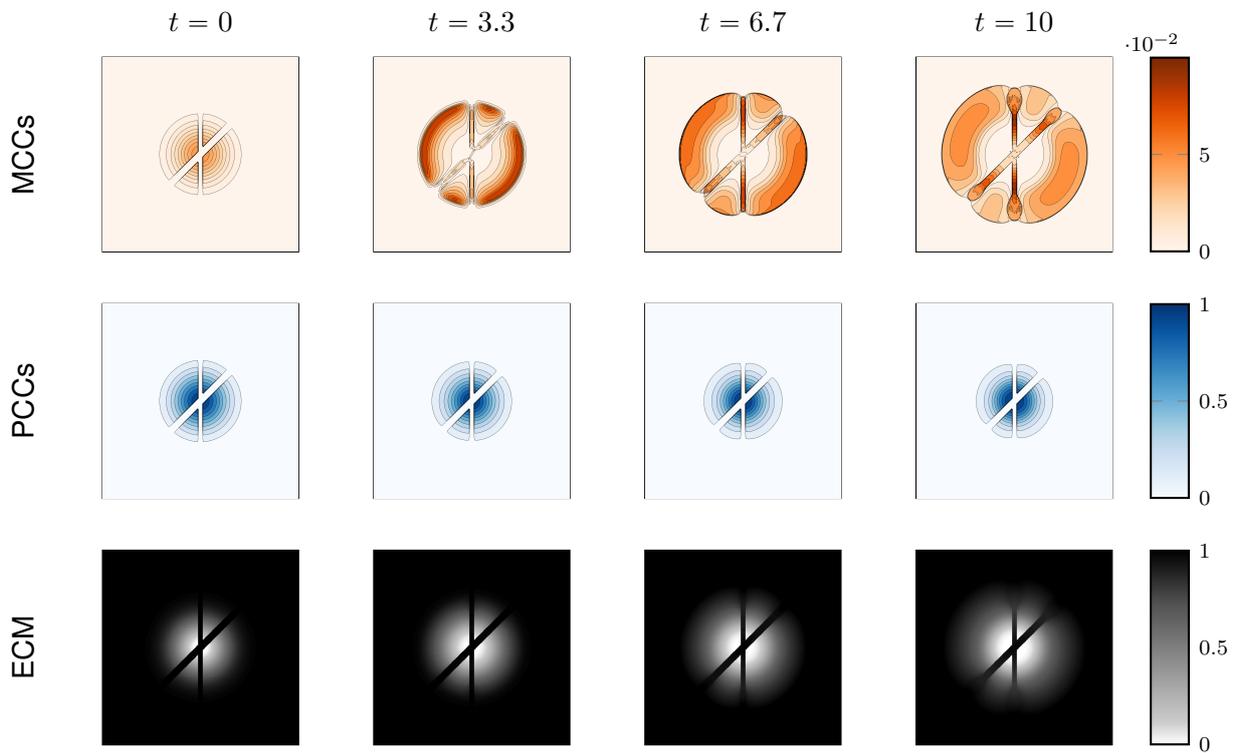

\begin{figure}[t] % epx4 - stripes - closeup
    \figurewidth=1.2\linewidth
    \def \expSetup {TODO3_ECM_stripes}
    \centering
    \begin{tikzpicture}
    \begin{groupplot}[
    /tikz/mark size=1.5pt,
    group style={
        group name=my plots,
        group size=3 by 1,
        horizontal sep=2cm,      % <-- default: 1cm
        vertical sep=0.7cm,        % <-- default: 1cm
    },
    xmin=-0.5,
    xmax=0.5,
    ymin=-0.5,
    ymax=0.5,
    ylabel shift = 1 em,
    xtick = \empty,
    ytick= \empty,
    ticklabel style = {font=\scriptsize},
    axis line style = thick,
    axis background/.style={fill=white},
    width=.25\figurewidth,
    height=.25\figurewidth, 
    ]
    \nextgroupplot[ylabel = MCCs ]
    \addplot [forget plot] graphics [xmin=-2, xmax=2, ymin=-2, ymax=2] {\expSetup4};
    \nextgroupplot[ylabel = PCCs]
    \addplot [forget plot] graphics [xmin=-2, xmax=2, ymin=-2, ymax=2] {\expSetup8};
    \nextgroupplot[ylabel = ECM ]
    \addplot [forget plot] graphics [xmin=-2, xmax=2, ymin=-2, ymax=2] {\expSetup12};
    \end{groupplot}
    \end{tikzpicture}%
    
    \caption{Closeup of the densities at $t=10$ in \textbf{Experiment~\ref{exp:degenerate} --- degenerate diffusion},
        with the ECM-with-stripes initial conditions \eqref{eq:ICrnd}. Compare to Figure~\ref{fig:exp4_stripes}. MCCs form localized aggregates while the PCC density remains similar as in Experiment~\ref{exp:dynEMT} shown in Figure~\ref{fig:exp2_stripes}.
    }\label{fig:exp4_stripes_closeup}	
\end{figure}
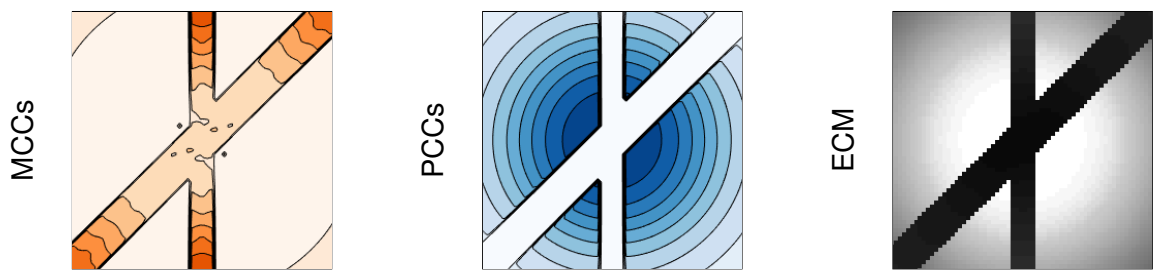

\begin{figure}[t] % exp4 - rnd
	\figurewidth=\linewidth
	\def \expSetup {TODO3_ECM_rnd}
	\begin{tikzpicture}
	\begin{groupplot}[
	/tikz/mark size=1.5pt,
	group style={
		group name=my plots,
		group size=4 by 3,
		horizontal sep=1cm,      % <-- default: 1cm
		vertical sep=0.7cm,        % <-- default: 1cm
	},
	xmin=-2,
	xmax=2,
	ymin=-2,
	ymax=2,
	ylabel shift = 2 em,
	xtick = \empty,
	ytick= \empty,
	ticklabel style = {font=\scriptsize},
	axis line style = thick,
	axis background/.style={fill=white},
	width=.25\figurewidth,
	height=.25\figurewidth, 
	]
	\nextgroupplot[ylabel = MCCs, title = {$t=0$}]
	\addplot [forget plot] graphics [xmin=-2, xmax=2, ymin=-2, ymax=2] {\expSetup1};
	\nextgroupplot[title = {$t=3.3$}]
	\addplot [forget plot] graphics [xmin=-2, xmax=2, ymin=-2, ymax=2] {\expSetup2};
	\nextgroupplot[title = {$t=6.7$}]
	\addplot [forget plot] graphics [xmin=-2, xmax=2, ymin=-2, ymax=2] {\expSetup3};
	\nextgroupplot[
	title = {$t=10$},
	point meta min=0,
	point meta max=0.34,
	colormap name=MCCs,
	colorbar,
	colorbar style={at={(1.2,1)},anchor=north west}
	]
	\addplot [forget plot] graphics [xmin=-2, xmax=2, ymin=-2, ymax=2] {\expSetup4};
	\nextgroupplot[ylabel = PCCs]
	\addplot [forget plot] graphics [xmin=-2, xmax=2, ymin=-2, ymax=2] {\expSetup5};
	\nextgroupplot
	\addplot [forget plot] graphics [xmin=-2, xmax=2, ymin=-2, ymax=2] {\expSetup6};
	\nextgroupplot
	\addplot [forget plot] graphics [xmin=-2, xmax=2, ymin=-2, ymax=2] {\expSetup7};
	\nextgroupplot[
	point meta min=0,
	point meta max=1,
	colormap name=PCCs,
	colorbar,
	colorbar style={at={(1.2,1)},anchor=north west}
	]
	\addplot [forget plot] graphics [xmin=-2, xmax=2, ymin=-2, ymax=2] {\expSetup8};
	\nextgroupplot[ylabel = ECM]
	\addplot [forget plot] graphics [xmin=-2, xmax=2, ymin=-2, ymax=2] {\expSetup9};
	\nextgroupplot
	\addplot [forget plot] graphics [xmin=-2, xmax=2, ymin=-2, ymax=2] {\expSetup10};
	\nextgroupplot
	\addplot [forget plot] graphics [xmin=-2, xmax=2, ymin=-2, ymax=2] {\expSetup11};
	\nextgroupplot[
	point meta min=0,
	point meta max=1,
	colormap name=ECM,
	colorbar,
	colorbar style={at={(1.2,1)},anchor=north west}
	]
	\addplot [forget plot] graphics [xmin=-2, xmax=2, ymin=-2, ymax=2] {\expSetup12};
	\end{groupplot}
	\end{tikzpicture}%
	
	\caption{Simulation results of \textbf{Experiment~\ref{exp:degenerate} --- degenerate diffusion}
		with the randomly structured ECM \eqref{eq:ICrnd}. When comparing with the non-degenerate diffusion
		in Experiment~\ref{exp:dynEMT}, shown in Figure~\ref{fig:exp2_rnd}, we note that the degenerate case leads to higher, more localized MCC densities, mainly near the invasion front.}\label{fig:exp4_rnd}	
\end{figure}
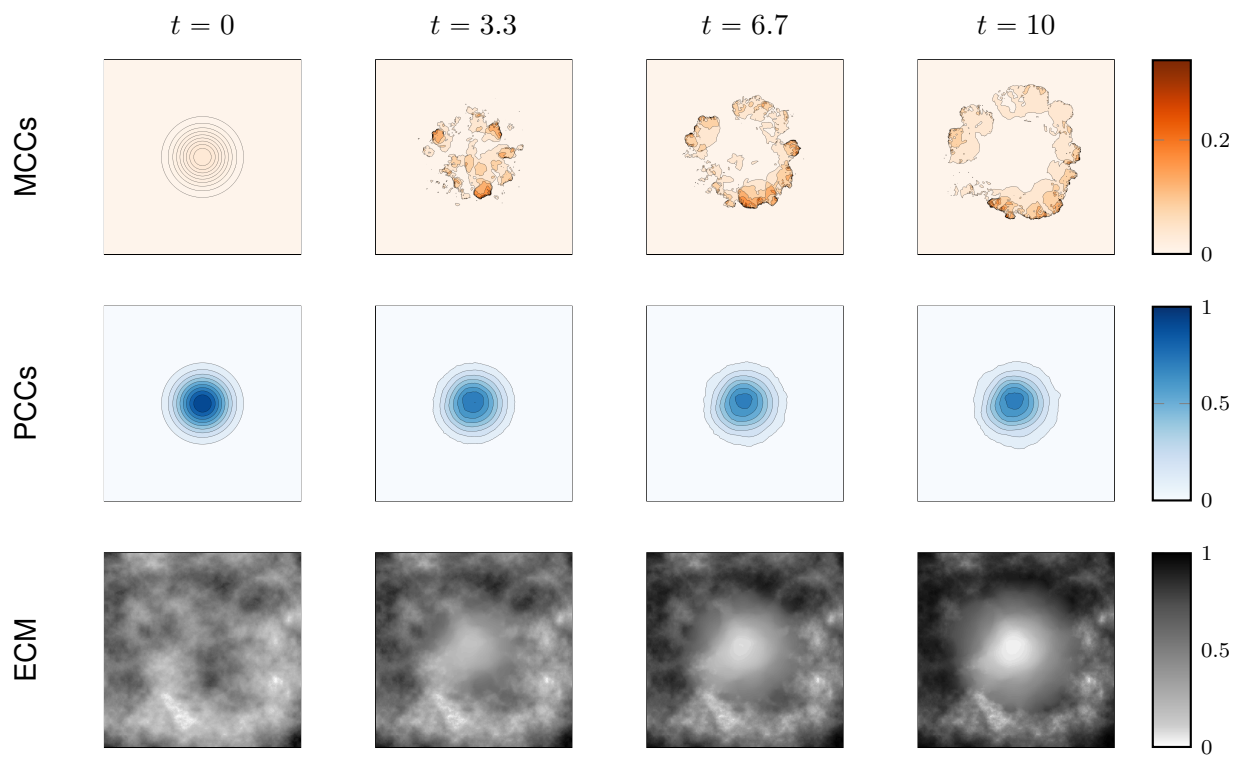

\begin{figure}[t] % exp6 -stripes
	\figurewidth=\linewidth
	\def \expSetup {TODO5_ECM_stripes}
	\begin{tikzpicture}
	\begin{groupplot}[
	/tikz/mark size=1.5pt,
	group style={
		group name=my plots,
		group size=4 by 3,
		horizontal sep=1cm,      % <-- default: 1cm
		vertical sep=0.7cm,        % <-- default: 1cm
	},
	xmin=-2,
	xmax=2,
	ymin=-2,
	ymax=2,
	ylabel shift = 2 em,
	xtick = \empty,
	ytick= \empty,
	ticklabel style = {font=\scriptsize},
	axis line style = thick,
	axis background/.style={fill=white},
	width=.25\figurewidth,
	height=.25\figurewidth, 
	]
	\nextgroupplot[ylabel = MCCs, title = {$t=0$}]
	\addplot [forget plot] graphics [xmin=-2, xmax=2, ymin=-2, ymax=2] {\expSetup1};
	\nextgroupplot[title = {$t=3.3$}]
	\addplot [forget plot] graphics [xmin=-2, xmax=2, ymin=-2, ymax=2] {\expSetup2};
	\nextgroupplot[title = {$t=6.7$}]
	\addplot [forget plot] graphics [xmin=-2, xmax=2, ymin=-2, ymax=2] {\expSetup3};
	\nextgroupplot[
	title = {$t=10$},
	point meta min=0,
	point meta max=0.09,
	colormap name=MCCs,
	colorbar,
	colorbar style={at={(1.2,1)},anchor=north west}
	]
	\addplot [forget plot] graphics [xmin=-2, xmax=2, ymin=-2, ymax=2] {\expSetup4};
	\nextgroupplot[ylabel = PCCs]
	\addplot [forget plot] graphics [xmin=-2, xmax=2, ymin=-2, ymax=2] {\expSetup5};
	\nextgroupplot
	\addplot [forget plot] graphics [xmin=-2, xmax=2, ymin=-2, ymax=2] {\expSetup6};
	\nextgroupplot
	\addplot [forget plot] graphics [xmin=-2, xmax=2, ymin=-2, ymax=2] {\expSetup7};
	\nextgroupplot[
	point meta min=0,
	point meta max=1,
	colormap name=PCCs,
	colorbar,
	colorbar style={at={(1.2,1)},anchor=north west}
	]
	\addplot [forget plot] graphics [xmin=-2, xmax=2, ymin=-2, ymax=2] {\expSetup8};
	\nextgroupplot[ylabel = ECM]
	\addplot [forget plot] graphics [xmin=-2, xmax=2, ymin=-2, ymax=2] {\expSetup9};
	\nextgroupplot
	\addplot [forget plot] graphics [xmin=-2, xmax=2, ymin=-2, ymax=2] {\expSetup10};
	\nextgroupplot
	\addplot [forget plot] graphics [xmin=-2, xmax=2, ymin=-2, ymax=2] {\expSetup11};
	\nextgroupplot[
	point meta min=0,
	point meta max=1,
	colormap name=ECM,
	colorbar,
	colorbar style={at={(1.2,1)},anchor=north west}
	]
	\addplot [forget plot] graphics [xmin=-2, xmax=2, ymin=-2, ymax=2] {\expSetup12};
	\end{groupplot}
	\end{tikzpicture}%
	
	\caption{Simulation results of \textbf{Experiment~\ref{exp:anoikis} --- anoikis effect} on an
	ECM with initial condition \eqref{eq:ICstrps}. Compared to Experiment~\ref{exp:dynEMT},
	shown in Figure~\ref{fig:exp2_stripes}, the results are almost identical; no particular anoikis effect is visible.}\label{fig:exp6_stripes}	
\end{figure}
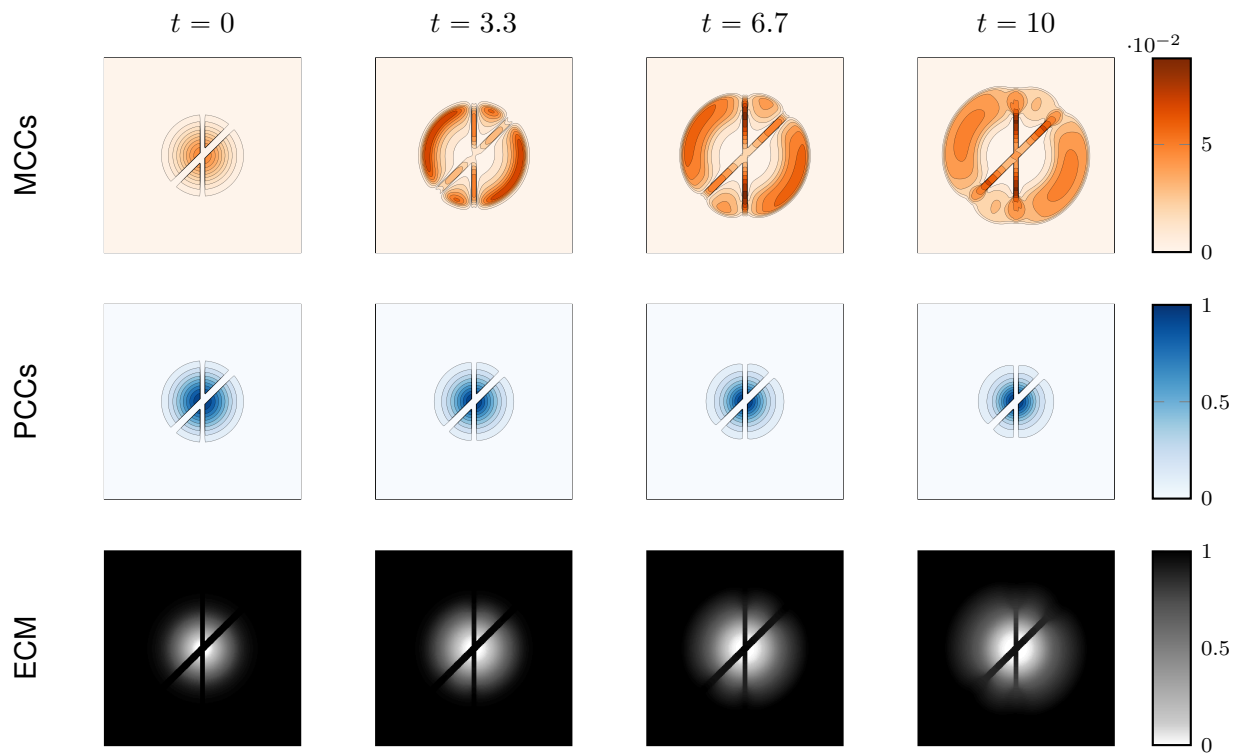

\section{Supplementary material}\label{app:supplement}

This part of the Appendix refers to video files of the simulation scenarios in Experiments \ref{exp:dynEMT}, \ref{exp:acidity}, and \ref{exp:remod} performed using the randomly-structured ECM initial conditions \eqref{eq:ICrnd}. The videos are available at \url{https://github.com/nklb/multi-taxis-supplementary} and can be downloaded using the links below.
 We refer to the corresponding sections of the manuscript for thorough discussions of the experiments.

\begin{tabular}{l}
	\href{https://github.com/nklb/multi-taxis-supplementary/blob/master/Experiment-2---Dynamic-phenotypic-switch-rates.mp4?raw=true}{\texttt{Experiment-2---Dynamic-phenotypic-switch-rates.mp4}}\\
	\href{https://github.com/nklb/multi-taxis-supplementary/blob/master/Experiment-3---Acidity-driven-migration.mp4?raw=true)}{\texttt{Experiment-3---Acidity-driven-migration.mp4}}\\
	\href{https://github.com/nklb/multi-taxis-supplementary/blob/master/Experiment-5---ECM-remodeling-by-the-cancer-cells.mp4?raw=true}{\texttt{Experiment-5---ECM-remodeling-by-the-cancer-cells.mp4}}
\end{tabular}\\
%\comm{If the journal supports the upload of supplementary material, links are to be replaced}
\end{document}